\RequirePackage{fix-cm}
\documentclass[twocolumn]{svjour3}          
\smartqed  
\usepackage{graphicx}
\RequirePackage{amsmath,amsfonts,amssymb}
\RequirePackage[authoryear]{natbib}
\RequirePackage[colorlinks,citecolor=blue,urlcolor=blue]{hyperref}
\RequirePackage{graphicx}

\usepackage[utf8]{inputenc}
\usepackage{dsfont}
\usepackage{stmaryrd}

\usepackage[caption=false,font=footnotesize,labelfont=sf,textfont=sf]{subfig}
\usepackage{booktabs}
\usepackage{multirow}
\usepackage{array}
\usepackage{xcolor}
\usepackage[shortlabels]{enumitem}
\usepackage[export]{adjustbox}

\usepackage{rotating}

\usepackage{orcidlink}

\usepackage{newtxtext,newtxmath}

\newcommand{\norm}[1]{\left\lVert#1\right\rVert}
\newcommand{\parent}[1]{\left(#1\right)}

\newcommand{\dotprod}[1]{\left<#1\right>}
\renewcommand{\brack}[1]{\left[#1\right]}
\renewcommand{\brace}[1]{\left\{#1\right\}}
\renewcommand{\det}[1]{\left|#1\right|}

\newcommand{\E}[1]{\mathbb{E}\brack{#1}}
\newcommand{\argmin}[1]{\underset{#1}{\text{argmin}} \,}

\renewcommand{\P}[1]{\mathbb{P}\parent{#1}}

\newcommand{\N}[1]{\mathcal{N}\parent{#1}}

\newcommand{\R}{\mathbb{R}}

\journalname{SN Computer Science}
\begin{document}

\title{Mixture of Conditional Gaussian Graphical Models for Unlabelled Heterogeneous Populations in the Presence of Co-factors}

\titlerunning{Mixture of Conditional Gaussian Graphical Models} 

\author{Thomas Lartigue$^{1, 2}$ \orcidlink{0000-0002-7820-0032}
        \and Stanley Durrleman$^1$ \orcidlink{0000-0002-9450-6920}
        \and St\'ephanie Allassonni\`ere$^3$ \orcidlink{0000-0002-5692-4945}
}

\authorrunning{Thomas Lartigue et al.} 

\institute{Thomas Lartigue$^*$ \at
            thomas.lartigue@inria.fr       \\
             \emph{Present e-mail:} thomas.lartigue@dzne.de  
           \and
          Stanley Durrleman \at
          stanley.durrleman@inria.fr
              \and
        St\'ephanie Allassonni\`ere \at
        stephanie.allassonniere@u-paris.fr
        \and
        $^*$ Corresponding author\\
        \\
        $^1$ Aramis project-team, INRIA, ICM, Paris, France\\
        \\
        $^2$ CMAP, CNRS, \'Ecole polytechnique, I.P. Paris, Palaiseau, France\\
        \\
        $^3$ Centre de Recherche des Cordeliers, Universit\'e de Paris, Inserm, HEKA Project team, INRIA Paris, Sorbonne Universit\'e, F-75006, Paris, France
}

\date{Received: date / Accepted: date}

\maketitle
\begin{abstract}
Conditional correlation networks, within Gaussian Graphical Models (GGM), are widely used to describe the direct interactions between the components of a random vector. In the case of an unlabelled Heterogeneous population, Expectation Maximisation (EM) algorithms for Mixtures of GGM have been proposed to estimate both each sub-population's graph and the class labels. However, we argue that, with most real data, class affiliation cannot be described with a Mixture of Gaussian, which mostly groups data points according to their geometrical proximity. In particular, there often exists external co-features whose values affect the features' average value, scattering across the feature space data points belonging to the same sub-population. Additionally, if the co-features' effect on the features is Heterogeneous, then the estimation of this effect cannot be separated from the sub-population identification. In this article, we propose a Mixture of Conditional GGM (CGGM) that subtracts the heterogeneous effects of the co-features to regroup the data points into sub-population corresponding clusters. We develop a penalised EM algorithm to estimate graph-sparse model parameters. We demonstrate on synthetic and real data how this method fulfils its goal and succeeds in identifying the sub-populations where the Mixtures of GGM are disrupted by the effect of the co-features.

\end{abstract}

\section{Introduction} \label{sect:intro}
The conditional correlation networks are a popular tool to describe the co-variations between the components of a random vector. Within the Gaussian Graphical Model (GGM) framework, introduced in \cite{dempster1972covariance}, the random vector of interest is modelled as a gaussian vector $\mathcal{N}(\mu, \Sigma)$, and the conditional correlation networks can be recovered from the sparsity of the precision matrix $\Lambda:=\Sigma^{-1}$. In this article, we consider the case of an unlabelled heterogeneous population, in which different sub-populations (or "classes") are described by different networks. Additionally, we take into account the presence of observed co-features (discrete and/or continuous) that have a heterogeneous (class-dependent) impact on the values of the features. The absence of known class labels turns the analysis of the population into an unsupervised problem. As a result, any inference method will have to tackle the problem of cluster discovery in addition to the parameter estimation. The former is a crucial task, especially since the relevance of the estimated parameters is entirely dependent on the clusters identified. The co-features, if their effects on the features are consequent, can greatly disrupt the clustering. Indeed, any unsupervised method will then be more likely to identify clusters correlated with the values of the co-features than with the hidden sub-population labels. This occurs frequently when analysing biological or medical features. To provide a simple illustration, if one runs an unsupervised method on an unlabelled population containing both healthy and obese patients, using the body fat percentage as a feature, then the unearthed clusters are very likely to be more correlated with the gender of the patients (a co-feature) rather than with the actual diagnostic (the hidden variable). Additionally, the fact that the effect of the gender on the average body fat is also dependent on the diagnostic (class-dependent effect) makes the situation even more complex.\\
The seminal work of \cite{dempster1972covariance}, and its later hierarchical extension introduced in \cite{honorio2010multi} and \cite{varoquaux2010brain}, generated tremendous interest in GGM and hierarchical GGM analyses respectively. See \cite{lartigue2020mixture} for an extensive literature review on these topics. Recently, unsupervised hierarchical GGM have received increasing attention, with works such as \cite{gao2016estimation} and \cite{hao2017simultaneous} adapting the popular supervised joint Hierarchical GGM methods of \cite{mohan2014node} and \cite{danaher2014joint} to the unsupervised case. When the labels are known in advance, these joint Hierarchical GGM are useful models to estimate several sparse conditional correlation matrices and are modular enough to allow for the recovery of many different forms of common structure between classes. However, we argue that they are not designed for efficient cluster identification in the unsupervised scenario, and will very likely miss the hidden variable and find clusters correlated to the most influential co-features instead. Which in turn will result in the estimation of irrelevant parameters. Even when there are no pre-existing hidden variables to recover, and the unsupervised method is run "blindly", it is uninteresting to recover clusters describing the values to already known co-features. Instead, one would rather provide beforehand the unsupervised method with the information of the co-features' values and encourage it to recover new information from the data.\\
In order to take into account the effect of co-features on features, \cite{yin2011sparse} and \cite{wytock2013sparse} introduced the Conditional Gaussian Graphical Models (CGGM). Within this model, the average effect of the co-features is subtracted from the features using a "linear transition" matrix $\Theta$, which leaves the features with only their residual behaviour. Both \cite{yin2011sparse} and \cite{sohn2012joint} worked with homogeneous populations, but the Hierarchical form of the CGGM was introduced by \cite{chun2013joint} to study labelled heterogeneous populations, with heterogeneous effects of the co-features on the features. Recent works such as \cite{huang2018joint} and \cite{ou2019differential} have adapted the state of the art supervised joint Hierarchical GGM methods to the CGGM, in order to also recover elements of common structure across all classes. However, to the best of our knowledge, there has been no effort to make use of the CGGM in the unsupervised case to identify groups with distinct sparse correlation structures with common elements.\\
In this article, we introduce a Mixture of Conditional GGM that models the class-dependent effect of the co-features on the features with different transition parameters $(\Theta_k)_{k=1}^K$. We propose an Expectation-Maximisation (EM) procedure to estimate this model without prior knowledge of the class labels. This EM algorithm can be regularised with all the structure-inducing penalties introduced for the supervised joint Hierarchical CGGM. Hence, the recovered sparse precision matrices $(\Lambda_k)_{k=1}^K$ and corresponding conditional correlation graphs can present any of the desired form of common structure. Moreover, with an additional penalty term, we can also enforce sparse and common structure within the transition parameters $\Theta_k$. We provide a very detailed computational scheme for our algorithm in the specific case of the Group Graphical Lasso (GGL) penalty of \cite{danaher2014joint}. Moreover, we make the code for our algorithm, as well as a toy example that reproduces some of the results of this paper, publicly available at: \href{https://github.com/tlartigue/Mixture-of-Conditional-Gaussian-Graphical-Models}{https://github.com/tlartigue/Mixture-of-Conditional-Gaussian-Graphical-Models}.\\
Thanks to the inclusion of the co-features within the model, our EM algorithm is able to avoid trivial clusters correlated with the co-features' values, and instead unearths clusters providing new information on the population. Additionally, since our model takes into account heterogeneous effects of the co-features, our EM can handle the more complex scenarios, where the co-features act differently on the features in each sub-population.\\
A similar method was introduced by \cite{kim2016sparse}, who also uses a Conditional Gaussian Mixture model between co-features and features, as well as sparse GGM on the features. In their case however, the prediction of the features by the co-features takes centre stage. Indeed, the Mixture model acts as a Mixture of Experts, see \cite{jordan1994hierarchical}, that improves the prediction by allowing non-linear relations. Whereas in this article, the co-features are not meant to predict the features, but instead to make the clustering of the population more sensible. Another difference is that our algorithm performs at each M step a joint estimation of the model parameters in order to identify the common covariance structure between the different groups. Another domain of research, the "Finite Mixture Regression models" (FMR) \cite{desarbo1988maximum, khalili2007variable}, exhibit some superficial similarities with the Mixture of CGGM, but are very different in scope. The FMR consist of several parallel linear regressions between co-features and features. Where the features are usually uni- or low-dimensional and their covariance structure is not in question. Prediction is once again the main focus. Any clustering derived from an FMR model has its different classes solely defined by different linear relations (different "slopes") between co-features and features. This is very different from our GGM approach, where a multidimensional feature vector is described with a hierarchical graphical model, and the co-feature are a tool to improve the clustering within the feature space.\\
We demonstrate the performance of our method on synthetic and real data. First with a 2-dimensional toy example, where we show the importance of taking into consideration the (heterogeneous) effects of co-features for the clustering. Then, in higher dimension, we demonstrate that our EM with Mixture of CGGM consistently outperforms, both in terms of classification and parameter reconstruction, the EM with a Mixture of GGM (used in \cite{gao2016estimation} and \cite{hao2017simultaneous}), as well as an improved Mixture of GGM EM, that takes into consideration a homogeneous co-feature effect. Finally, on real Alzheimer's Disease data, we show that our method is the better suited to recover clusters correlated with the diagnostic, from both MRI and Cognitive Score features. 

\section{Supervised Hierarchical GGM and CGGM}
In this section, we summarise the whys and wherefores of Gaussian Graphical Modelling: the simple models for homogeneous populations, as well as the hierarchical models for heterogeneous populations. First, we explore the classical Gaussian Graphical Models techniques to describe a vector of features $Y \in \R^p$, then we discuss the Conditional Gaussian Graphical Models implemented in the presence of additional co-features $X \in \R^{q}$.  For every parametric model, we call $\theta$ the full parameter, and $p_{\theta}$ the probability density function. Hence, in the example of a gaussian model $\theta = (\mu, \Sigma)$. For hierarchical models with $K$ classes, we will have $K$ parameters $(\theta_1, ..., \theta_K)$.
\subsection{Basics of Hierarchical Gaussian Graphical Models}
In the classical GGM analysis introduced by \cite{dempster1972covariance}, the studied features $Y\in \R^p$ are assumed to follow a Multivariate Normal distribution: $Y \sim \N{\mu, \Sigma}$. The average $\mu$ is often ignored and put to 0. With $\Lambda:=\Sigma^{-1}$, the resulting distribution is:
\begin{equation} \label{eq:GGM_density}
    p_{\theta}(Y) = (2 \pi)^{-p/2} \det{\Lambda}^{1/2} exp\parent{-\frac{1}{2} Y^T \Lambda Y }\, .
\end{equation}
In this case $\theta = \Lambda$. Using the property that $corr(Y_u, Y_v | (Y_w)_{w \neq u,v}) = - \frac{(\Lambda)_{uv}}{\sqrt{(\Lambda)_{uu} (\Lambda)_{vv}}}$, the conditional correlation network is obtained using a sparse estimation of the precision matrix $\Lambda$. Heterogeneous populations, where different correlation networks may exist for each sub-population (or "class"), can be described with the Hierarchical version of the GGM \eqref{eq:GGM_density}. With $K$ classes, Let $\theta := (\theta_1, ..., \theta_k)$ be the parameter for each class and $z \in \llbracket 1, K \rrbracket$ the categorical variable corresponding to the class label of the observation $Y$. With $\theta_k:=\Lambda_k$ and $z$ known, the Hierarchical density can be written:
\begin{equation} \label{eq:hierarchical_density}
\begin{split}
    p_{\theta}(Y|z) &= \sum_{k=1}^K \mathds{1}_{z=k} p_{\theta_k}(Y) \\
                     &= \sum_{k=1}^K \mathds{1}_{z=k} (2 \pi)^{-p/2} \det{\Lambda_k}^{1/2} exp\parent{-\frac{1}{2} Y^T \Lambda_k Y } \, .
\end{split}
\end{equation}
Mirroring the famous Graphical LASSO (GLASSO) approach introduced by \cite{yuan2007model} and \cite{banerjee2006convex} for homogeneous populations, many authors have chosen to estimate sparse $\widehat{\Lambda}_k$ as penalised Maximum Likelihood Estimator (MLE) of $\Lambda_k$. For $i=1, ..., n$, let $Y^{(i)}$ be independent identically distributed (iid) feature vectors and $z^{(i)}$ their labels. These MLE are computed from the simple convex optimisation problem
\begin{equation} \label{eq:supervised_joint_GGM}
    \hat{\theta} = \argmin{\theta} - \frac{1}{n} \sum_{k=1}^K \sum_{i=1}^n \mathds{1}_{z^{(i)}=k} \, ln\,  p_{\theta_k}(Y^{(i)}) + pen(\theta)\, .
\end{equation}
Where the convex penalty $pen(\theta)$ is usually designed to induce sparsity within each individual $\hat{\Lambda}_k$ as well as to enforce a certain common structure between the $\hat{\Lambda}_k$. This common structure is a desirable outcome when the different sub-populations are assumed to still retain core similarities. Following in the footsteps of \cite{guo2011joint}, most authors propose such a joint estimation of the matrices $\Lambda_k$. In the case of the penalised MLE estimation \eqref{eq:supervised_joint_GGM}, the form of the resulting common structure is dependent on the penalty. For instance, \cite{danaher2014joint} propose the "Fused Graphical LASSO"" and "Group Graphical LASSO" penalties that encourage shared values and shared sparsity pattern across the different $\Lambda_k$ respectively. Likewise, \cite{yang2015fused} propose another fused penalty to incentivise common values across matrices. With their node based penalties, \cite{mohan2014node} can encourage the recovery of common hubs in the graphs. 
\begin{remark}
    Within a hierarchical model, one can also take $\theta_k:=(\mu_k, \Lambda_k)$, and adapt $p_{\theta_k}(Y)$ accordingly, since it is natural to allow each sub-population to have different average levels $\mu_k$.
\end{remark}

\subsection{Conditional GGM in the presence of co-features}
In some frameworks, additional variables, noted $X\in R^q$ and called "co-features" or "cofactors" can be observed alongside the regular features within the gaussian vector $Y\in \R^p$. In all generality, $X$ can be a mix of finite, discrete and continuous random variables. In the GGM analysis, these co-features are not included as nodes of the estimated conditional correlation graph. Instead, they serve to enrich the conditioning defining each edge: in the new graph, there is an edge between the nodes $Y_u$ and $Y_v$ if and only if $cov(Y_u, Y_v | (Y_w)_{w \neq u,v}, X) \neq 0$. The Conditional Gaussian Graphical Models (CGGM) were introduced by \cite{yin2011sparse} and \cite{sohn2012joint} in order to properly take into account the effect of $X$ on $Y$ and easily identify the new conditional correlation network in-between the $Y$. They propose a linear effect, expressed by the conditional probability density function (pdf):
\begin{equation} \label{eq:CGGM_density}
    p_{\theta}(Y | X) := \frac{\det{\Lambda}^{\frac{1}{2}}}{(2 \pi)^{\frac{p}{2}}} exp\bigg(-\frac{1}{2} (Y+ \Lambda^{-1} \Theta^T X)^T \Lambda (Y+ \Lambda^{-1} \Theta^T X) \bigg)\, ,
\end{equation}
with $\Theta\in \R^{q \times p}$ and $\theta = \brace{\Lambda, \Theta}$. In other words: $Y|X \sim \N{- \Lambda^{-1} \Theta^T X , \Lambda^{-1}}$. To simplify later equations, we define the function $L(Y, X, \Lambda, \Theta) := -\frac{1}{2} (Y+ \Lambda^{-1} \Theta^T X)^T \Lambda (Y+ \Lambda^{-1} \Theta^T X)$. Eq~\eqref{eq:CGGM_density} can be re-written:
\begin{equation*}
    p_{\theta}(Y | X) = \frac{\det{\Lambda}^{\frac{1}{2}}}{(2 \pi)^{\frac{p}{2}}} exp \Big(L(Y, X, \Lambda, \Theta) \Big)\,.
\end{equation*}
Two main branches of CGGM exist, depending on whether the pdf of $X$ is also modelled. In this work, we chose to impose no model on $X$. The lack of assumption on the density of $X$ provides far more freedom than the joint gaussian assumption. In particular, X can have categorical and even deterministic components. This allows us to integrate any observed variables without restriction to the model.\\
To tackle heterogeneous populations, works such as \cite{chun2013joint} have introduced the Hierarchical version of the CGGM pdf: 
\begin{equation} \label{eq:CGGM_density_hierarchical}
    p_{\theta}(Y | X, z) := \sum_{k=1}^K \mathds{1}_{z=k} \frac{\det{\Lambda_k}^{\frac{1}{2}}}{(2 \pi)^{\frac{p}{2}}} exp\Big(L(Y, X, \Lambda_k, \Theta_k) \Big)\, .
\end{equation}
In particular, \cite{huang2018joint} have adapted the penalised MLE \eqref{eq:supervised_joint_GGM} to the Hierarchical CGGM density for some of the most popular GGM penalties. With an iid observed sample $(Y^{(i)}, X^{(i)}, z^{(i)})_{i=1}^n$, the corresponding penalised CGGM MLE can be written;
\begin{equation} \label{eq:supervised_joint_CGGM}
    \hat{\theta} = \argmin{\theta} - \frac{1}{n} \sum_{k=1}^K \sum_{i=1}^n \mathds{1}_{z^{(i)}=k} \, ln\,  p_{\theta_k}(Y^{(i)}|X^{(i)}) + pen(\theta)\, .
\end{equation}

\begin{remark}
    To include a regular average value for $Y$, independent of the values of $X$, one can simply add a constant component equal to "1" in $X$.
\end{remark}

\section{Mixtures of CGGM for unlabelled heterogeneous population}\label{sect:Mixture_CGGM}
In this section, we tackle the problem of an unlabelled heterogeneous population. We introduce a Mixture of Conditional Gaussian Graphical Model to improve upon the state of the art unsupervised methods by taking into consideration the potent co-features that can drive the clustering. We develop a penalised EM algorithm to both identify data clusters and estimate sparse, structured, model parameters. We justify that our algorithm is usable with a wide array of penalties and provide detailed algorithmic steps for the Group Graphical LASSO (GGL) penalty.
\subsection{Presentation and motivation of the model}
When the labels of a heterogeneous population are missing, supervised parameter estimation methods like \eqref{eq:supervised_joint_GGM} have to be replaced by unsupervised approaches that also tackle the problem of cluster discovery. When $z$ is unknown, the Hierarchical model \eqref{eq:hierarchical_density} can easily be replaced by a Mixture model with observed likelihood:
\begin{equation} \label{eq:GGM_observed_density_mixture}
     p_{\theta, \pi}(Y) = \sum_{k=1}^K \pi_k p_{\theta_k}(Y) \, ,
\end{equation}
and complete likelihood:
\begin{equation} \label{eq:GGM_complete_density_mixture}
     p_{\theta, \pi}(Y, z) = \sum_{k=1}^K \mathds{1}_{z=k} \pi_k p_{\theta_k}(Y) \,.
\end{equation}
Where $\pi_k := \mathbb{P}(z=k)$ and $\pi:=(\pi_1, ..., \pi_k)$ . Then, the supervised penalised log-likelihood maximisation \eqref{eq:supervised_joint_GGM} can be adapted into the penalised observed log-likelihood optimisation:
\begin{equation} \label{eq:unsupervised_joint_GGM}
    \hat{\theta}, \hat{\pi} = \argmin{\theta, \pi} - \frac{1}{n}\sum_{i=1}^n ln\parent{\sum_{k=1}^K  \pi_k \,   p_{\theta_k}\parent{Y^{(i)}}} + pen(\theta, \pi)\, .
\end{equation}
This is a non-convex problem, and authors such as \cite{zhou2009penalized} and \cite{krishnamurthy2011high} have proposed EM algorithms to find local solutions to \eqref{eq:unsupervised_joint_GGM}. For a more recent example, \cite{hara2018network} used such a Mixture of GGM with an $l_1$ regularisation on each class with each M step solved by Graphical LASSO, \cite{banerjee2008model, friedman2008sparse}. Likewise, \cite{fop2019model} also consider problems of the form \eqref{eq:unsupervised_joint_GGM}, solved with an EM algorithm, although their regularisation functions belong to another family entirely, that enforces sparsity on the covariance matrix instead of the precision matrix. However, all these works omit the common structure inducing penalties that are the signature of the supervised joint GGM methods. The works of \cite{gao2016estimation} and \cite{hao2017simultaneous} correct this by proposing EM algorithms that solve \eqref{eq:unsupervised_joint_GGM} for some of the joint-GGM penalties, such as the Fused and Group Graphical LASSO penalties.\\
By design, the EM algorithm must handle the cluster identification jointly with the mixture parameters estimation. The underlying assumption is that the different sub-populations can be identified as different clusters in the feature space. With real data, and especially medical data, this is generally untrue, as many factors other than the class label can have a larger impact on the position of the data points in the feature space. Even when there are no specific sub-populations to recover, and the EM is ran "blindly" in order to observe which data points are more naturally grouped together by the method, the unearthed clusters have every chance to be very correlated with very influential but trivial external variables, such as the age group or the gender. In order to guide the cluster discovery of the EM algorithm, we propose a Mixture of Conditional Gaussian Graphical Models with which the overbearing effect of trivial external variables can be removed. By placing all external observed variable into $X$, we define the Mixture of CGGM with its observed likelihood:
\begin{equation} \label{eq:CGGM_density_mixture}
\begin{split}
    p_{\theta, \pi}(Y | X) :=& \sum_{k=1}^K \pi_k p_{\theta_k}(Y|X) \\
    =& \sum_{k=1}^K \pi_k \frac{\det{\Lambda_k}^{\frac{1}{2}}}{(2 \pi)^{\frac{p}{2}}} exp\Big(L(Y, X, \Lambda_k, \Theta_k) \Big)\, .
\end{split}
\end{equation}
Within this model, the position of each feature vector $Y$ is corrected by its, class-dependent, linear prediction by the co-features $X$: $\mathbb{E}[Y | X, z=k] = - \Lambda_k^{-1} \Theta_k^T X$. In other words the "Mixture of Gaussians" type clustering is done on the residual vector $Y - \mathbb{E}[Y | X, z=k] = Y +  \Lambda_k^{-1} \Theta_k^T X$. Hence, even if the co-features $X$ have a class-dependent impact on the average level of the features $Y$, the Mixture of CGGM model is still able to regroup in the feature space the observations $Y^{(i)}$ that belong to the same class, $z^{(i)}=k$. We illustrate this dynamic in section \ref{sect:experiment_2D}.\\
Like the previous works on joint-GGM estimation, our goal is to estimate the parameters of model \eqref{eq:CGGM_density_mixture} with sparse precision matrices $\Lambda_k$ and common structure across classes. Sparsity in the matrices $\Theta_k$ is also desirable for the sake of interpretation. Hence, we define the following penalised Maximum Likelihood problem:
\begin{equation} \label{eq:unsupervised_joint_CGGM}
\begin{split}
    \hat{\theta}, \hat{\pi} = \argmin{\theta, \pi} \Bigg\{&-\frac{1}{n}\sum_{i=1}^n ln\parent{\sum_{k=1}^K  \pi_k \,   p_{\theta_k}\parent{Y^{(i)} | X^{(i)}}}\\ 
    &+ pen(\theta, \pi) \Bigg\}\, .
\end{split}
\end{equation}
As with \eqref{eq:unsupervised_joint_GGM}, this is a non-convex problem, and we define an EM algorithm to find local minima of the optimised function.

\subsection{Penalised EM for the Mixture of CGGM}
In this section, we provide the detailed steps of a penalised EM algorithm to find local solution of the non-convex penalised MLE \eqref{eq:unsupervised_joint_CGGM} in order to estimate the parameters of the mixture model \eqref{eq:CGGM_density_mixture} with precision matrix sparsity as well as common structure. First we provide the different steps of the algorithm and justify that it can be run with a wide array of penalty functions. Then, we provide a detailed optimisation scheme for the Group Graphical Lasso (GGL) penalty specifically. The full code for this GGL-penalised CGGM EM can be found at \href{https://github.com/tlartigue/Mixture-of-Conditional-Gaussian-Graphical-Models}{https://github.com/tlartigue/Mixture-of-Conditional-Gaussian-Graphical-Models}.

\subsubsection{EM algorithm for Mixtures of CGGM} With $n$ fixed $\brace{X^{(i)}}_{i=1}^n$ and $n$ iid observations $\brace{Y^{(i)}}_{i=1}^n$ following the mixture density $p_{\theta, \pi}(Y | X)$ given in \eqref{eq:CGGM_density_mixture}, the penalised observed negative log-likelihood to optimise is:
\begin{equation} \label{eq:CGGM_penalised_observed_likelihood}
     - \frac{1}{n}\sum_{i=1}^n ln\parent{\sum_{k=1}^K  \pi_k \,   p_{\theta_k}\parent{Y^{(i)} | X^{(i)}}} + pen(\theta, \pi)\, .
\end{equation}
The EM algorithm is an iterative procedure updating the current parameter $(\theta^{(t)}, \pi^{(t)})$ with two steps. For a Mixture model, the Expectation (E) step is:
\begin{equation*}
\begin{split}
    p_{i, k}^{(t)} &:= \mathbb{P}_{\theta^{(t)}, \pi^{(t)}}(z^{(i)} = k | Y^{(i)}, X^{(i)})\\
    &= \frac{ p_{\theta_k^{(t)}}(Y^{(i)}|X^{(i)}) \pi_k^{(t)} }{\sum_{l=1}^K   p_{\theta_l^{(t)}}(Y^{(i)}|X^{(i)}) \pi_l^{(t)}} \, .
\end{split}
\end{equation*}
More explicitly, by replacing $p_{\theta_k}(Y|X)$ by its formula \eqref{eq:CGGM_density}:
\begin{equation} \label{eq:E_step_CGGM}
     (E) \quad p_{i, k}^{(t)} = \frac{\det{\Lambda_k}^{\frac{1}{2}} exp\parent{(L(Y^{(i)}, X^{(i)}, \Lambda_k, \Theta_k)} \pi_k^{(t)}}{\sum_{l=1}^K \det{\Lambda_l}^{\frac{1}{2}} exp\parent{L(Y^{(i)}, X^{(i)}, \Lambda_l, \Theta_l)} \pi_l^{(t)} } \, .
\end{equation}
Likewise, for Mixture models, the M step is:
\begin{equation*}
\begin{split}
    (\theta, \pi)^{(t+1)} = \argmin{\theta, \pi}\bigg\{&-\frac{1}{n} \sum_{k=1}^K \sum_{i=1}^n  p_{i, k}^{(t)}\,  ln \big(\pi_k  p_{\theta_k}(Y^{(i)}|X^{(i)})\big)\\
    &+ pen(\theta, \pi) \bigg\} \, .
\end{split}
\end{equation*}
Assuming that there is no coupling between $\pi$ and $\theta$ in the penalty, i.e. $pen(\pi, \theta) = pen_{\pi}(\pi) + pen_{\theta}(\theta)$, then the two optimisations can be separated:
\begin{equation*} 
\begin{split}
    \theta^{(t+1)} = \argmin{\theta} \bigg\{&-\frac{1}{n} \sum_{k=1}^K \sum_{i=1}^n  p_{i, k}^{(t)} \, ln\,  p_{\theta_k}(Y^{(i)}|X^{(i)})\\ 
    &+ pen_{\theta}(\theta) \bigg\} \, , \\
    \pi^{(t+1)} = \argmin{\pi} \bigg\{&-\frac{1}{n} \sum_{k=1}^K \sum_{i=1}^n p_{i, k}^{(t)} \,  ln\, \pi_k + pen_{\pi}(\pi) \bigg\}\, .
\end{split}
\end{equation*}
We denote the sufficient statistics $n_k^{(t)} := \sum_{i=1}^n  p_{i, k}^{(t)}$, $S_{YY}^{k, (t)} := \frac{1}{n}\sum_{i=1}^n  p_{i, k}^{(t)} Y^{(i)}Y^{(i)\, T}$, $S_{YX}^{k, (t)} := \frac{1}{n} \sum_{i=1}^n  p_{i, k}^{(t)} Y^{(i)}X^{(i) \, T}$ and $S_{XX}^{k, (t)} :=  \frac{1}{n} \sum_{i=1}^n  p_{i, k}^{(t)} X^{(i)}X^{(i \,T)}$. Then, the M step can be formulated as:
\begin{equation}\label{eq:M_step_CGGM}
    \begin{split}
    \theta^{(t+1)} = \argmin{\theta} \Bigg\{& \sum_{k=1}^K \Bigg(\dotprod{\Lambda_k, \frac{S_{YY}^{k, (t)}}{2}} + \dotprod{\Theta_k, S_{YX}^{k, (t)}}\\
    & \hspace{0.87cm} + \dotprod{\Theta_k \Lambda_k^{-1} \Theta_k^T,  \frac{S_{XX}^{k, (t)}}{2}} \\ 
    (M)\hspace{2.5cm} & \hspace{0.87cm} - \frac{n_k^{(t)}}{n} ln(\det{\Lambda_k}) \Bigg) +pen_{\theta}(\theta) \Bigg\} \, ,  \\ 
    \pi^{(t+1)} = \argmin{\pi} \bigg\{&
    - \sum_{k=1}^K \frac{n_k^{(t)}}{n} \,  ln\, \pi_k + pen_{\pi}(\pi) \bigg\}\, .
    \end{split}
\end{equation}
The E step in Eq~\eqref{eq:E_step_CGGM} is in closed form. Likewise, the optimisation problem that defines $\pi^{(t+1)}$ in Eq~\eqref{eq:M_step_CGGM} is very simple, and with any reasonable penalty $pen_{\pi}$, solving it is not issue. For instance, usual regularisation of the penalty take either the form $pen_{\pi}(\pi) = 0$ or $pen_{\pi}(\pi) \propto \sum_k \log(\pi_k)$, which both result in closed form solutions for $\pi^{(t+1)}$. The M step update of the model parameter $\theta$ is a harder optimisation problem. However, as is usually the case with EM algorithms, the update of $\theta$ in the M step \eqref{eq:M_step_CGGM} has the same form as the supervised MLE \eqref{eq:supervised_joint_CGGM}. As a result, when the supervised Hierarchical CGGM problem \eqref{eq:supervised_joint_CGGM} is tractable, then our M step is tractable as well. Conveniently, this supervised CGGM estimation problem has already been studied and solved in \cite{huang2018joint} for the Group Graphical Lasso (GGL) penalty of \cite{danaher2014joint}. Furthermore, \cite{huang2018joint} show that the supervised negative log-likelihood is a convex function of $\theta$. As a consequence, as long as the penalty $pen_{\theta}$ is convex and differentiable, then our M step \eqref{eq:M_step_CGGM} is a convex optimisation problem with a differentiable objective function. Since those problems are solvable with gradient descent algorithm, this means that the proposed EM algorithm is tractable for a very wide array of penalties $pen_{\theta}$. In particular, this includes all the convex differentiable penalties, which covers a lot of popular GGM penalties for sparsity and joint structure, such as the ones from the works of \cite{danaher2014joint},  \cite{mohan2014node} and \cite{yang2015fused} for instance.\\
In order to provide an algorithm with more specific and detailed steps, we consider in the rest of the section the special case of the GGL penalty. The GGL penalty was noticeably used in the supervised case by \cite{huang2018joint}, who proposed a proximal gradient algorithm. Likewise, we can use a proximal gradient algorithm to compute the M step \eqref{eq:M_step_CGGM} of our EM algorithm. 

\subsubsection{Proximal gradient algorithm to solve the M step with the GGL penalty} The GGL penalty, introduced in \cite{danaher2014joint} and adapted to the hierarchical CGGM by \cite{huang2018joint}, can be written:
\begin{equation} \label{eq:GGL_penalty}
\begin{split}
    pen_{\theta}(\theta) := &\sum_{1 \leq i \neq j \leq p} \parent{\lambda_1^{\Lambda} \sum_{k=1}^K \det{\Lambda_k^{(ij)}} + \lambda_2^{\Lambda} \sqrt{\sum_{k=1}^K \parent{\Lambda_k^{(ij)}}^2}}\\
    &+ \hspace{-0.1cm} \sum_{\substack{1 \leq i \leq q\\ 1 \leq j \leq p}} \parent{\lambda_1^{\Theta} \sum_{k=1}^K \det{\Theta_k^{(ij)}} + \lambda_2^{\Theta} \sqrt{\sum_{k=1}^K \parent{\Theta_k^{(ij)}}^2}} \, .
\end{split}
\end{equation}
Unlike in \cite{huang2018joint}, where$\lambda_1^{\Lambda}=\lambda_1^{\Theta}$ and $\lambda_2^{\Lambda}=\lambda_2^{\Theta}$, we use different levels of penalisation for the parameters $\Lambda$ and $\Theta$, since both their scales and their desired sparsity level can be very different. This penalty borrows its design from the Group Lasso, see \cite{yuan2006model}, where the $l_1$ norm induces individual sparsity of each coefficient, and the $l_2$ induces simultaneous sparsity of groups of coefficients. In Eq.~\eqref{eq:GGL_penalty}, for each pair $(i,j)$ belonging to the relevant space, $\brace{\Lambda_k^{(ij)}}_{k=1}^K$ constitutes a group that can be entirely put to 0. This incites the algorithm to set a certain matrix coefficient to 0 over all $K$ classes. These common zeros constitute the common structure sought after by the GGL approach. In our CGGM case, the same can be said for the group $\brace{\Theta_k^{(ij)}}_{k=1}^K$. Regarding the theoretical analysis, we underline that the $l_2$ part of the penalty is not separable in a sum of $K$ different penalties, which forces a joint optimisation problem to be solved, even in the supervised framework.\\
We detail here how to solve the M step \eqref{eq:M_step_CGGM} with $pen_{\theta}(\theta)$ defined as in Eq~\eqref{eq:GGL_penalty}. We assume, as usual, that the optimisation in $\pi$ is both independent from the optimisation in $\theta = \brace{\Lambda_k, \Theta_k}_{k=1}^K$ and trivial. The function to minimise in $\theta$ at the M step is:
\begin{equation*}
\begin{split}
    f(\theta) := \sum_{k=1}^K\Bigg( &-\frac{n_k^{(t)}}{n} ln(\det{\Lambda_k}) + \dotprod{\Lambda_k, S_{YY}^{k, (t)}} + \dotprod{2\Theta_k, S_{YX}^{k, (t)}}\\
    &+ \dotprod{\Theta_k \Lambda_k^{-1} \Theta_k^T,  S_{XX}^{k, (t)}} \Bigg)+pen_{\theta}(\theta)\, .
\end{split}
\end{equation*}
As shown in \cite{huang2018joint}, this function is convex and infinite on the border of its set of definition and as a unique global minimum. We note $f(\theta) =: g(\theta) + pen_{\theta}(\theta)$ for the sake of simplicity. The proximal gradient algorithm, see \cite{combettes2011proximal}, is an iterative method, based on a quadratic approximation on $g(\theta)$, that benefits from theoretical convergence guarantees. If $\theta^{(s-1)}$ is the current state of the parameter within the proximal gradient iterations, then the next stage, $\theta^{(s)}$, is found by optimising the approximation:
\begin{equation}\label{eq:proximal_gradient_quadratic_approximation}
\begin{split}
        f\parent{\theta^{(s)}} &= f\parent{\theta^{(s-1)} + \theta^{(s)}-\theta^{(s-1)}}\\ 
        &\approx g\parent{\theta^{(s-1)}} + \nabla g\parent{\theta^{(s-1)}}^T . \parent{\theta^{(s)}-\theta^{(s-1)}}\\
        & \hspace{1cm} + \frac{1}{2 \alpha} \norm{\theta^{(s)}-\theta^{(s-1)}}_2^2 + pen_{\theta}\parent{\theta^{(s)}} \\
        &\equiv \frac{1}{2 \alpha} \norm{\theta^{(s)} -\parent{\theta^{(s-1)} - \alpha \nabla g\parent{\theta^{(s-1)}}  } }_2^2\\
        & \hspace{1cm} +  pen_{\theta}\parent{\theta^{(s)}} \, .
\end{split}
\end{equation}
Where we removed in the last line the constants irrelevant to the optimisation in $\theta^{(s)}$ and $\alpha$ denotes the step size of the gradient descend. Note that we use the superscript $(s)$ to indicate the current stage of the proximal gradient iteration, to avoid confusion with the superscript $(t)$ used for the EM iterations (which are one level above). We underline that, in addition to $g(\theta)$ itself, the second order term in the Taylor development of $g(\theta)$ is also approximated. Using $ \frac{1}{2 \alpha} \norm{\theta^{(s)}-\theta^{(s-1)}}_2^2$ instead of $\frac{1}{2} \parent{\theta^{(s)}-\theta^{(s-1)}}^T . H_g\parent{\theta^{(s-1)}} . \parent{\theta^{(s)}-\theta^{(s-1)}}$ spares us from computing the Hessian $H_g\parent{\theta^{(s-1)}}$ and simplifies the calculations to come. The approximated formulation in Eq~\eqref{eq:proximal_gradient_quadratic_approximation} leads to the definition of the proximal optimisation problem:
\begin{equation} \label{eq:proximal_problem}
	{prox}_{\alpha}(x) := \argmin{\theta} \frac{1}{2 \alpha} \norm{\theta-x}_2^2 + pen_{\theta}(\theta) \, .
\end{equation}
So that the proximal gradient step can be written:
\begin{equation} \label{eq:proximal_gradient_step}
	\theta^{(s)} = {prox}_{\alpha_s}\parent{\theta^{(s-1)} - \alpha_s \nabla g\parent{\theta^{(s-1)}}  } \, .
\end{equation}
Where the step size $\alpha_s$ is determined by line search. The usual proximal gradient heuristic is to take a initial step size $\alpha^{0}$, a coefficient $\beta\in  (0, 1)$, and to reduce the step size, $\alpha \longleftarrow \beta \alpha$, as long as:
\begin{equation*}
\begin{split}
    g\parent{\theta^{(s-1)}-\alpha G_{\alpha}\big(\theta^{(s-1)}\big) } >& g\parent{\theta^{(s-1)}}\\
    &- \alpha \nabla g\parent{\theta^{(s-1)}}^T. G_{\alpha}\big(\theta^{(s-1)}\big)\\
    &+ \frac{\alpha}{2} \norm{G_{\alpha}\big(\theta^{(s-1)}\big)}_2^2 \, , 
\end{split}
\end{equation*}
with $G_{\alpha}\parent{\theta^{(s-1)}} := \frac{\theta^{(s-1)} - {prox}_{\alpha} \parent{\theta^{(s-1)} - \alpha \nabla g\parent{\theta^{(s-1)}}  } }{\alpha}$ the generalised gradient.\\
To apply the proximal gradient algorithm, we need to be able to solve the proximal \eqref{eq:proximal_problem} with the CGGM likelihood and the GGL penalty. Thankfully, \cite{danaher2014joint} found an explicit solution to this problem in the GGM case, which \cite{huang2018joint} adapted to the CGGM. The proximal optimisation is separable in $\Lambda$ and $\Theta$, and the solutions $\Lambda^{(prox)}$ and $\Theta^{(prox)}$ share the same formula. As a result, we use $D$ as a placeholder name for either $\Lambda$ or $\Theta$, i.e. depending on the context either $D_k^{ij}=\Lambda_k^{ij}$ or $D_k^{ij}=\Theta_k^{ij}$. Let $S$ be the soft thresholding operator: $S(x, \lambda):=sign(x)\,  max(|x|-\lambda, 0)$, and $\widetilde{D}_{k, \alpha}^{ij} := D_k^{ij, (s-1)} - \alpha \frac{\partial g}{\partial D_k^{ij}}\parent{\theta^{(s-1)}}$. The solution of \eqref{eq:proximal_problem}, with $x=\theta^{(s-1)} - \alpha \nabla g\parent{\theta^{(s-1)}} $, is given coefficient-by-coefficient in Eq~\eqref{eq:proximal_solution}:
\begin{equation} \label{eq:proximal_solution}
\begin{split}
    D_k^{ij, (prox)} = &S\parent{\widetilde{D}_{k, \alpha}^{ij}, \lambda_1^{D} \alpha}\\
    &\times max\parent{1 - \frac{ \lambda_2^D \alpha}{\sqrt{\sum_k S(\widetilde{D}_{k, \alpha}^{ij}, \lambda_1^{D} \alpha)^2}}, 0} \, .
\end{split}
\end{equation}
Note that the partial derivatives $\frac{\partial g}{\partial D_k^{ij}}\parent{\theta^{(s-1)}}$, necessary to get $\widetilde{D}_{k, \alpha}^{ij}$, are easily calculated in closed form from the likelihood formula. With the proximal problem \eqref{eq:proximal_problem} and the line search easily solvable, the proximal gradient steps can be iterated until convergence to find the global minimum of $f(\theta)$. With $f(\theta)$ optimised, the M step \eqref{eq:M_step_CGGM} is solved.

\section{Experiments}
In this section, we demonstrate the performances of our EM with Mixture of CGGM. First on a visual toy example in 2 dimension, then on a higher dimensional synthetic example and finally on real Alzheimer's Disease data. We compare the Mixture of CGGM to the regular Mixture of GGM which ignores co-features and to a Mixture of GGM that assumes a uniform linear effect of the co-features on the features.

\subsection{An illustration of co-features with class-dependent effect} \label{sect:experiment_2D}
\begin{figure*}[tbhp]
    \centering
        \subfloat{\includegraphics[width = 0.49\textwidth]{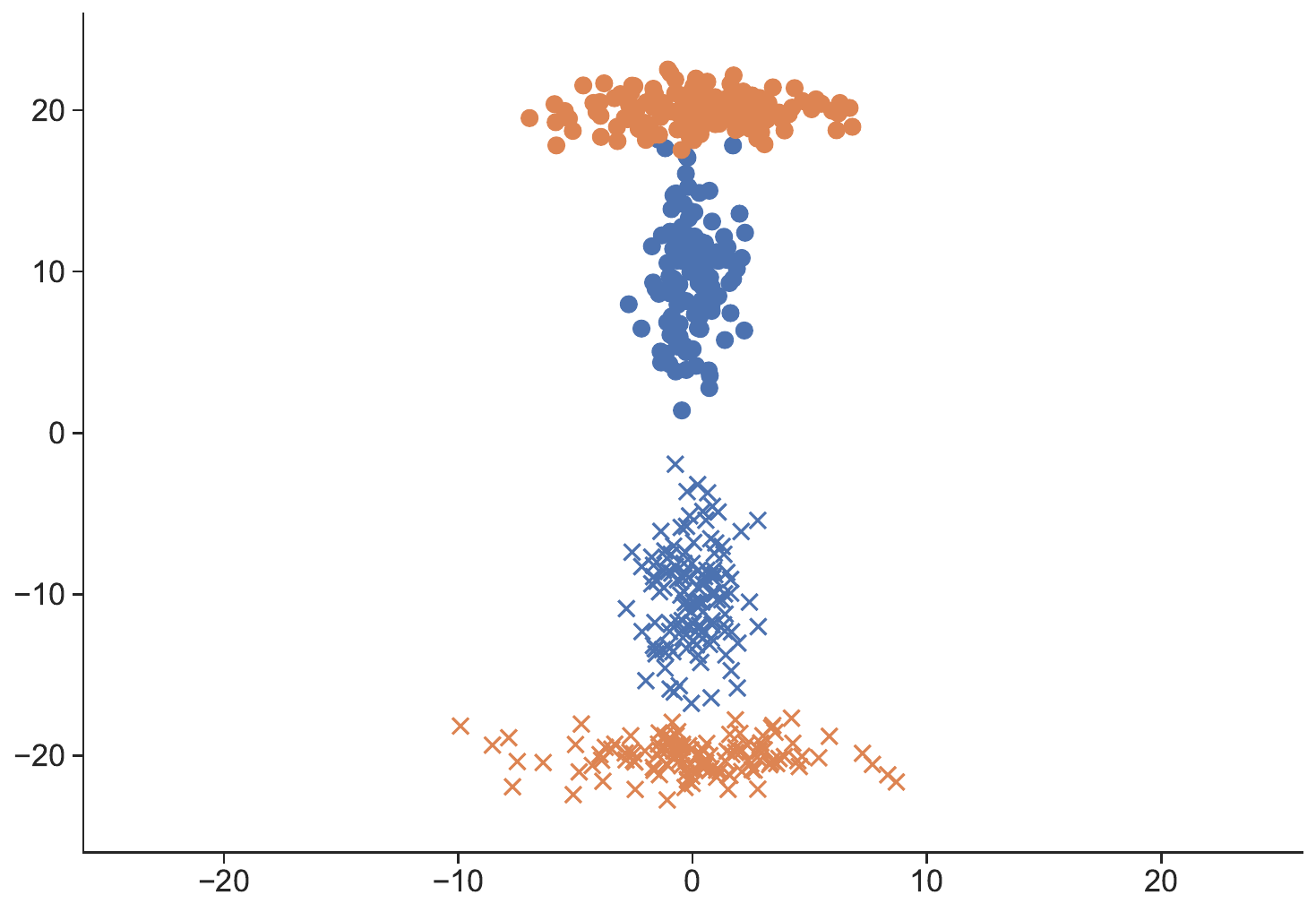}}
        \subfloat{\includegraphics[width=0.49\linewidth]{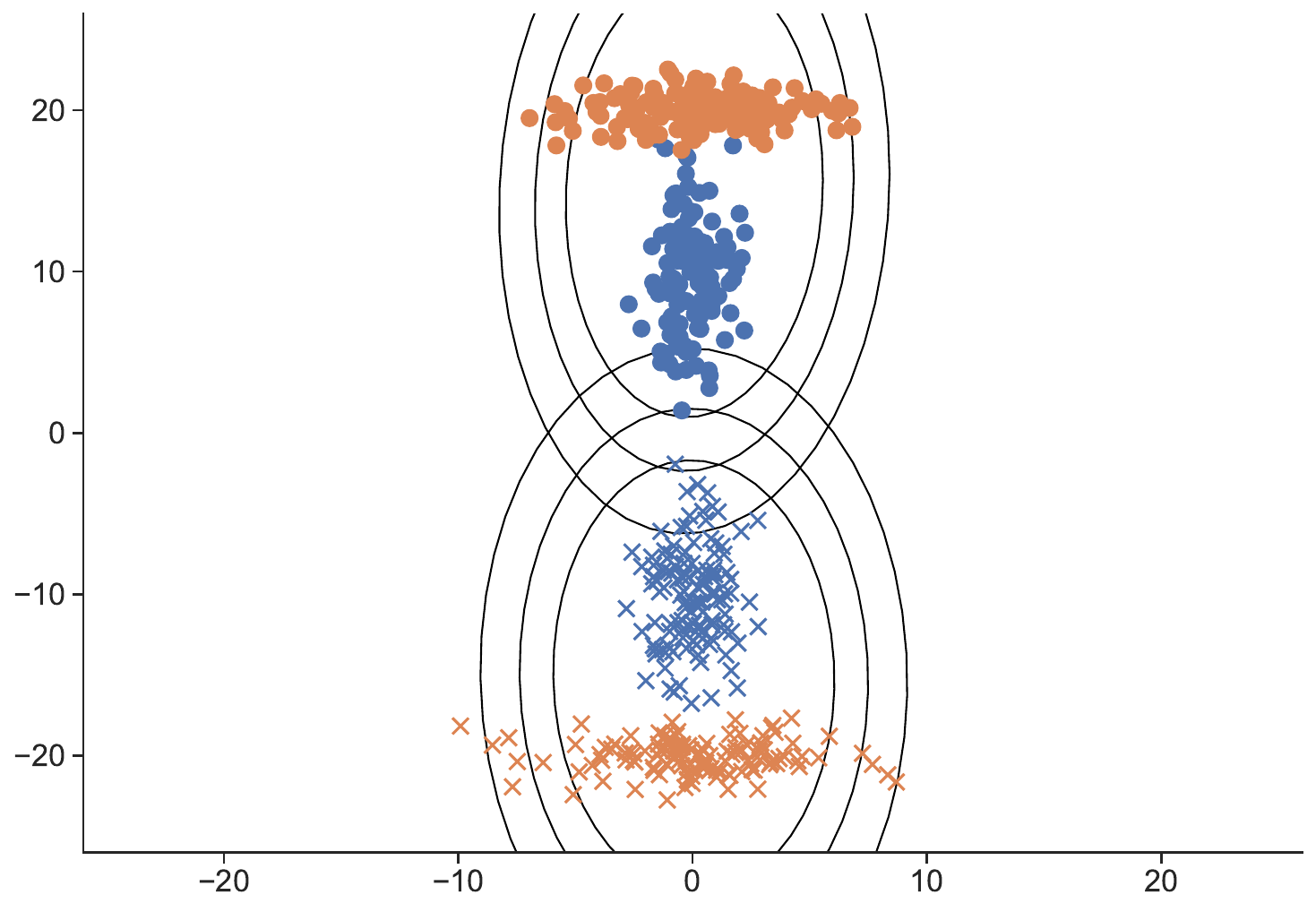}}
        \caption{(Left) Observed data $Y$ in the 2D space. The observed conditioning variable $X$ is binary. Data points with $X=-1$ are represented as crosses, and the ones with $X=1$ are represented as dots. In addition, there is an unknown "class" variable $z$. Class 1 is in blue, class 2 in orange. $Y | X, z$ follows the hierarchical conditional model. As a result, the two classes (orange and blue) are hard to separate in two clusters. (Right) Typical clusters estimated by an EM that fits a GGM mixture on $Y$} \label{fig:2D_EM}
\end{figure*}

\begin{figure*}[tbhp]
    \centering
    \subfloat{\includegraphics[width=0.49\linewidth]{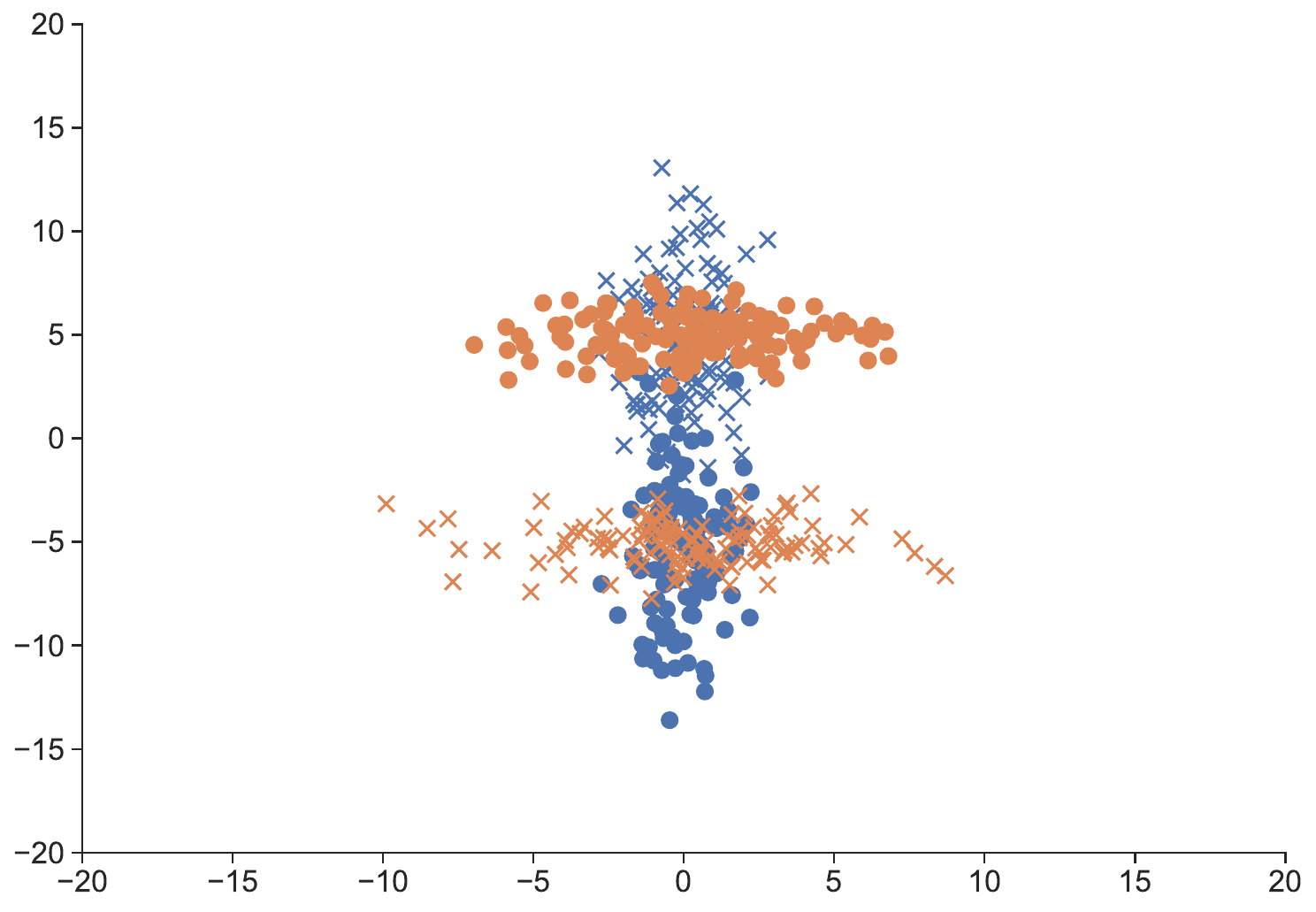}}
    \subfloat{\includegraphics[width=0.49\linewidth]{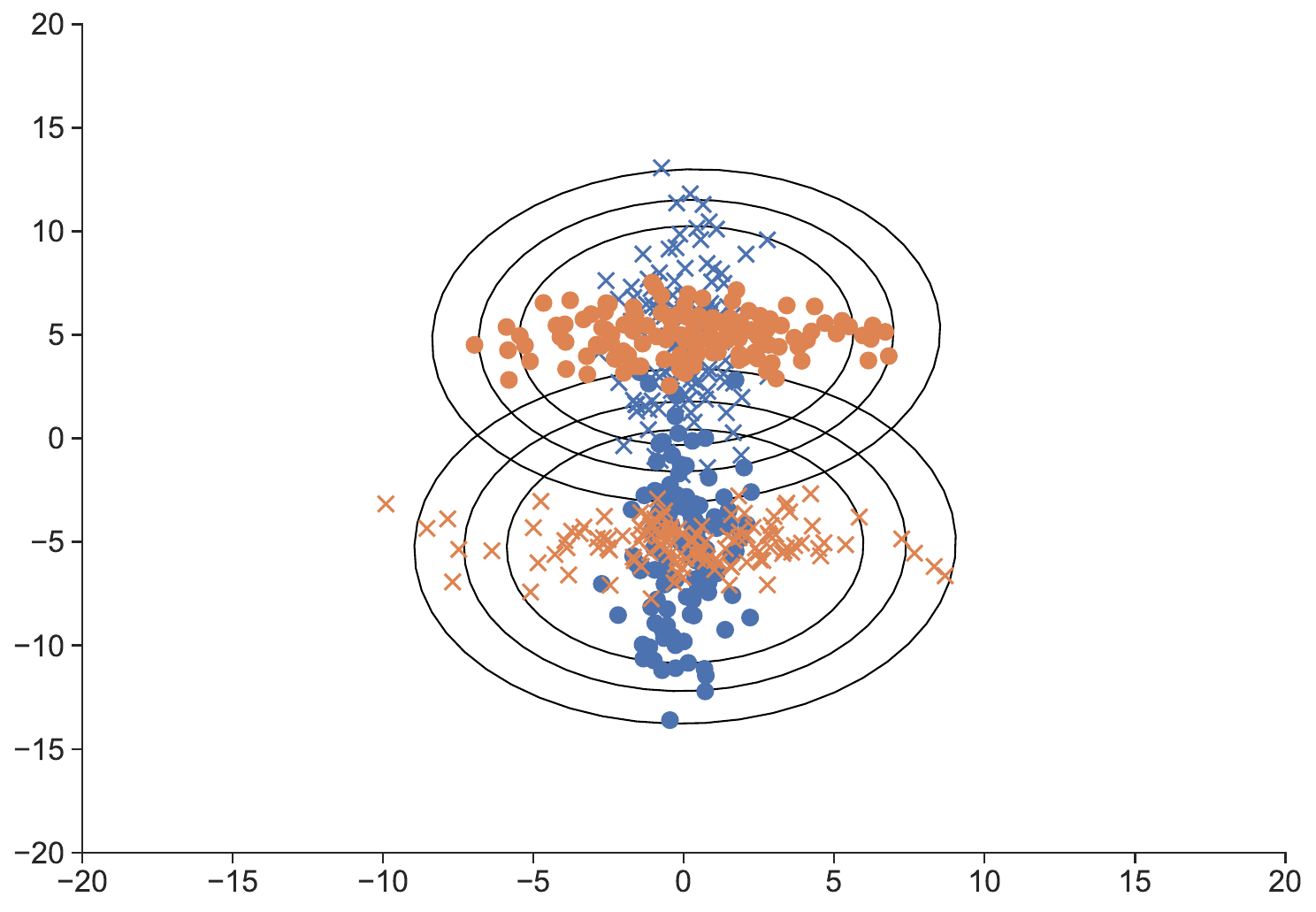}}
\caption{(Left) Residual $\tilde{Y} = Y-\hat{\beta} Y$ data after taking into account the estimated effect of $X$. Since the effect had different intensities on class 1 and 2, only the average effect was subtracted, and two classes are still not well separated. (Right) Typical clusters estimated by the "residual EM", that fits a GGM mixture on $\tilde{Y}$ }
\label{fig:2D_EMres}
\end{figure*}

\begin{figure*}[tbhp]
    \centering
    \subfloat{\includegraphics[width=0.32\linewidth]{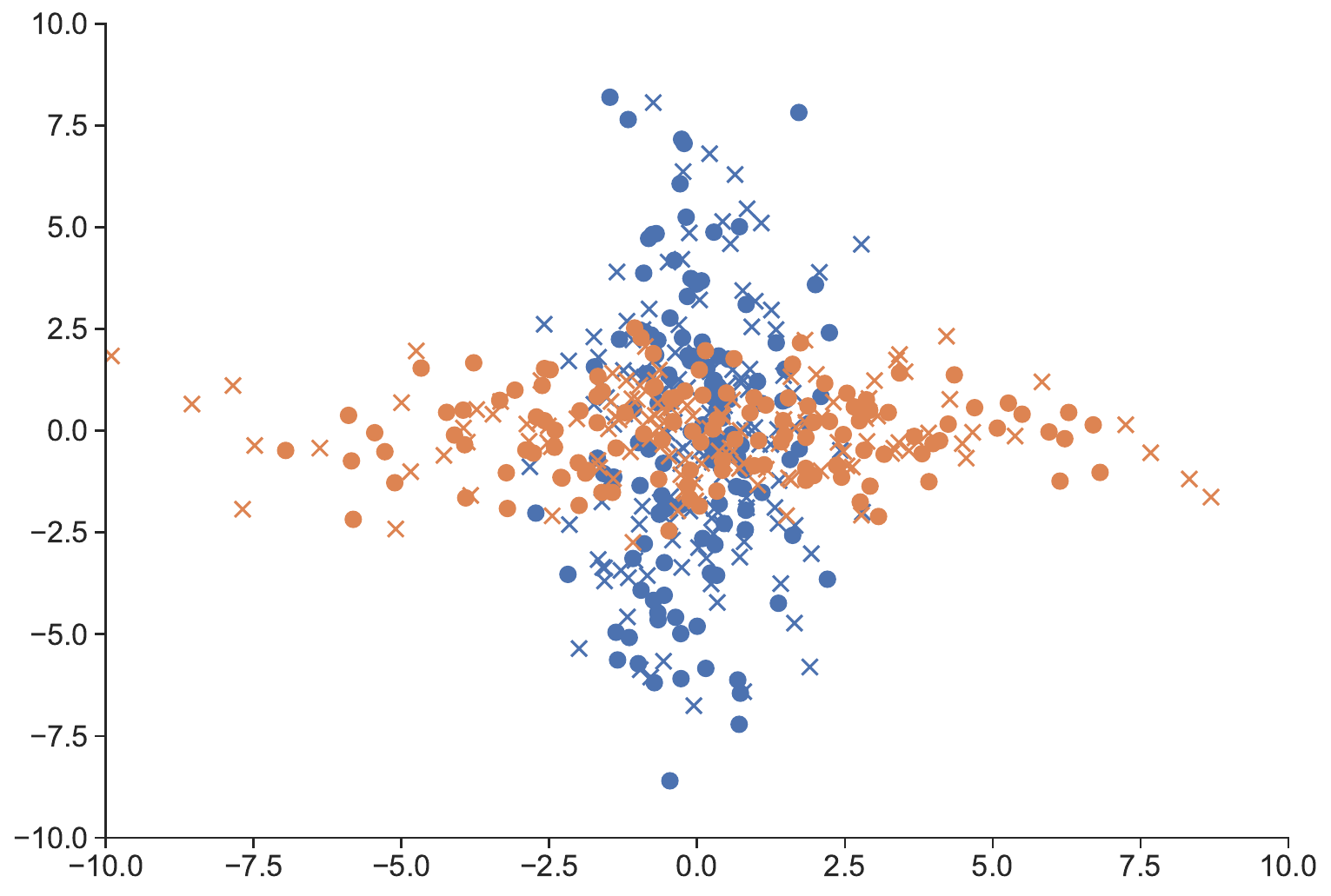}}
     \subfloat{\includegraphics[width=0.32\linewidth]{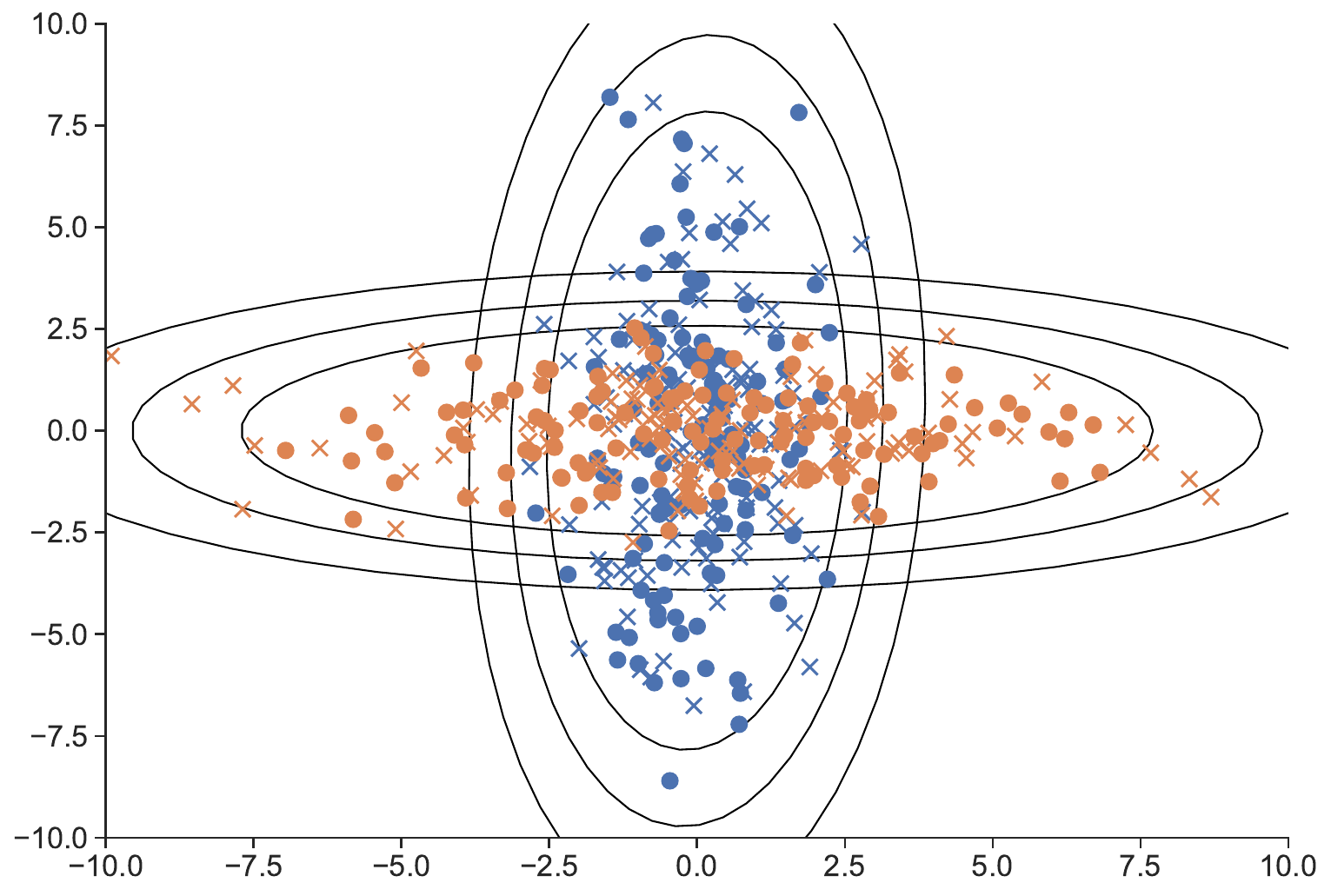}}
    \subfloat{\includegraphics[width=0.32\linewidth]{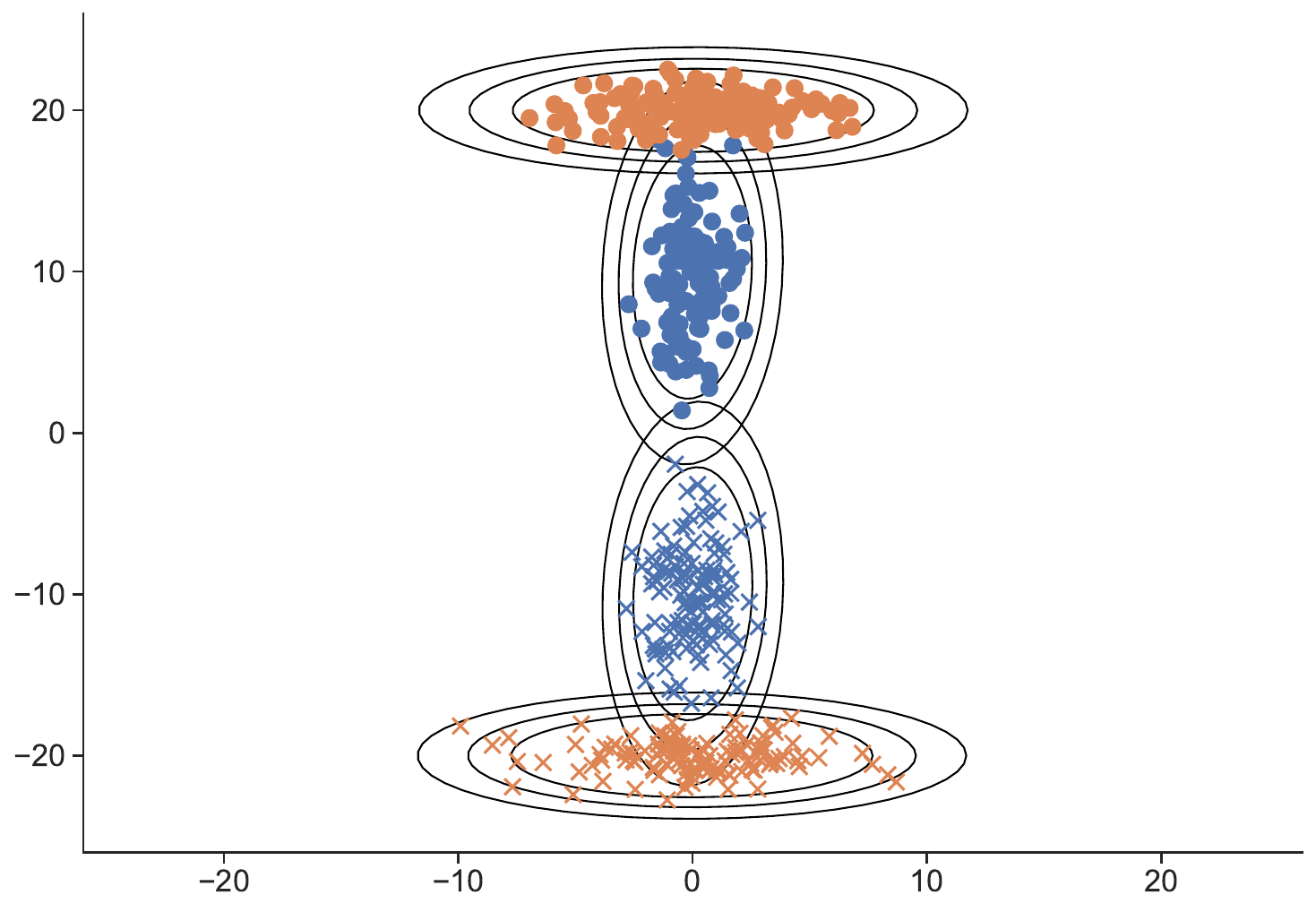}}
\caption{(Left) Observations $\tilde{Y}' = Y - \beta_1 X \mathds{1}_{z=1} - \beta_2 X \mathds{1}_{z=2}$ exactly corrected for the class-dependent effect of $X$. In this state the two classes appear as two distinct clusters. The Conditional-EM is designed to transform the data in this manner. (Middle) One possible representation of the CEM results. The corrected observations $\tilde{Y}'$ are displayed alongside centered normal distributions with the two estimated covariance matrices: $\mathcal{N}(0, \widehat{\Lambda}_k^{-1})$. (Right) Another possible representation of the same CEM results. The original observations $Y$ are displayed, alongside with the four \textit{de facto} estimated distributions $\mathcal{N}(\pm\widehat{\Sigma}_k \widehat{\Theta}_k^T, \widehat{\Lambda}_k^{-1})$.}
\label{fig:2D_CEM}
\end{figure*}

In this section, we present a simple visual example to illustrate the importance of taking into account heterogeneous co-feature effects. We show that even with a single binary co-feature, and with low dimensional features, the state of the art unsupervised GGM techniques are greatly disrupted by the co-features. Whereas our EM with Mixture of CGGM (which we call "Conditional EM" or "C-EM") achieves near perfect classification.\\
Under the Mixture of Gaussians (MoG) model, the observed data, $Y \sim \sum_{k=1}^K \pi_k \N{\mu_k, \Sigma_k}$, belongs to $K$ classes which can directly be represented as $K$ clusters in the feature space $\R^p$. Each cluster centered around a centroid at position $\mu_k$ and with an ellipsoid shape described by $\Sigma_k$. However, when there exists conditioning variables $X \in \R^q$ that have an effect on $Y$, this geometric description becomes more complex. Typically, the value of $Y$ could depend linearly on the value of $X$, with  $\E{Y | X, z = k} = \beta_k^T X$ for some $\beta_k \in \R^{p \times q}$. In this case, the average position in class $k$ is not a fixed $\mu_k$ but a function of $X$. If $X$ contains categorical variables, this creates as many different centroid positions as there are possible category combinations in $X$. The number of these \textit{de facto} clusters geometrically increases with the dimension $q$, which deters from simply running a clustering method with an increased number of clusters $K'$ to identify all of them. Moreover, if $X$ contains continuous variables, there is a continuum of positions for the centroid, not a finite number of \textit{de facto} clusters. If $X$ mixes the two types of variables, the two effects coexist. This shatters any hope to run a traditional MoG-based EM clustering algorithm, since its success is heavily dependent on its ability to identify correctly the $K$ distinct cluster centroids $\mu_k$.\\
Since the $X$ are observed, a possible solution is to run the linear regression $\hat{Y} = \hat{\beta} X$ beforehand, and run the EM algorithm on the residual $Y-\hat{Y}$ to remove the effect of $X$. This is what we call the "residual EM" or "residual Mixture of GGM". However this does not take into account the fact that this effect can be different for each class $k$, $\beta_1 \neq \beta_2 \neq ... \neq \beta_K$. Since the label is not known beforehand in the unsupervised context, the linear regression $\hat{Y} = \hat{\beta} X$ can only be run on all the data indiscriminately, hence is insufficient in general. On the other hand, the hierarchical CGMM \eqref{eq:CGGM_density_hierarchical}, which verifies: $\E{Y | X, z = k} = - \Lambda_k^{-1} \Theta_k^T X$, is designed to capture heterogeneous co-feature effects. We design a simple experiment to substantiate this intuition.\\
\\
In this example, $Y \in \R^2, X \in \brace{-1, 1}$ and $z \in \brace{1,2}$. $Y | X, z$ follows the hierarchical conditional model of \eqref{eq:CGGM_density_hierarchical}. In this simple case, this can be written as $Y = \parent{\beta_1 X +  \epsilon_1} \mathds{1}_{z=1} + \parent{\beta_2 X + \epsilon_2} \mathds{1}_{z=2}$. With $\P{z=1}=\P{z=2}=0.5$, $\epsilon_1 \sim \mathcal{N}(0, \Lambda_1^{-1})$ and $\epsilon_2 \sim \mathcal{N}(0, \Lambda_2^{-1})$. A typical iid data sample $(Y^{(i)}_{i=1})^n$ is represented on the left sub-figure of \figurename~\ref{fig:2D_EM}. The hidden variable $z$ is represented by the colour (blue or orange). The covariance of the cluster changes with $z$. The observable co-feature $X$ is represented by the shape of the data point (dot or cross). The centroid of the cluster is translated differently according to the values of $X$ and $z$. It is clear from the figure that a Mixture of Gaussians model with $K=2$ cannot properly separate the blue and orange points in two clusters. Indeed, on the right sub-figure of \figurename~\ref{fig:2D_EM}, we observe the final state of an EM that fits a Mixture of Gaussians on $Y$. The two recovered clusters are more correlated with the co-feature $X$ than the hidden variable $z$. However, this method did not take advantage of the knowledge of the co-feature $X$. As previously mentioned, one could first subtract the effect of $X$ from $Y$ before running the EM. On the left sub-figure of \figurename~\ref{fig:2D_EMres}, we represent the residual data $\tilde{Y} :=Y - \hat{\beta} X$. Where $\hat{\beta}$ is the Ordinary Least Square estimator of the linear regression between $X$ and $Y$ over all the dataset ($\hat{\beta} \approx \frac{\beta_1+\beta_2}{2}$ if $n$ is large enough). Since the linear effect between $X$ and $Y$ is not uniform over the dataset, but class dependent, the correction is imperfect, and the two class clusters remain hardly separable. This is why the residual EM, that fits a Mixture of GGM on $\tilde{Y}$ is also expected to fail to identify clusters related to the hidden variable. Which is shown by the right sub-figure of \figurename~\ref{fig:2D_EMres}, where we see a typical final state of the residual EM.\\
On the leftmost sub-figure of \figurename~\ref{fig:2D_CEM}, we display the proper correction for the co-features' effect $\tilde{Y}' = Y - \beta_1 X \mathds{1}_{z=1} - \beta_2 X \mathds{1}_{z=2} = \epsilon_1 \mathds{1}_{z=1}+\epsilon_2\mathds{1}_{z=2}$. Under this form, a Mixture of Gaussian can separate the data by colour. This is precisely the kind of translation that each data point undergoes within a Hierarchical CGGM. Hence a Mixture of CGGM can succeed in identifying the hidden variable $z$, provided that it estimates correctly the model parameters. To illustrate this point, the two next sub-figures in \figurename~\ref{fig:2D_CEM} represent the same final state of the EM fitting a Mixture of CGGM on $Y$. The middle sub-figure represents $\tilde{Y}'$ as well as the two estimated centered distributions $\mathcal{N}(0, \widehat{\Lambda}_k^{-1})$ for $k=1, 2$. We can see the two formally identified clusters after removing the effect of $X$. The rightmost sub-figure represents the original data $Y$ as well as the four estimated distributions $\mathcal{N}(\pm\widehat{\Sigma}_k \widehat{\Theta}_k^T, \widehat{\Lambda}_k^{-1})$ for $k=1, 2$. The four \textit{de facto} clusters present in the data $Y$ before removing the effect of $X$ are well estimated by the method.\\
We confirm these illustrative results by running several simulations. Under the previously introduced framework, we generate 50 datasets with $n=500$ data points each. For each simulation, we make 10 random initialisations from which we run the three EMs: with GGM, residualised GGM or CGGM. Table \ref{table:2D_example_CEM} summarises the results of these simulations. We follow the errors made by the estimated class probabilities or "soft labels", $\widehat{\mathbb{P}}(z_i=k)$, which we call the "soft misclassification error", as well as the error made by the "hard labels", $\mathds{1}_{\widehat{z}_i=k}$, which we call the "hard misclassification error". They can be expressed as $\frac{1}{2 n} \sum_{i, k} \det{\mathds{1}_{z_i=k}  - \widehat{\mathbb{P}}(z_i=k)}$ and $\frac{1}{2 n} \sum_{i, k}\det{\mathds{1}_{z_i=k}  - \mathds{1}_{\widehat{z}_i=k}}$ respectively. We see that the Mixture of CGGM performs much better, with less than $10 \%$ of misclassification in average, while the two GGM methods are both above $40\%$ of error, fairly close to the level of a random uniform classifier, $50\%$.
\begin{table}[tbhp]
\begin{center}
\caption{Average and standard deviation of the misclassification error achieved on the 2-dimensional example with the EMs on the Mixture of GGM, the Mixture of GGM with residualised data, and the Mixture of CGG. The two GGM methods are close to the threshold of random classification (0.50), while the Mixture of CGGM is in average below $10\%$ of error.}
\begin{tabular}{l|ccc}
        \toprule
         & EM GGM & EM res. GGM & EM CGGM \\
        \midrule \\
        soft misclassif. & 0.41 (0.11) & 0.47 \textbf{(0.05)} & \textbf{0.08} (0.17) \\
        hard misclassif. & 0.41 (0.12) & 0.46 \textbf{(0.06)} & \textbf{0.07} (0.17) \\ 
        \bottomrule
\end{tabular}
\end{center}
\label{table:2D_example_CEM}
\end{table}

\subsection{Experiments in high dimension}
\begin{figure}[tbhp]
    \centering
    \subfloat{\includegraphics[width = 0.22\textwidth]{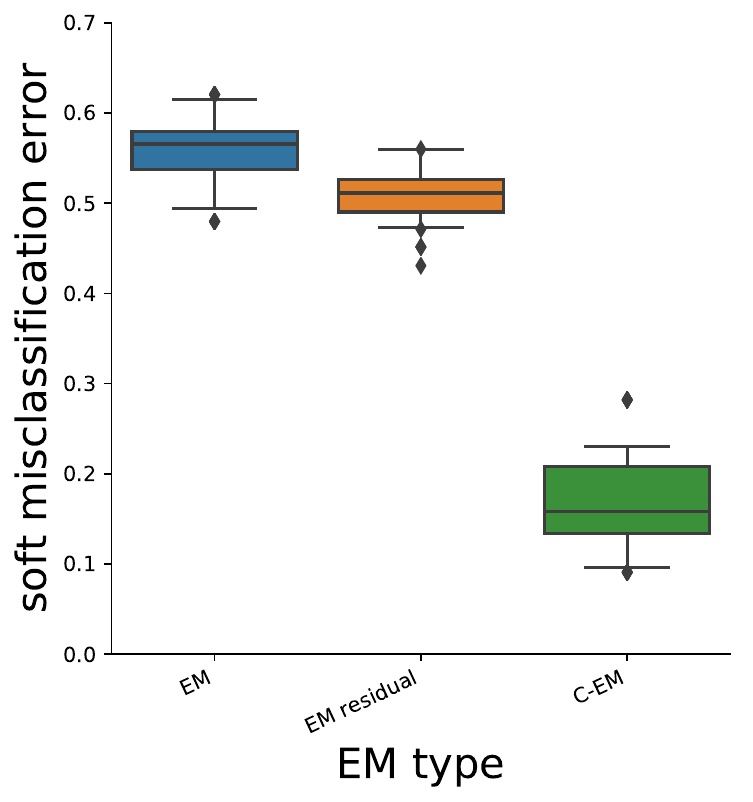}}
    \subfloat{\includegraphics[width = 0.22\textwidth]{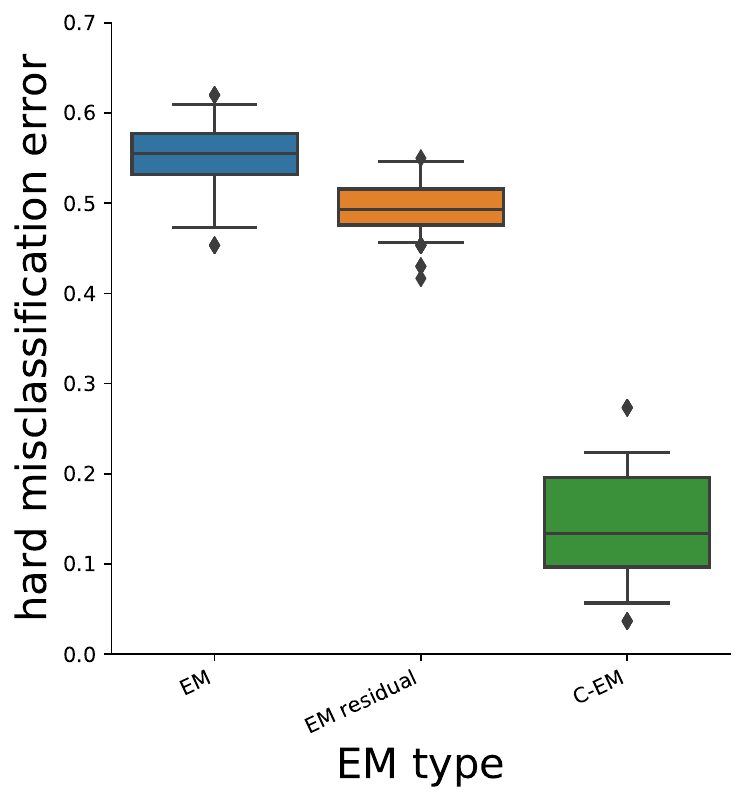}}\\
    \subfloat{\includegraphics[width = 0.22\textwidth]{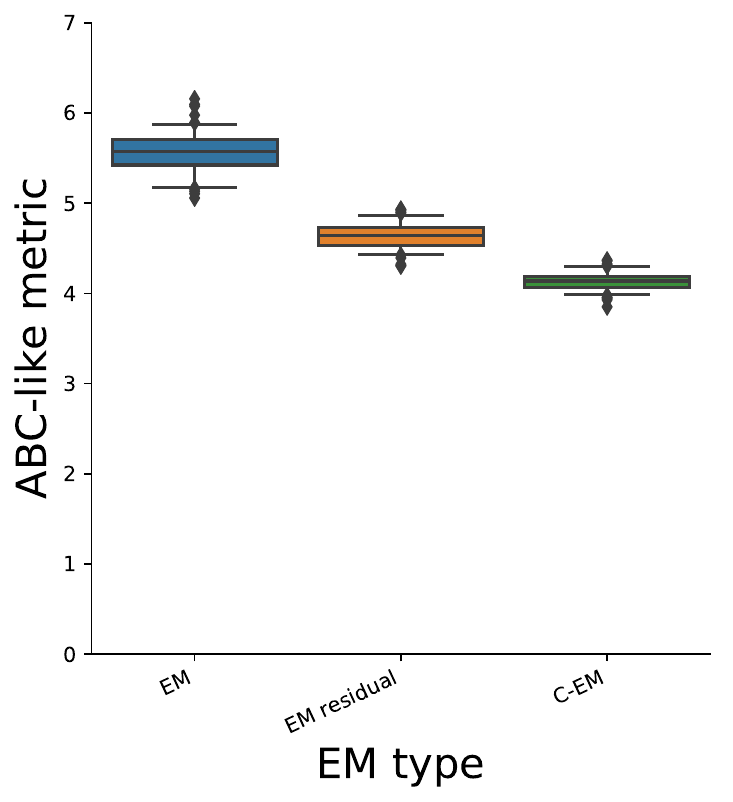}}
    \subfloat{\includegraphics[width = 0.22\textwidth]{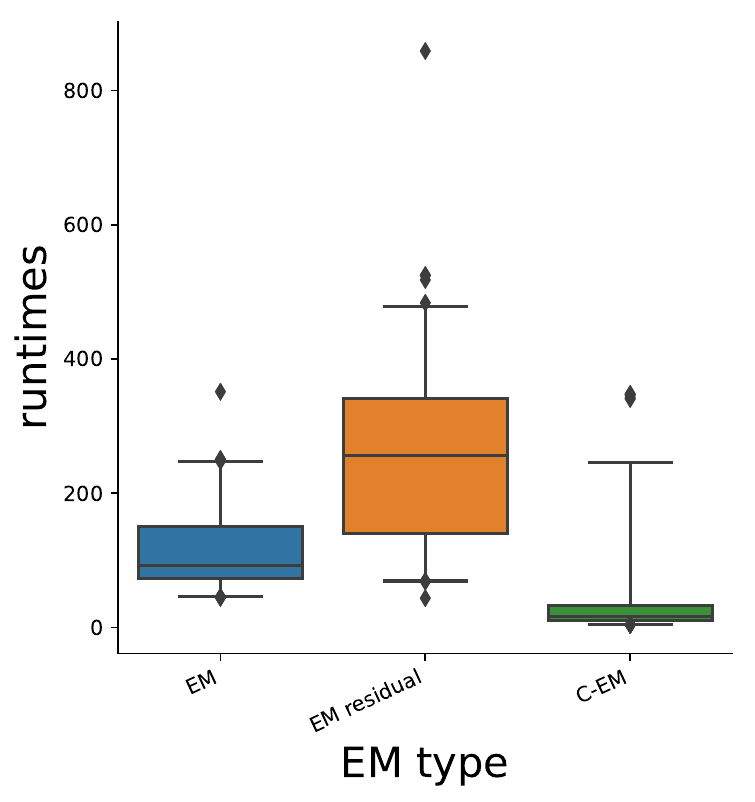}}
    \caption{Empirical distribution of several performance metrics measured over many simulations. The sample is made of 20 simulations with 5 different initialisations each. Three methods are compared. The EM and EM residual algorithms estimate a Mixture of GGM. The C-EM algorithm estimates a Mixture of CGGM. The C-EM is much better performing and faster. (Upper left) Soft mis-classification error $\det{\mathds{1}_{z_i=k}  - \widehat{\mathbb{P}}(z_i=k)}$. (Upper right) Hard mis-classification error $\det{\mathds{1}_{z_i=k}  - \mathds{1}_{\widehat{z}_i=k}}$. (Bottom left) ABC-like metric. (Bottom right) Run time.}
    \label{fig:HD_all_metrics}
\end{figure}

\begin{figure}[tbhp]
    \centering
    \begin{tabular}{cc}
    \adjustbox{valign=c, scale=0.8}{
    \begin{tabular}{lccc}
    & \footnotesize class 1 & \footnotesize class 2 & \footnotesize class 3\\
    \footnotesize EM &
    \adjustbox{valign=c, scale=0.8}{
    \subfloat{\includegraphics[width=0.32\linewidth]{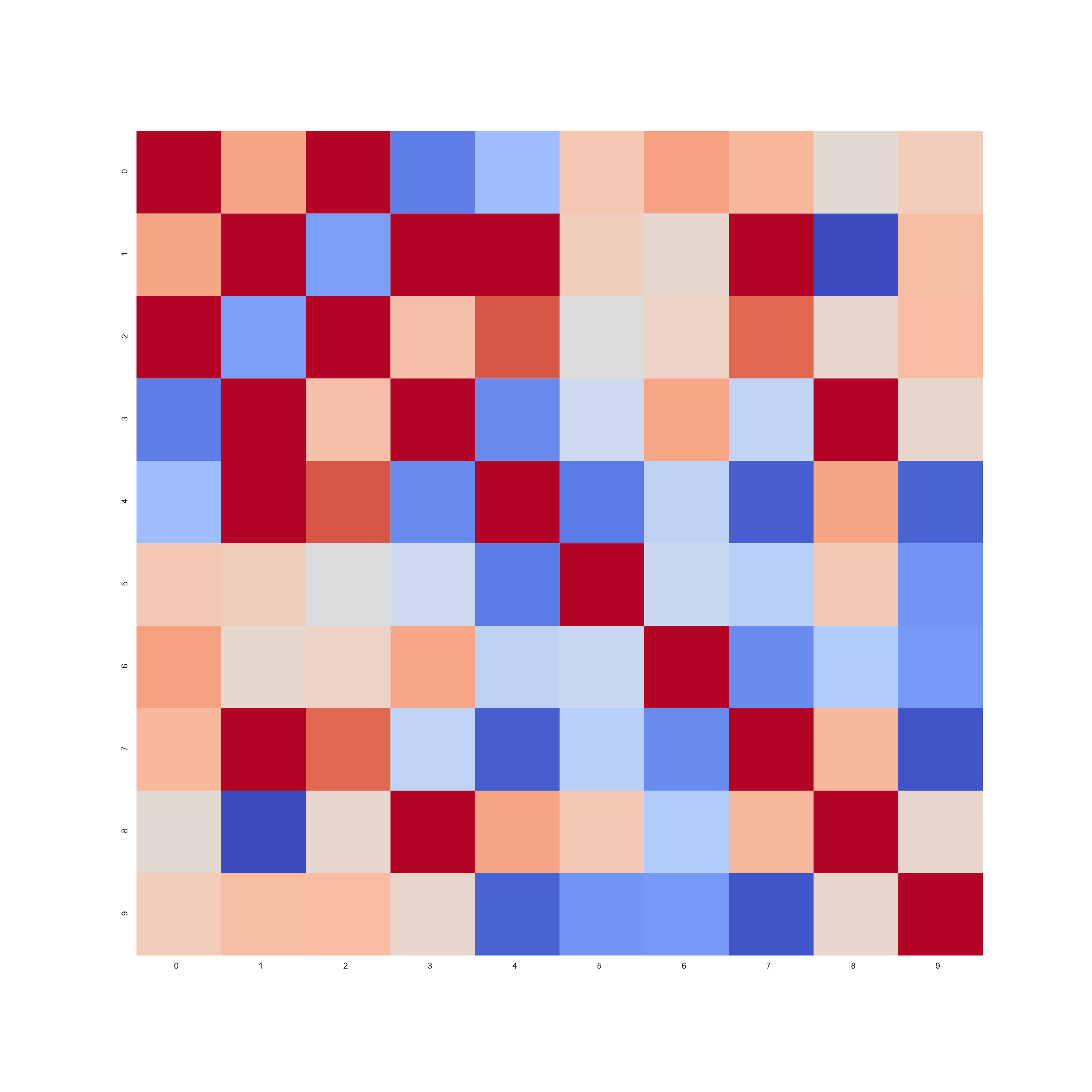}}
    }&
    \adjustbox{valign=c, scale=0.8}{
    \subfloat{\includegraphics[width=0.32\linewidth]{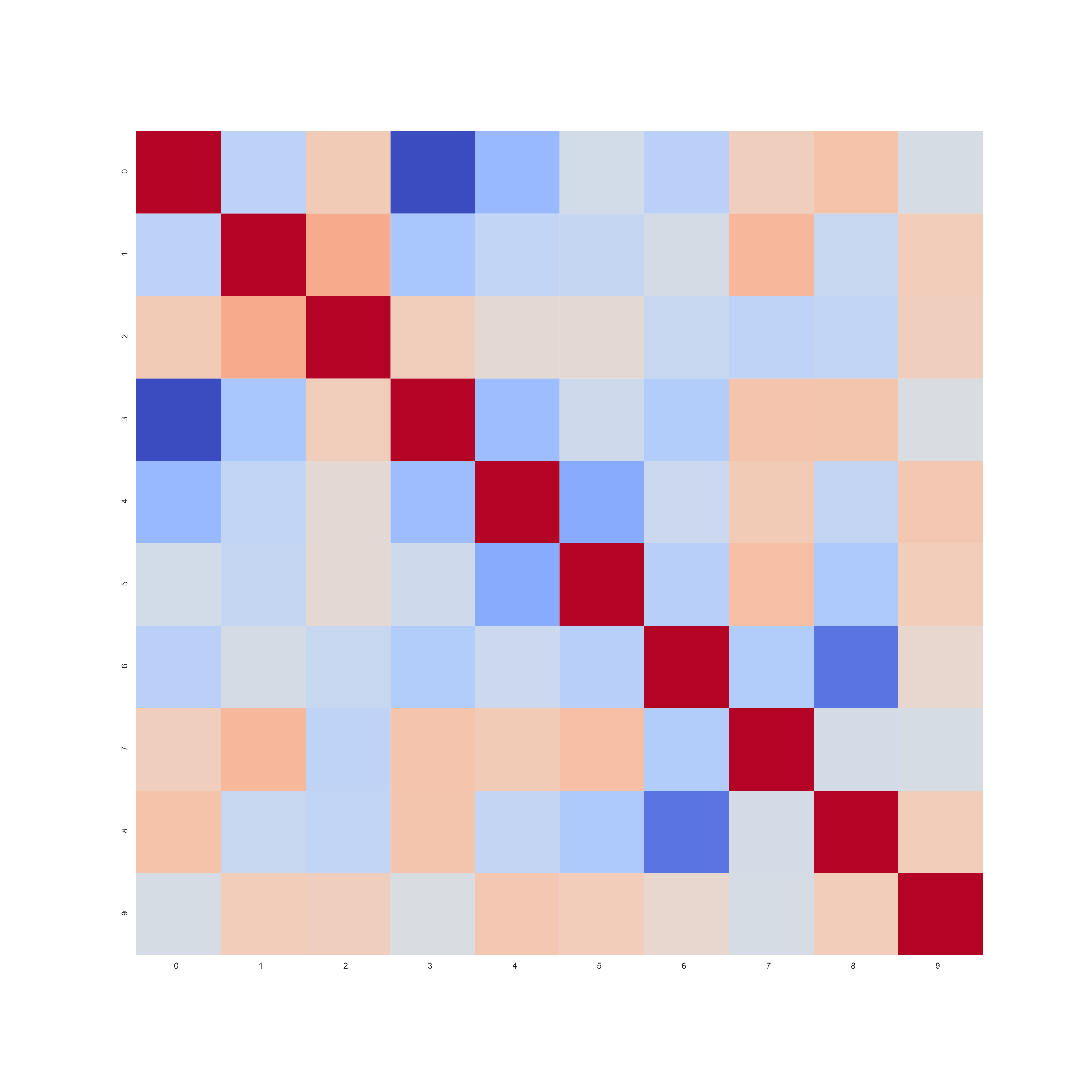}}
    }&
    \adjustbox{valign=c, scale=0.8}{
    \subfloat{\includegraphics[width=0.32\linewidth]{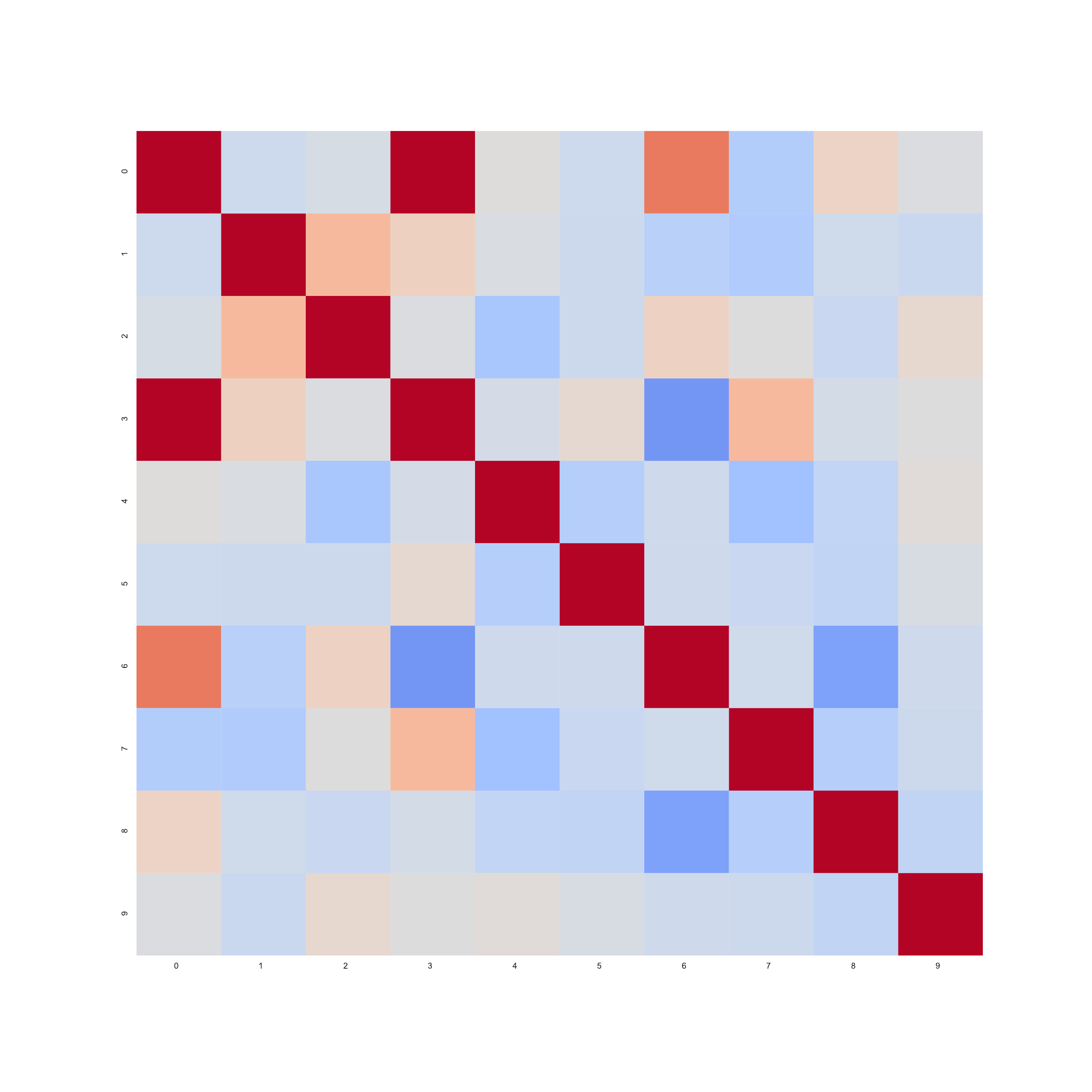}}
    }
    \\
    \footnotesize EM res. &
    \adjustbox{valign=c, scale=0.8}{
    \subfloat{\includegraphics[width=0.32\linewidth]{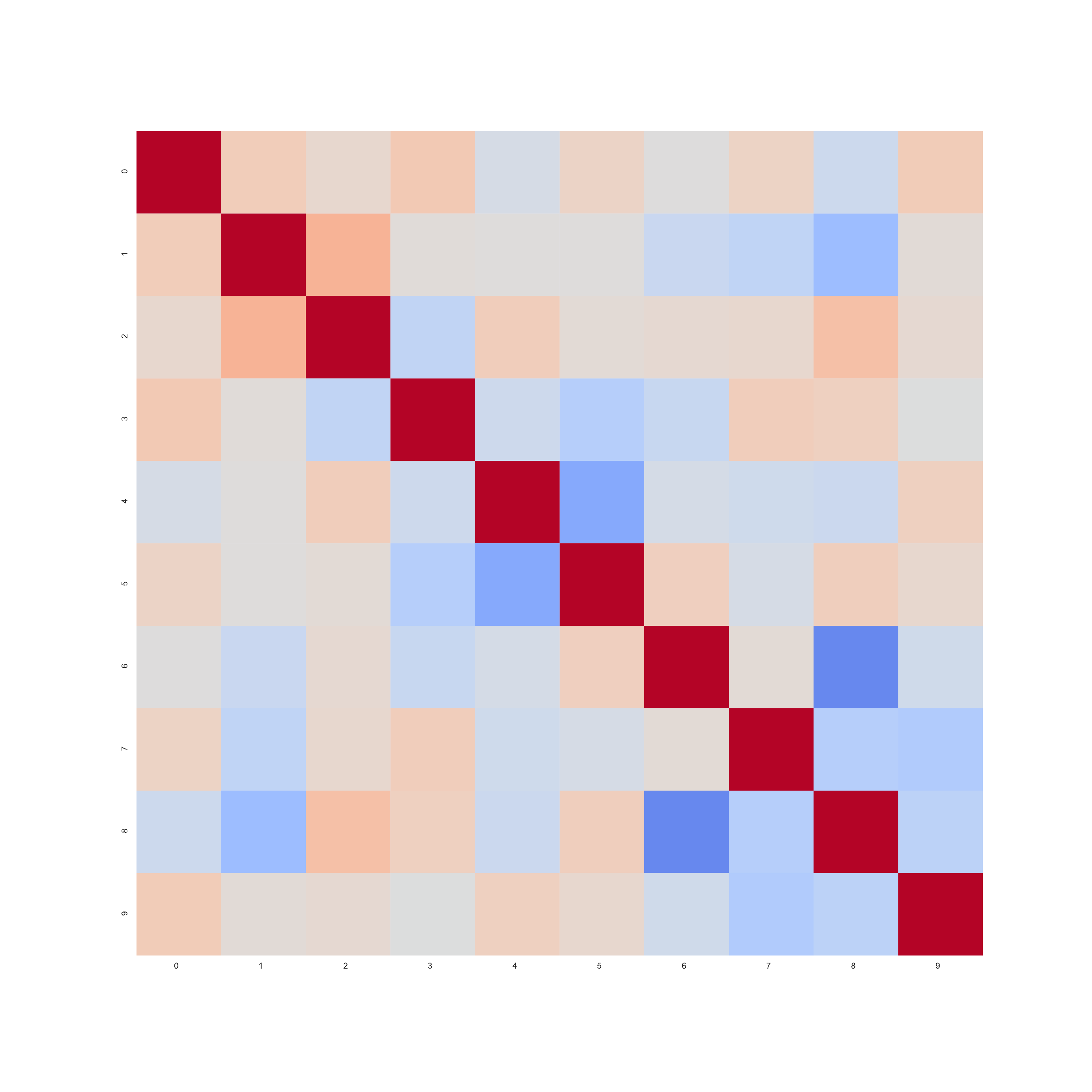}}
    }&
    \adjustbox{valign=c, scale=0.8}{
    \subfloat{\includegraphics[width=0.32\linewidth]{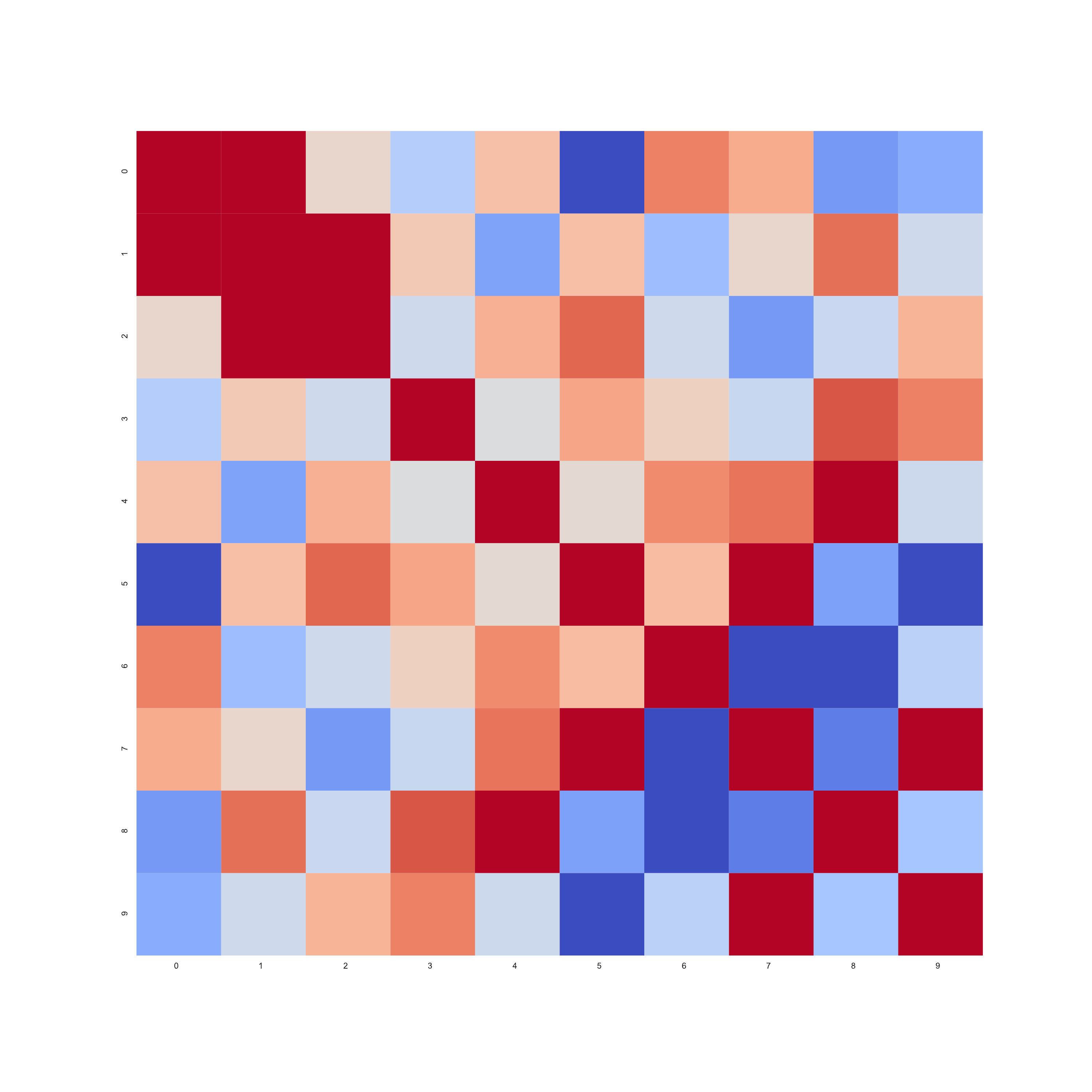}}
    }&
    \adjustbox{valign=c, scale=0.8}{
    \subfloat{\includegraphics[width=0.32\linewidth]{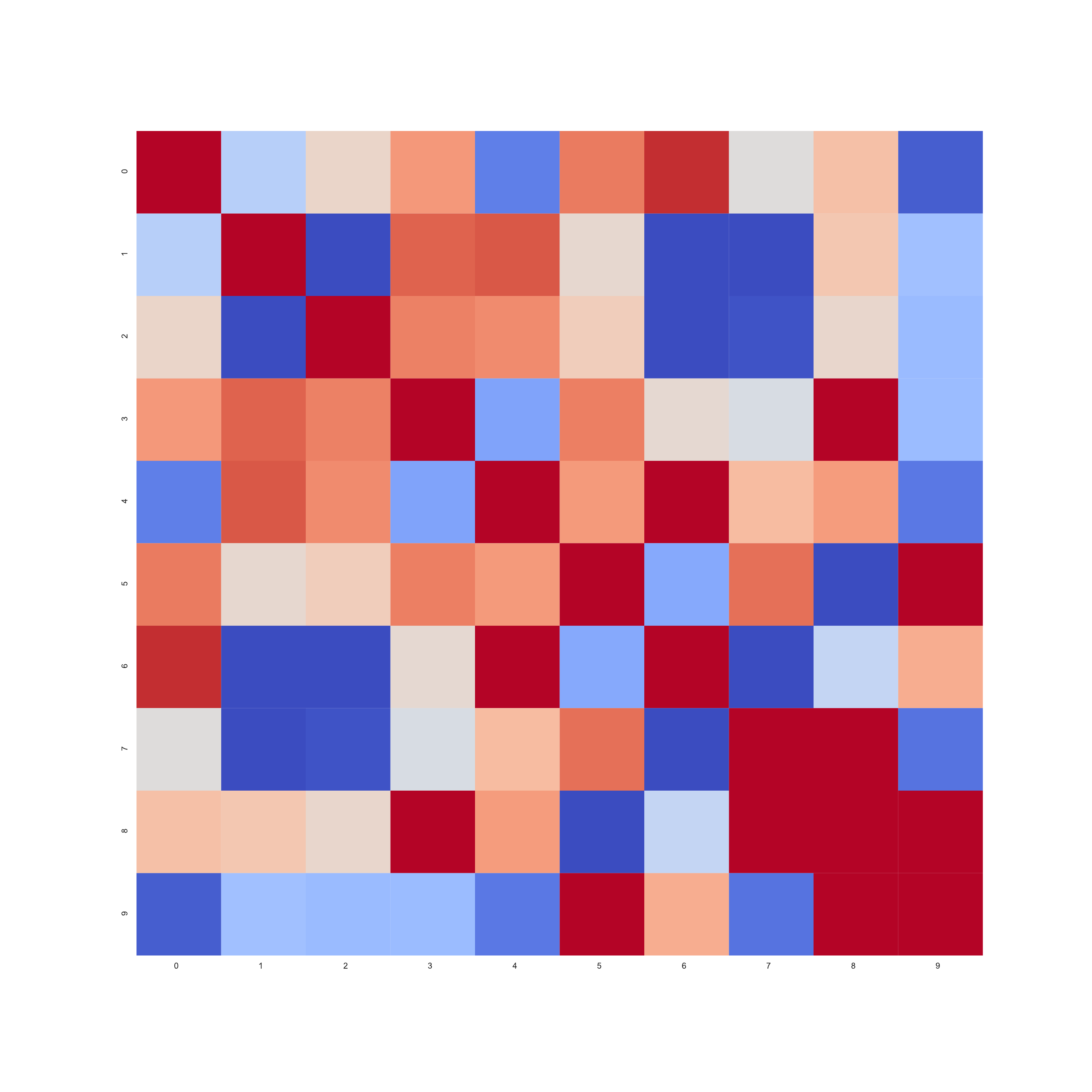}}
    }
    \\
    \footnotesize C-EM  &
    \adjustbox{valign=c, scale=0.8}{
    \subfloat{\includegraphics[width=0.32\linewidth]{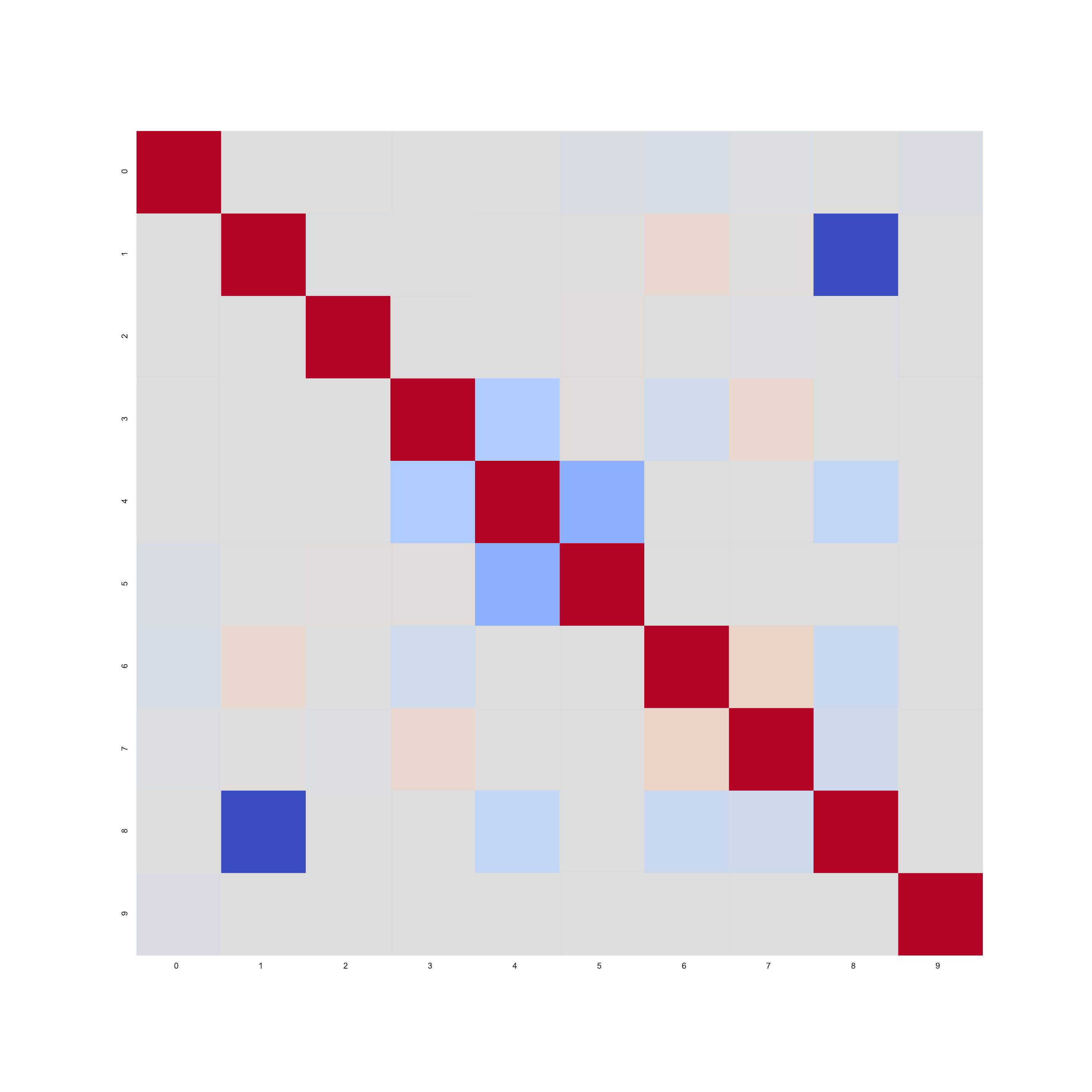}}
    }&
    \adjustbox{valign=c, scale=0.8}{
    \subfloat{\includegraphics[width=0.32\linewidth]{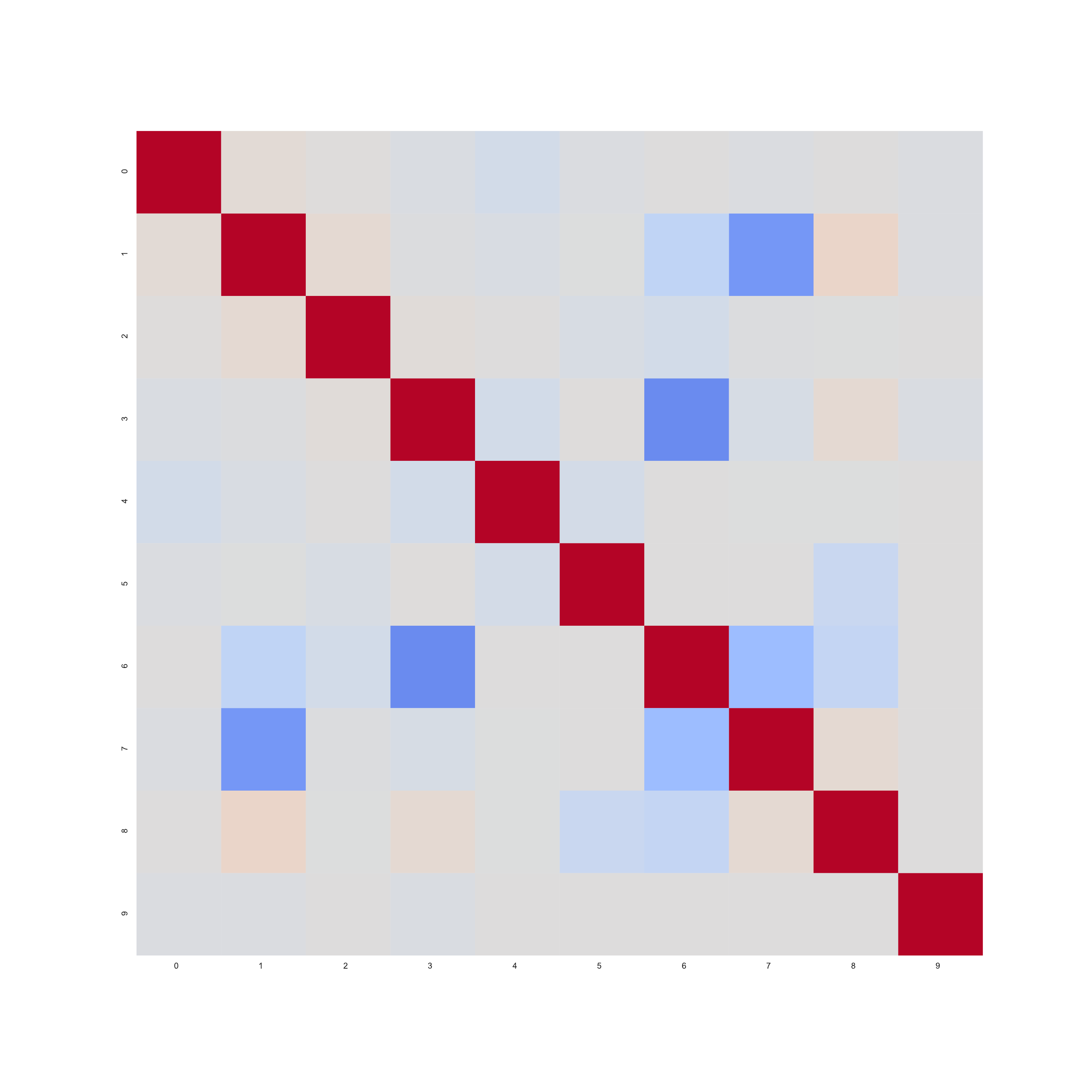}}
    }&
    \adjustbox{valign=c, scale=0.8}{
    \subfloat{\includegraphics[width=0.32\linewidth]{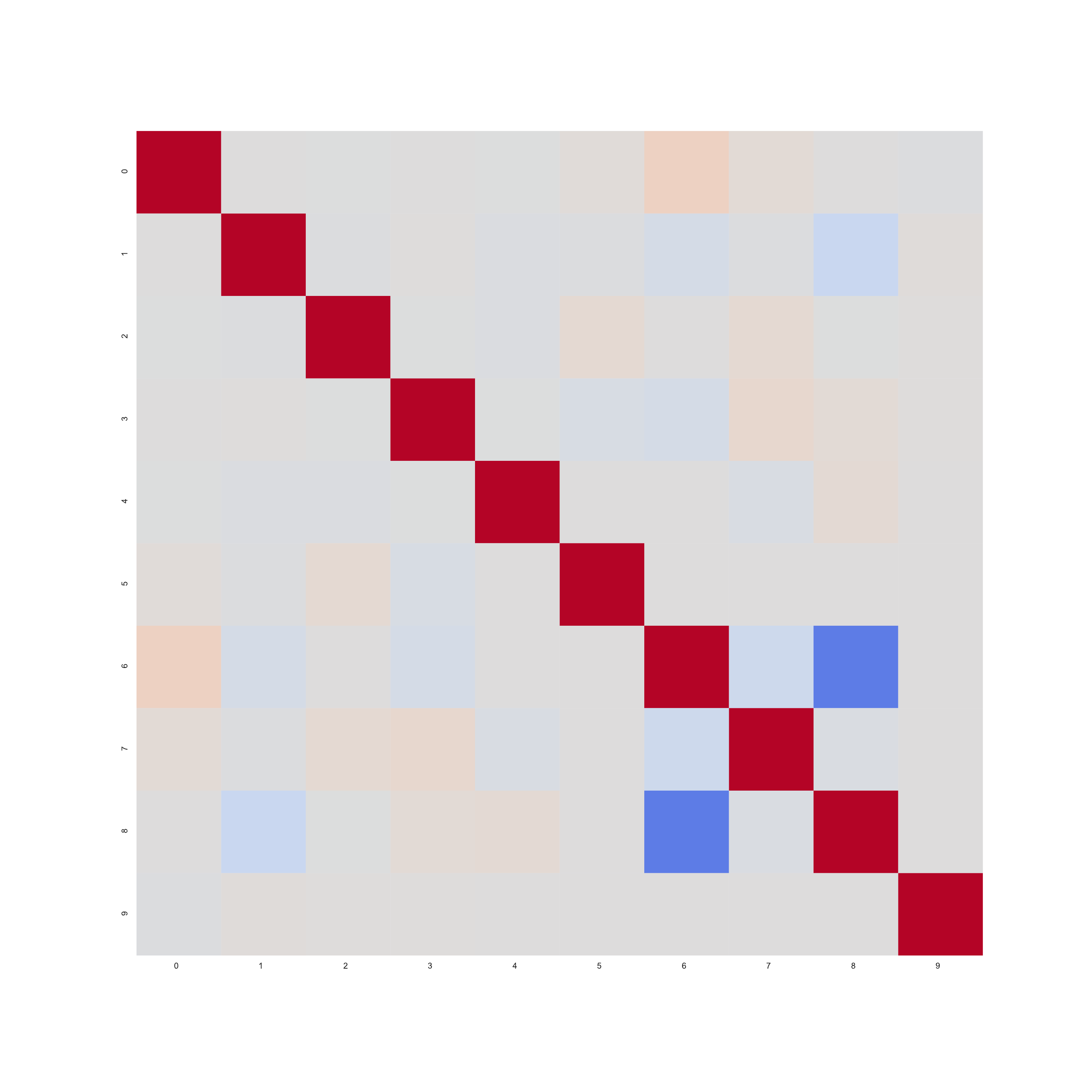}}
    }
    \\
    \footnotesize True &
    \adjustbox{valign=c, scale=0.8}{
    \subfloat{\includegraphics[width=0.32\linewidth]{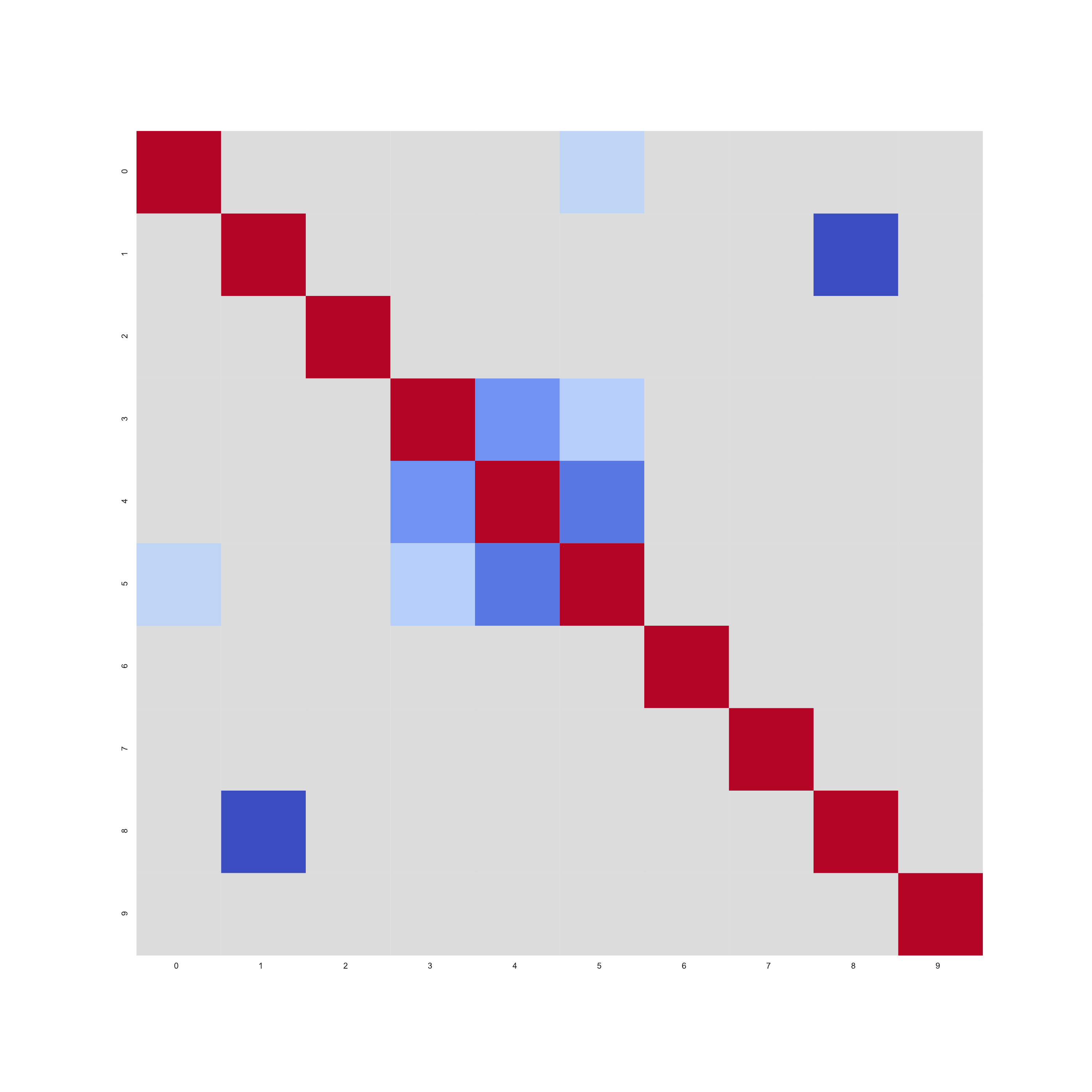}}
    }&
    \adjustbox{valign=c, scale=0.8}{
    \subfloat{\includegraphics[width=0.32\linewidth]{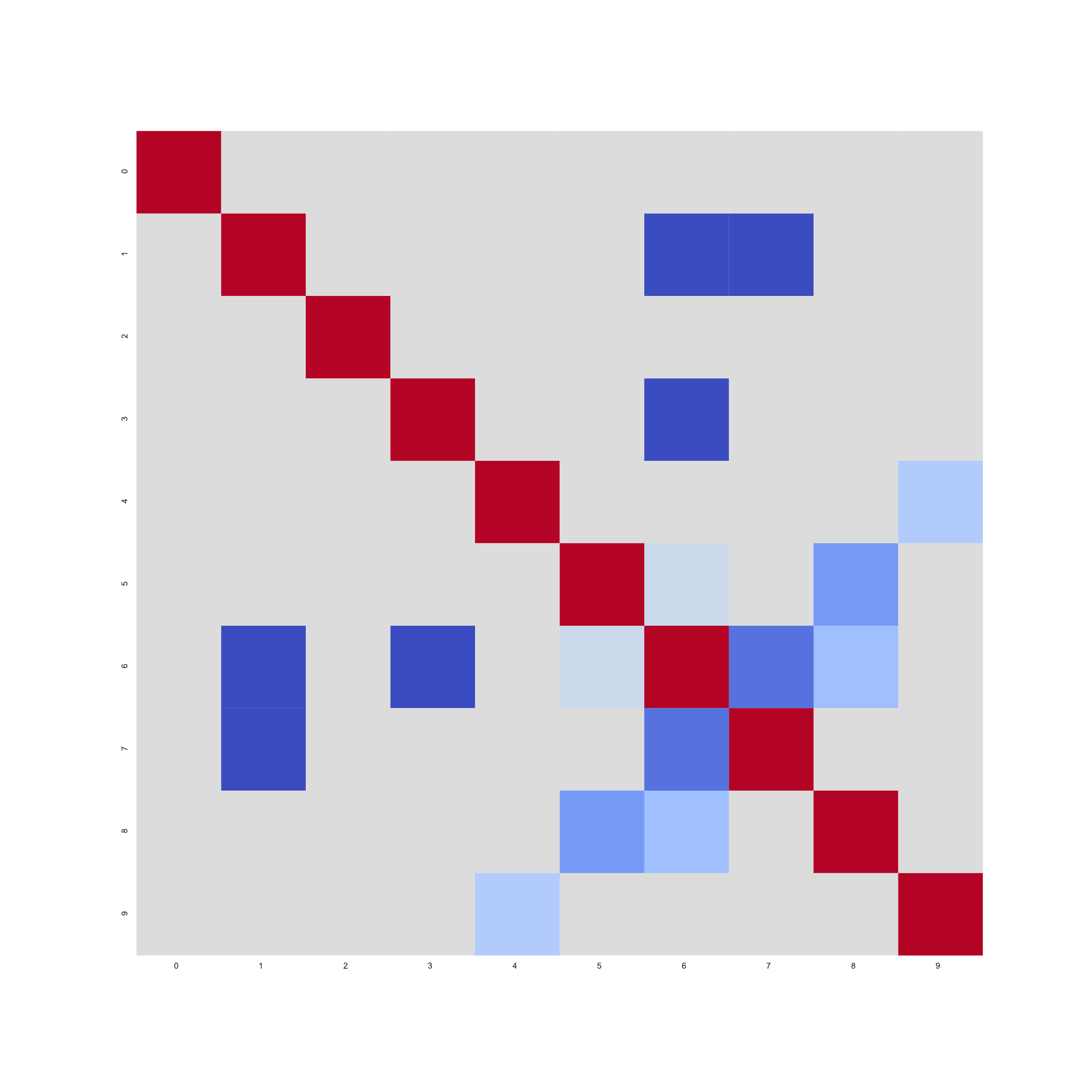}}
    }&
    \adjustbox{valign=c, scale=0.8}{
    \subfloat{\includegraphics[width=0.32\linewidth]{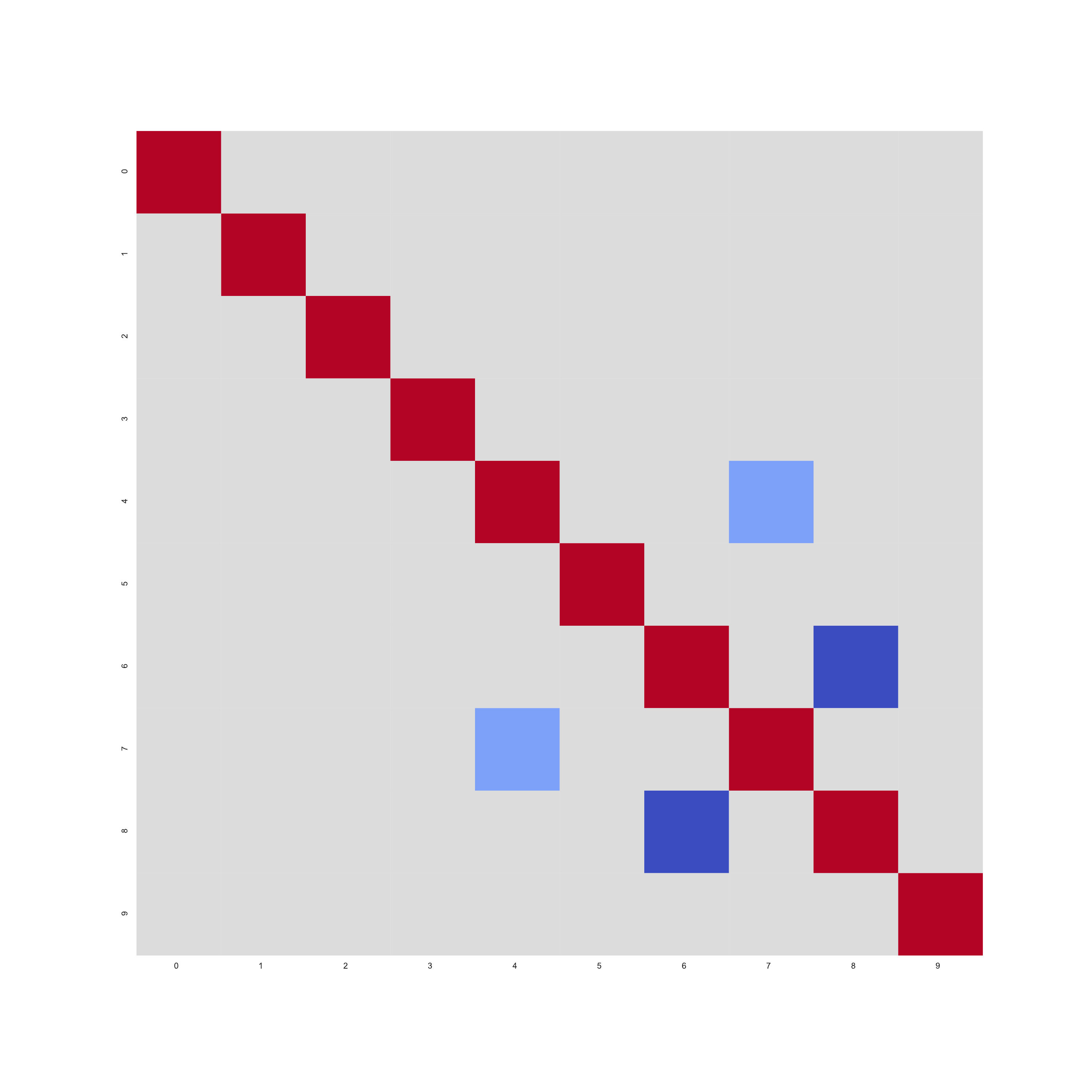}}
    }
    \end{tabular}
    }
    & 
    \hspace{-0.8cm}
    \adjustbox{valign=c, scale=0.8}{
    \subfloat{\includegraphics[width=0.1\linewidth]{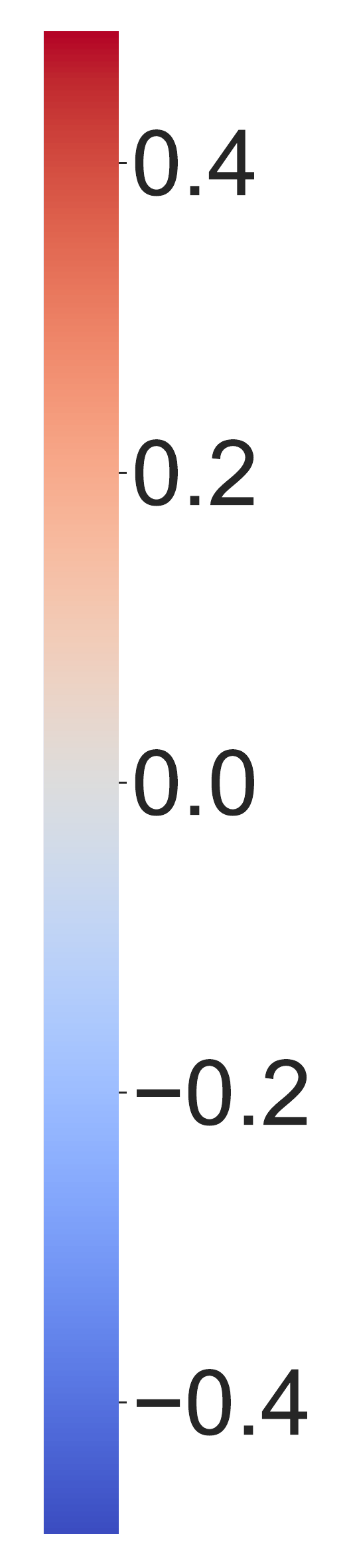}}
    }
    \end{tabular}
\caption{Comparison between the several estimated and the true conditional correlation matrices for each sub-population. The three columns of figures correspond to the three sub-populations. The first two rows of figures are the matrices estimated by the two Mixture of GGM methods, with and without residualisation with the co-features. The third row of figures correspond to the matrices estimated by the Mixture of CGGM. The final row displays the real conditional correlation matrices. Unlike the two GGM-based methods, the Mixture of CGGM recovers correct edges with very few False Positives.}
\label{fig:HD_illustrative_Sigma_cond}
\end{figure}

\begin{figure}[tbhp]
\centering
    \begin{tabular}{cl}
    \begin{tabular}{llll}
    & \footnotesize class 1 & \footnotesize class 2 & \footnotesize class 3\\
    \footnotesize C-EM  &
    \hspace{-0.62cm}
    \adjustbox{valign=c, scale=0.8}{
    \subfloat{\includegraphics[width=0.32\linewidth]{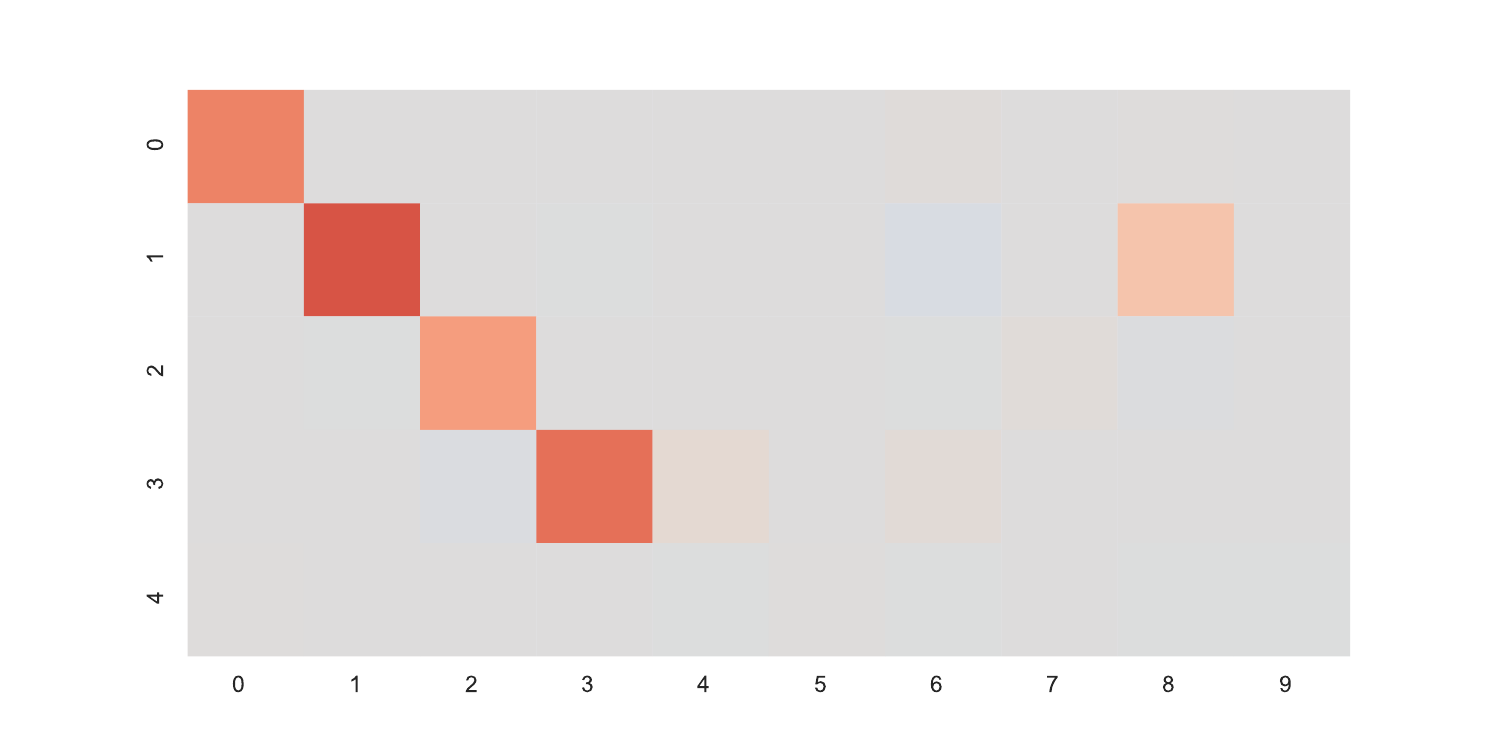}}
    }&
    \hspace{-0.62cm}
    \adjustbox{valign=c, scale=0.8}{
    \subfloat{\includegraphics[width=0.32\linewidth]{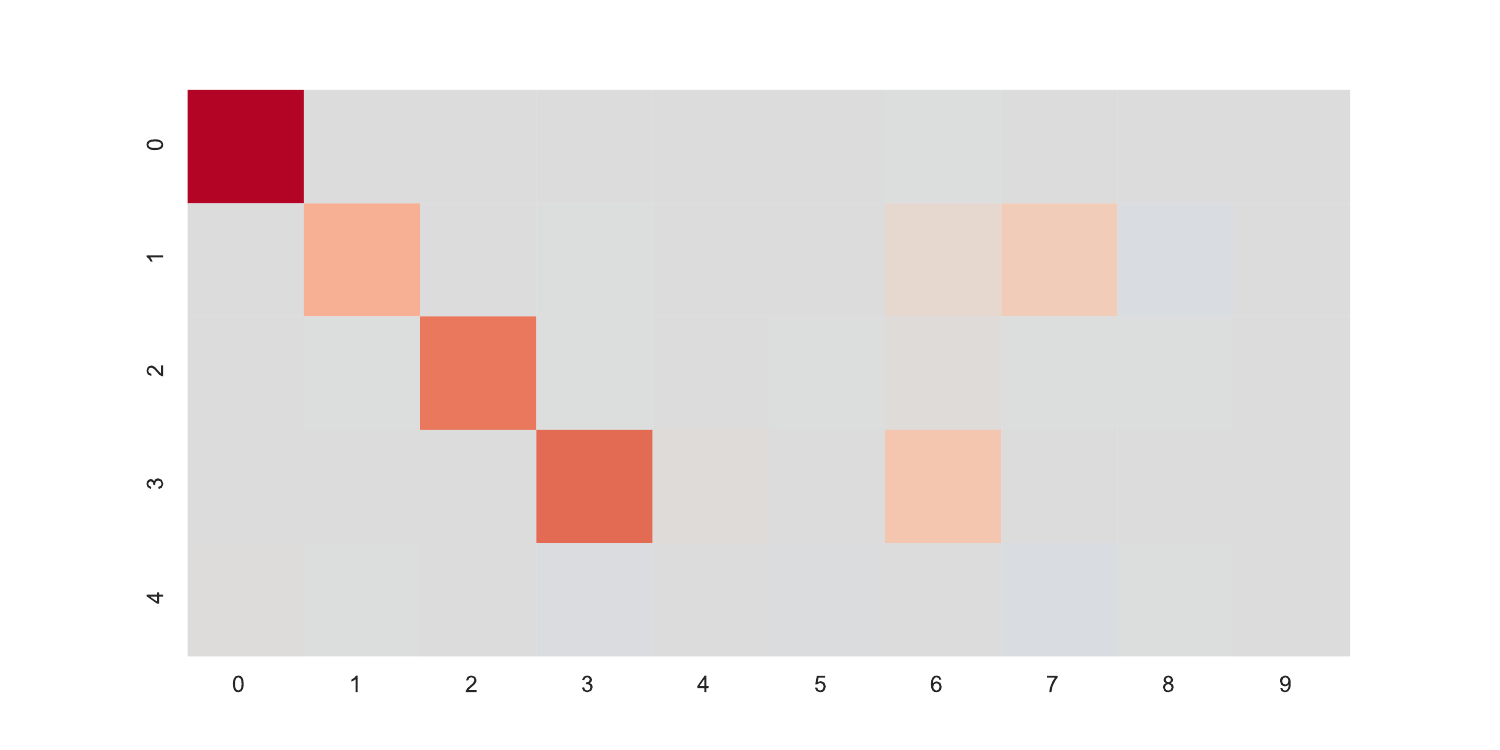}}
    }&
    \hspace{-0.62cm}
    \adjustbox{valign=c, scale=0.8}{
    \subfloat{\includegraphics[width=0.32\linewidth]{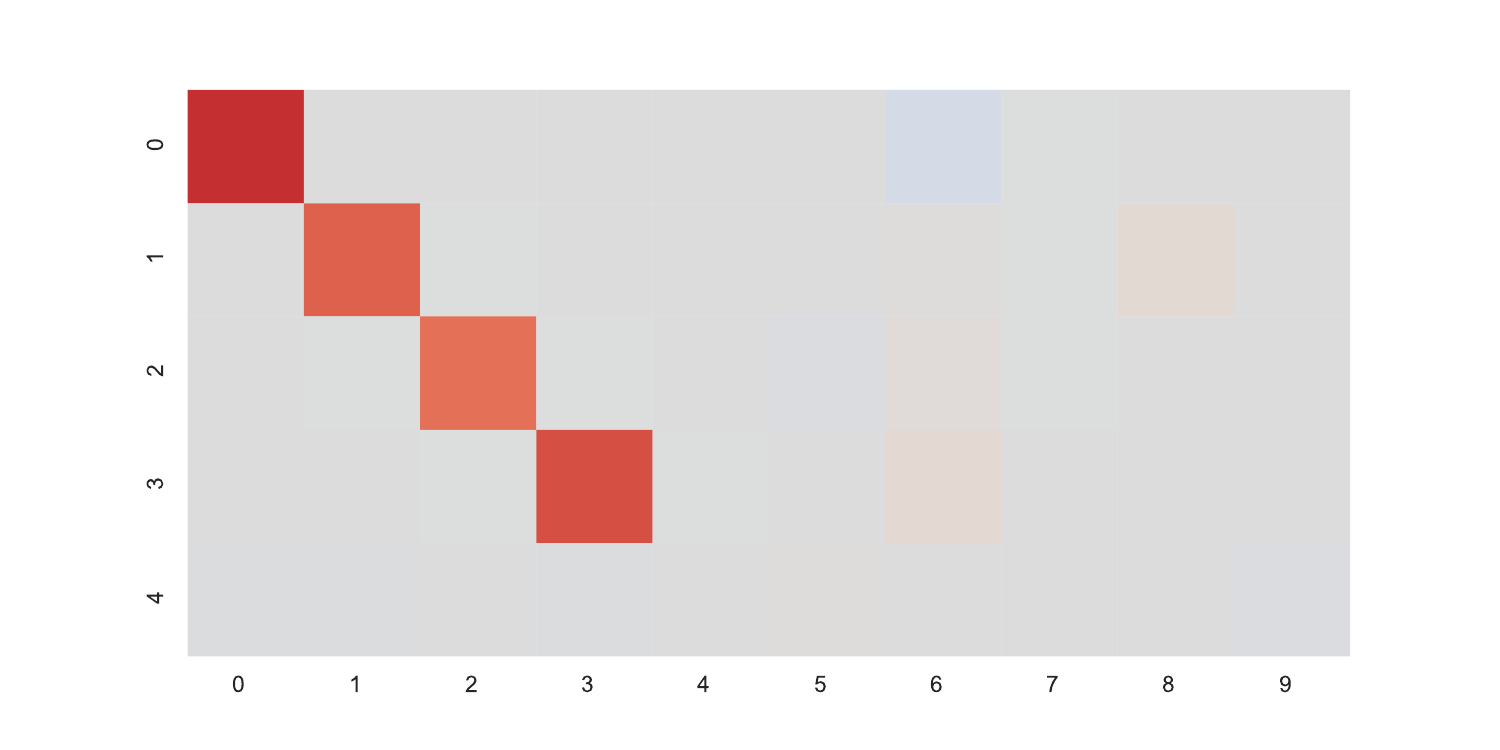}}
    }
    \\
    \footnotesize True &
    \hspace{-0.62cm}
    \adjustbox{valign=c, scale=0.8}{
    \subfloat{\includegraphics[width=0.32\linewidth]{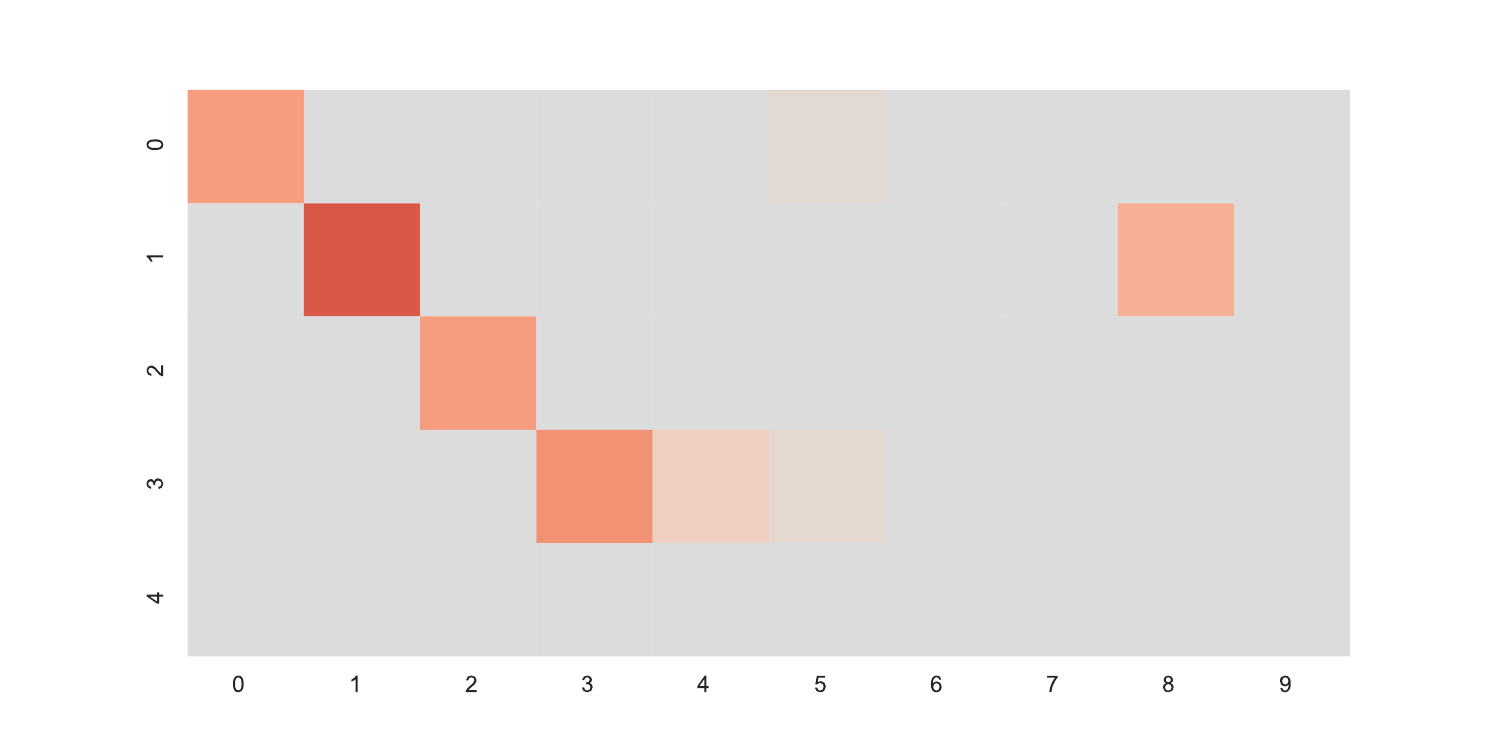}}
    }&
    \hspace{-0.62cm}
    \adjustbox{valign=c, scale=0.8}{
    \subfloat{\includegraphics[width=0.32\linewidth]{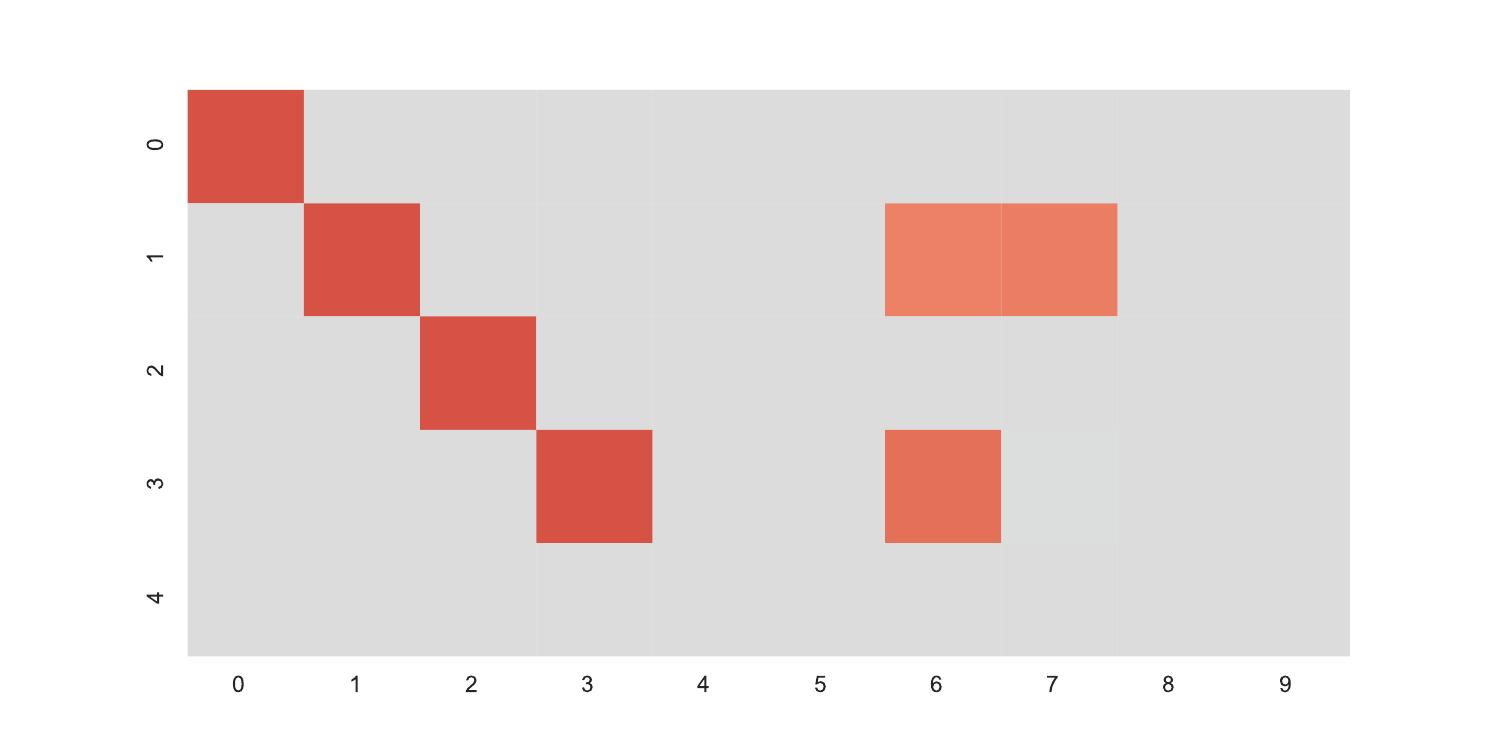}}
    }&
    \hspace{-0.62cm}
    \adjustbox{valign=c, scale=0.8}{
    \subfloat{\includegraphics[width=0.32\linewidth]{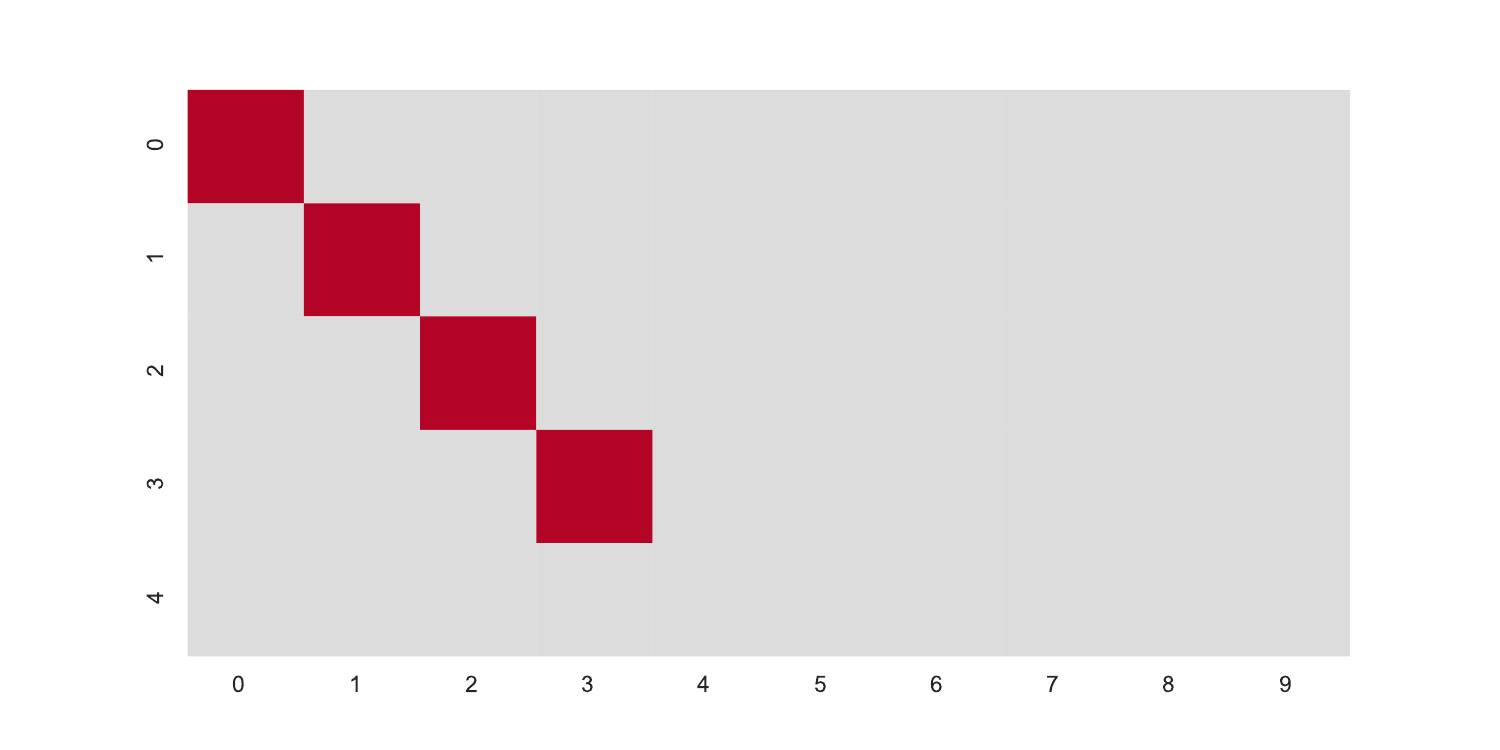}}
    }
    \end{tabular}
    &
    \hspace{-0.8cm}
    \adjustbox{valign=c, scale=0.4}{
    \subfloat{\includegraphics[width=0.15\linewidth]{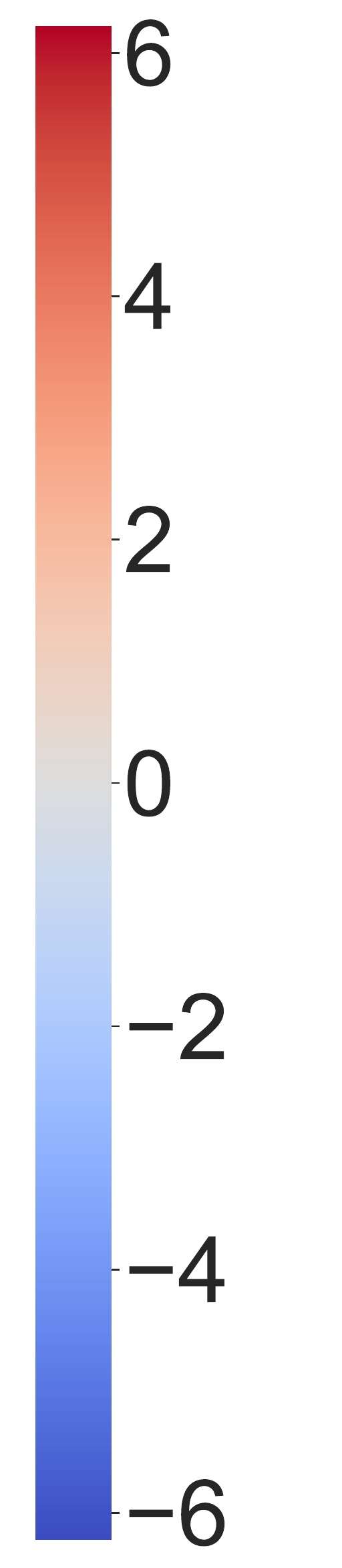}}
    }
    \end{tabular}
\caption{Reconstruction of the $\Theta_k$ by the EM on the Mixture of CGGM. The three columns of figures correspond to the three sub-populations. Almost all the edges are right, with no False Positive and almost no False Negative. Moreover, the intensities are also mostly correct.}
\label{fig:HD_illustrative_Theta}
\end{figure}

In this section, we perform a quantitative analysis of the algorithms in a higher dimension framework, where the matrix parameters $\Lambda$ and $\Theta$ are more naturally interpreted as sparse networks. We confirm that the Mixture of Conditional Gaussian Graphical Models is better suited to take into account the heterogeneous effects of co-features on the graph.\\
For this experiment, the observed data follows a mixture model with $K=3$ classes. Each class $k$ has the probability weight $\pi_k=\frac{1}{3}$. An observation $(Y, X) \in \R^p \times \R^q$ belonging to the class $k$ is described by the distribution: $Y | X \sim \N{- \Lambda_k^{-1} \Theta_k^T X, \Lambda_k^{-1}}$. In this example, the dimensions of $Y$ and $X$ are respectively $p=10$ and $q=5$. $X$ contains four independent Rademacher variables ($=+1$ or $-1$) and a constant variable always equal to 1. The precision matrix $\Lambda_k \in \R^{p \times p}$ and the transition matrices $\Theta_k \in \R^{q \times p}$ are both sparse. For $\Lambda$, over the $45$ off-diagonal coefficients, $\Lambda_1$ has $7$ non-zero coefficients, $\Lambda_2$ has $16$, and $\Lambda_3$ has $2$. Resulting in an average $18.5 \%$ rate of sparsity among the off-diagonal coefficients. The sparsity pattern of the $\Lambda_k$ are represented on \figurename~\ref{fig:HD_illustrative_Sigma_cond}, within the row of figures labelled "True". For $\Theta$, over $50$ possible coefficients, $\Theta_1$ has $11$ non-zero ones, $\Theta_2$ has $14$, and $\Theta_3$ has $4$. Resulting in an average sparsity rate of $19.3\%$. Additionally, there is a repeating pattern in the $\Theta_k$: the $4$ random components of $X$ each affect (at least) the corresponding component among the first $4$ of $Y$. That is to say that $\forall k, \forall i \leq 4, (\Theta_k)^{i,i} \neq 0$. The real $\Theta_k$ are represented on \figurename~\ref{fig:HD_illustrative_Theta}, within the row of figures labelled "True". \\
We run 20 simulations. A simulation consists of $n=100$ generated data points. On these data points, we run the compared methods, all initialised with the same random parameters. For all simulations, we make $5$ of these runs, each with a different random initialisation. We compared the same three algorithms as in section \ref{sect:experiment_2D}: the EM for the Mixture of GGM, the EM for the Mixture of GGM with average effect of $X$ subtracted, and the EM applied to the Mixture of CGGM. \\
We follow three metrics to assess the method's success in terms of cluster recovery. The classification error (both soft and hard labels versions) and an "ABC-like" metric. The "ABC-like" metric is meant to assess how well each of the estimated solutions is able to replicate the observed data. Since each solution is the parameter of a probability distribution, at the end of each EM, we generate new data following this proposed distribution. Then, for each synthetic data point, we compute the distance to the closest neighbour among the real data points. These minimal distances constitute our "ABC-like" metric. Finally, we also compute the execution time of each EM, knowing that they all have the same stopping criteria. We represent on \figurename~\ref{fig:HD_all_metrics} the empirical distribution of these four metrics and we quantify with Table \ref{tab:HD_all_metrics} the key statistics (mean, standard deviation, median) that characterise them. With $K=3$ and balanced classes, a uniform random classifier would guess the wrong label $66.7\%$ of the time. We observe that the two Mixture of GGM method are dangerously close to this threshold, with more than $50\%$ hard misclassification. The EM on the Mixture of CGGM (C-EM) on the other hand, achieves a much better classification with less than $15\%$ hard misclassifcation. This demonstrate that, even when faced with a more complex situation, in higher dimension, the Mixture of CGGM is better suited to correct for the effect of the co-features and discover the right clusters of data points. This also underlines once more the importance of allowing different values of the effect of $X$ on $Y$ for each class. Indeed, the residual Mixture of GGM - which took into account the average effect of $X$ on $Y$ over the entire population - was unable to achieve better performances than the EM that did not even use the co-features $X$. In terms of reconstruction of the observed data by the estimated model (ABC-like metric), we see that the synthetic data points generated from the estimated Mixture of CGGM model have closer nearest neighbours than the data points generated by the other estimated models. In addition to all these observations, the C-EM is also faster than the other two methods, reaching the convergence threshold faster.
\begin{table}[tbhp]
    \centering
    {\footnotesize
    \caption{\footnotesize Average, standard deviation and median (below) of the four followed performance metrics over the $30\times 5$ simulations. The best values are in \textbf{bold}. We can see that the classification performances with the Mixture of CGGM are much better than the two methods with Mixtures of GGM, and with faster computation times.}
    \label{tab:HD_all_metrics}
    \begin{tabular}{l|ccc}
        \toprule\\
                         & EM GGM                  &    EM GGM resid.         &     EM CGGM        \\
        \midrule \\
        soft misclassif. & 0.56 \textbf{(0.03)} & 0.51 \textbf{(0.03)}  & \textbf{0.17} (0.05) \\
                         &    0.57              &    0.51               &    \textbf{0.16}     \\ 
        \midrule \\
        hard misclassif. & 0.55 (0.04)          & 0.50 \textbf{(0.03)}  & \textbf{0.14} (0.06) \\
                         &    0.56              &    0.49               &     \textbf{0.13}    \\
        \midrule \\
        ABC-like metric  & 5.57 \textbf{(0.09)} & 4.64 (0.22)           & \textbf{4.13} (0.14) \\
                         &   5.58               &  4.64                 &    \textbf{4.14}     \\
        \midrule \\
        runtimes         & 115 \textbf{(61)}    & 253 (137)             & \textbf{58} (91)     \\
                         &   93                 &   256                 &    \textbf{16}       \\
        \bottomrule
        \end{tabular}
        }
\end{table}

\noindent
In addition to the cluster recovery, we can also assess the parameter reconstruction of each method. Since the three clustering methods estimate different parametric models over the data, they do not actually try to estimate the same parameters. Regardless, all the methods still estimate a certain precision matrix $\Lambda_k$ (conditional or not on the $X$ depending on the model) of each sub-population that they identify. Hence, we also record the differences between the true $\Lambda_k$ and their estimations by the methods. The two metrics followed are the Kullback–Leibler (KL) divergence \cite{kullback1951information} between the gaussian distribution $f_{\Lambda_k} \sim \mathcal{N}(0, \Lambda_k^{-1}))$ and $f_{\widehat{\Lambda}_k} \sim \mathcal{N}(0, \widehat{\Lambda}_k^{-1}))$, and the $l_2$ difference given by the Froebenius norm: $\norm{\Lambda_k-\widehat{\Lambda}_k}_F^2$. In Table \ref{tab:synthetic_parameter_reconstruction}, we can check that, according to these metrics, the $\widehat{\Lambda}_k$ estimated by with the Mixture of CGGM are indeed a much better fit for the real $\Lambda_k$ than the estimated matrices from the other models. This is expected, since the real $\Lambda_k$ actually correspond to the CGGM model.\\
To illustrate the different level of success concerning the conditional correlation graph recovery, we display on \figurename~\ref{fig:HD_illustrative_Sigma_cond} the conditional correlation matrix (i.e. the conditional correlation graph with weighted edges) estimated by each method. The three columns of figures correspond to the three sub-populations. The first two rows of figures are the matrices estimated by the two Mixtures of GGM methods, with and without residualisation with the co-features. The third row of figures correspond to the matrices estimated by the Mixture of CGGM. The final row displays the real conditional correlation matrices. We observe that the two Mixture of GGM recover way too many edges, with no particular fit with the real matrix. By contrast, the matrices from the CGGM Mixture exhibit the proper edge patterns, with few False Positive and False Negative. 
This is not an easy feat to achieve, since the method was run from a random initialisation on a totally unsupervised dataset, with heavily translated data points all over the 10 dimensional space. Moreover, the matrices in \figurename~\ref{fig:HD_illustrative_Sigma_cond} all result from the inversion of the empirical covariance matrix, which is neither a very geometrical nor a very stable operation.\\ 
In \figurename~\ref{fig:HD_illustrative_Theta}, we represent the regression parameter $\widehat{\Theta}_k$ estimated by with the Mixture of CGGM alongside the real $\Theta_k$. Once again, we see that the sparsity pattern is very well identified, with no False Positive. Moreover, in this case, there are also almost no False Negative, and all the edge intensities are correct. This is not a surprise. Indeed, the parameter $\Theta$ plays a huge role in the correct classification of the data, since it serves to define the expected position of each data point in the feature space (playing the part of the "average" parameter in Mixtures of GGM). Hence, a good estimation of $\Theta$ is mandatory to reach a good classification. Since the EM with Mixture of CGGM achieved very good classification results, it was expected that $\Theta$ would be well estimated.
\begin{table}[tbhp]
    \begin{center}
    {\footnotesize
    \caption{\footnotesize Average and standard deviation of the metrics describing the reconstruction of each precision matrix $\Lambda_k$. The matrices are consistently better reconstructed with the mixture of CGGM.} 
    \label{tab:synthetic_parameter_reconstruction}
    \begin{tabular}{lc|ccc}
    \toprule
    metric &     class   & EM GGM & EM res. GGM & EM CGGM\\
    \midrule\\
    \multirow{3}{*}{$KL(f_{\Lambda}, f_{\hat{\Lambda}})$}                  
    & 1  & 11.0 (3.0) & 7.5 (6.8)  & \textbf{0.8 (0.2)} \\
    & 2  & 10.3 (2.2) & 8.5 (5.0)  & \textbf{1.9 (0.3)}  \\
    & 3 &  13.6 (2.5) & 5.2 (2.3)  & \textbf{3.4 (1.1)}  \\
    \midrule\\
    \multirow{3}{*}{$\norm{\Lambda-\hat{\Lambda}}_F^2$}                  
    & 1  & 39.2 (48.4) & 44.2 (114)  & \textbf{2.2 (0.8)} \\
    & 2  & 15.1 (12.2) & 102 (73.9)  & \textbf{6.6 (0.9)}  \\
    & 3 &  14.2 (13.8) & 15.1 (25.7) & \textbf{5.8 (4.0)}  \\
    \bottomrule
    \end{tabular}
    }
    \end{center}
\end{table}

\subsection{Experiments on real data}
\begin{figure*}
    \centering
    \includegraphics[width=0.9\linewidth]{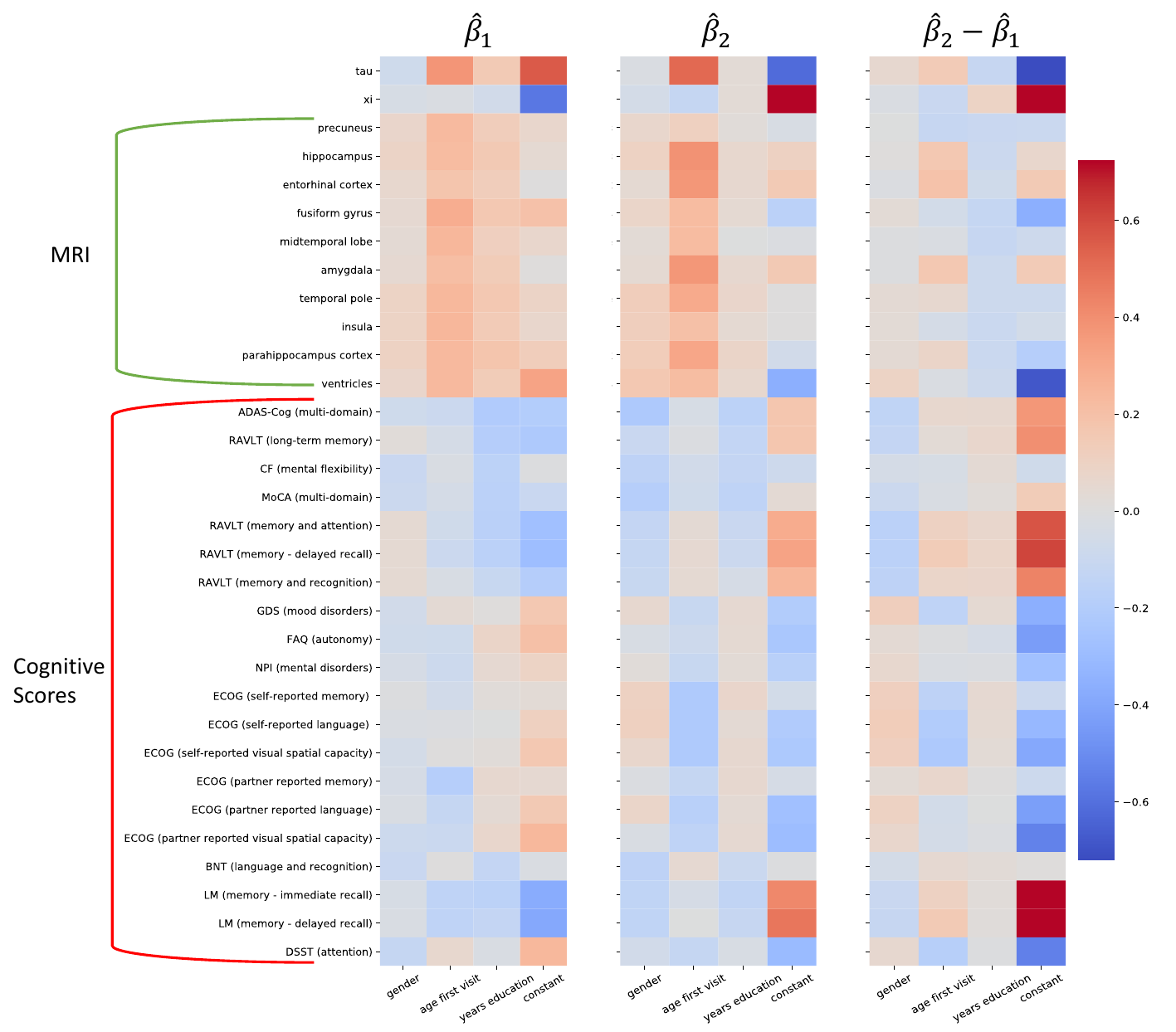}
    \caption{Average of $\widehat{\beta}_k := - \widehat{\Sigma}_k \widehat{\Theta}_k^T$ over 10 bootstrap sampling of the data. For each bootstrapped dataset, $\widehat{\Sigma}_k$ and $\widehat{\Theta}_k$ are estimated with the C-EM. We make 3 of these C-EM runs, each starting from a different KMeans initialisation of the labels. This results in $10 \times 3 = 30$ different estimates of $\widehat{\beta}_k$ to average over.  (Left) $\widehat{\beta}_1$, the cluster $k=1$ is always very correlated with the \textbf{Control patients} sub-population (less than 10\% deviation). (Middle) $\widehat{\beta}_2$, the cluster $k=2$ is likewise very correlated with the \textbf{AD patients}. In each figure, the last column is the constant coefficient. The largest inter-cluster differences are between the two constant terms. However there are some noticeable difference on the other regression coefficients as well.(Right) Average $\widehat{\beta}_2 - \widehat{\beta}_1$ over the 30 bootstrap runs of the C-EM. Here, the \textbf{differences in intensity} between AD ($k=2$) and Control ($k=1$) clusters are more explicit.}
    \label{fig:real_data_beta_CEM}
\end{figure*}

\begin{figure*}
    \centering
    \includegraphics[width=0.9\linewidth]{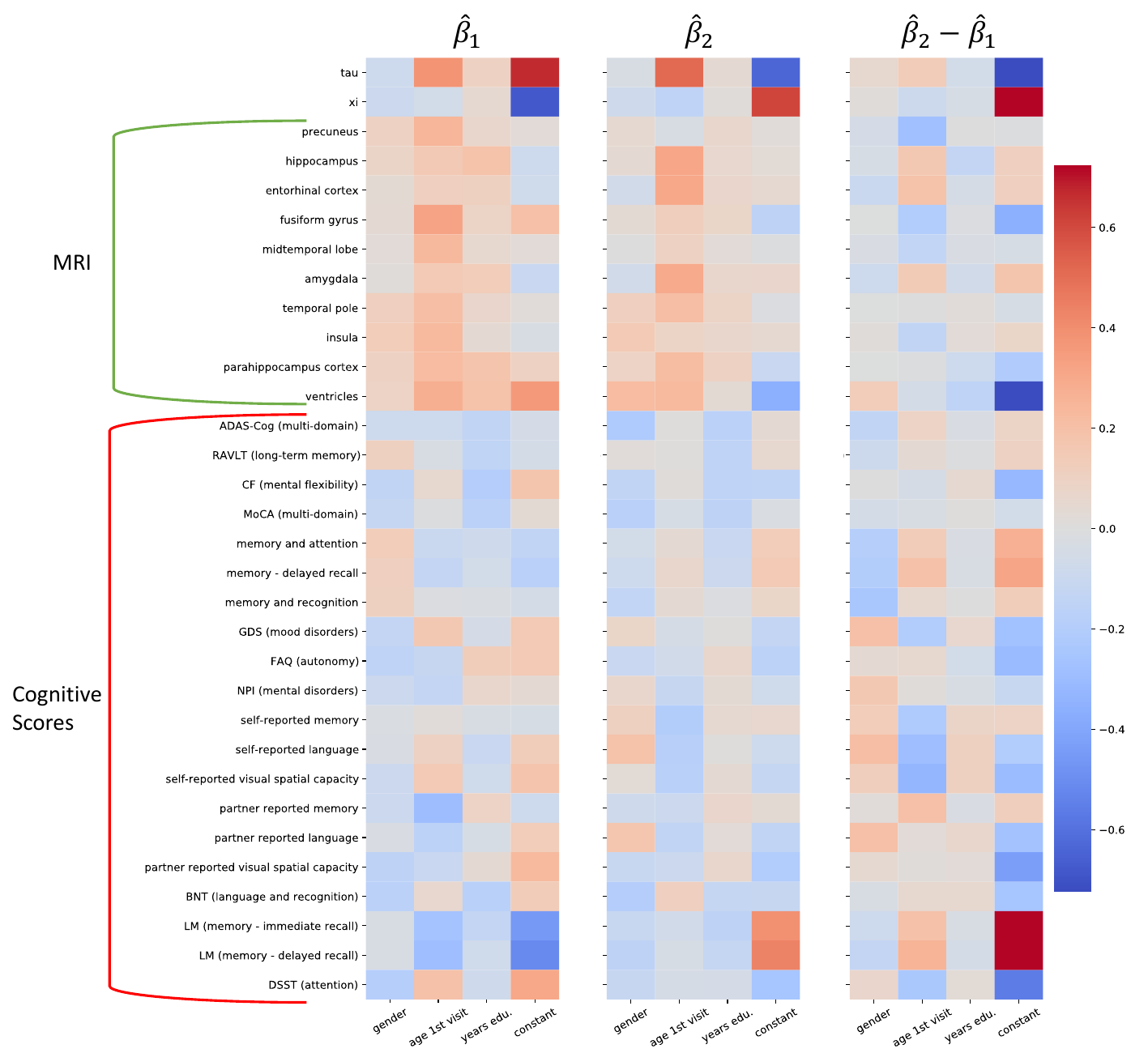}
    \caption{Oracle MLE $\widehat{\beta}_k := - \widehat{\Sigma}_k \widehat{\Theta}_k^T$ estimated with the knowledge of the real labels (supervised estimation). A comparison with \figurename~\ref{fig:real_data_beta_CEM} shows that the C-EM properly recoverd these structures. (Left) $\widehat{\beta}_1$, MLE of the \textbf{Control patients} sub-population. (Middle) $\widehat{\beta}_2$, MLE of the \textbf{AD patients}. (Right) $\widehat{\beta}_2 - \widehat{\beta}_1$, their \textbf{difference in intensity}.}
    \label{fig:real_data_beta_oracle}
\end{figure*}

\begin{figure*}[tbhp]
    \centering
    \subfloat{\includegraphics[width=0.49\linewidth]{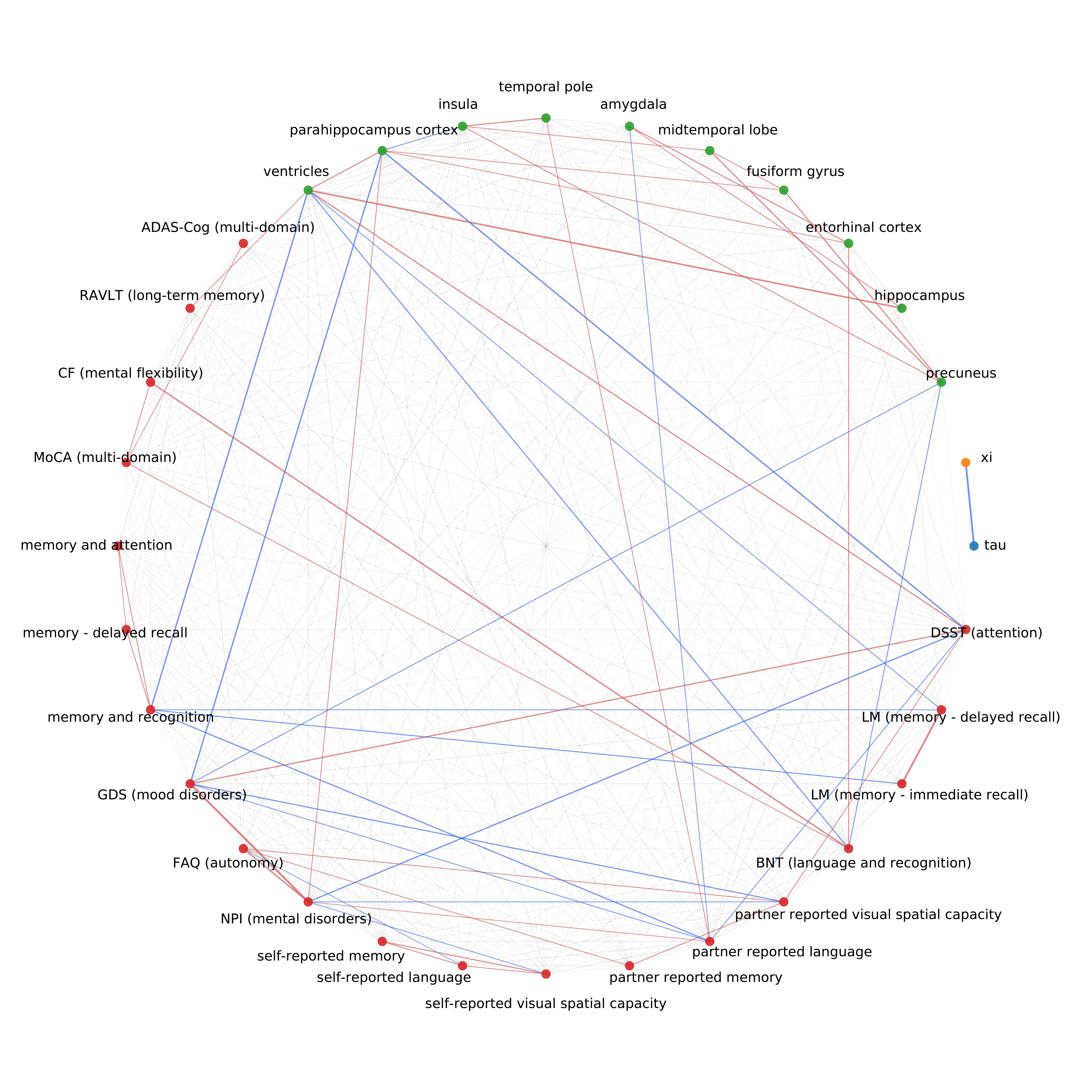}}
    \subfloat{\includegraphics[width=0.49\linewidth]{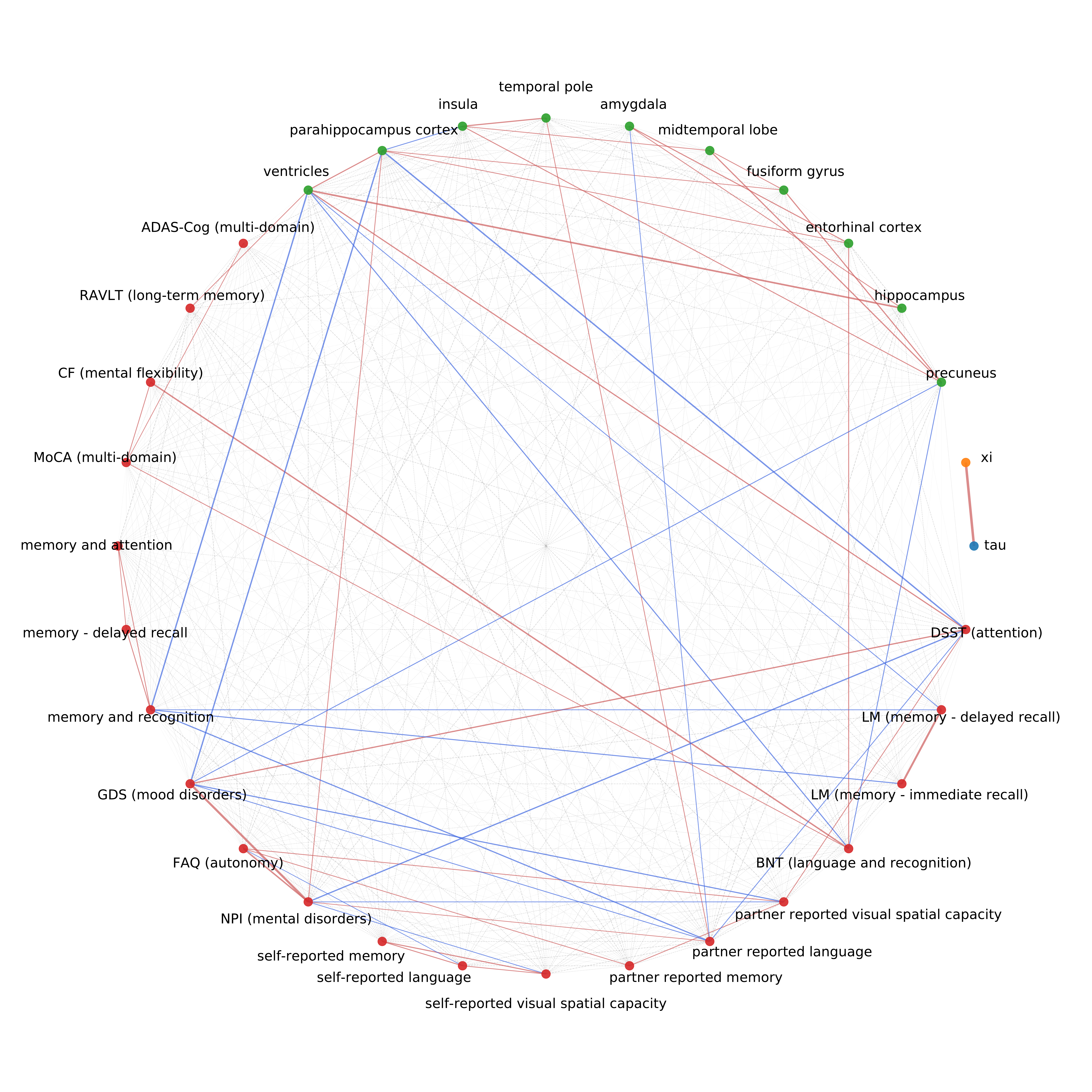}}
\caption{(Left) Conditional correlation graph of the estimated cluster most correlated with the \textbf{"Control" diagnosis}. (Right) Conditional correlation graph of the estimated cluster most correlated with the \textbf{"AD" diagnosis}. }
\label{fig:real_data_graph}
\end{figure*}

\begin{figure*}[tbhp]
    \centering
    \subfloat{\includegraphics[width=0.49\linewidth]{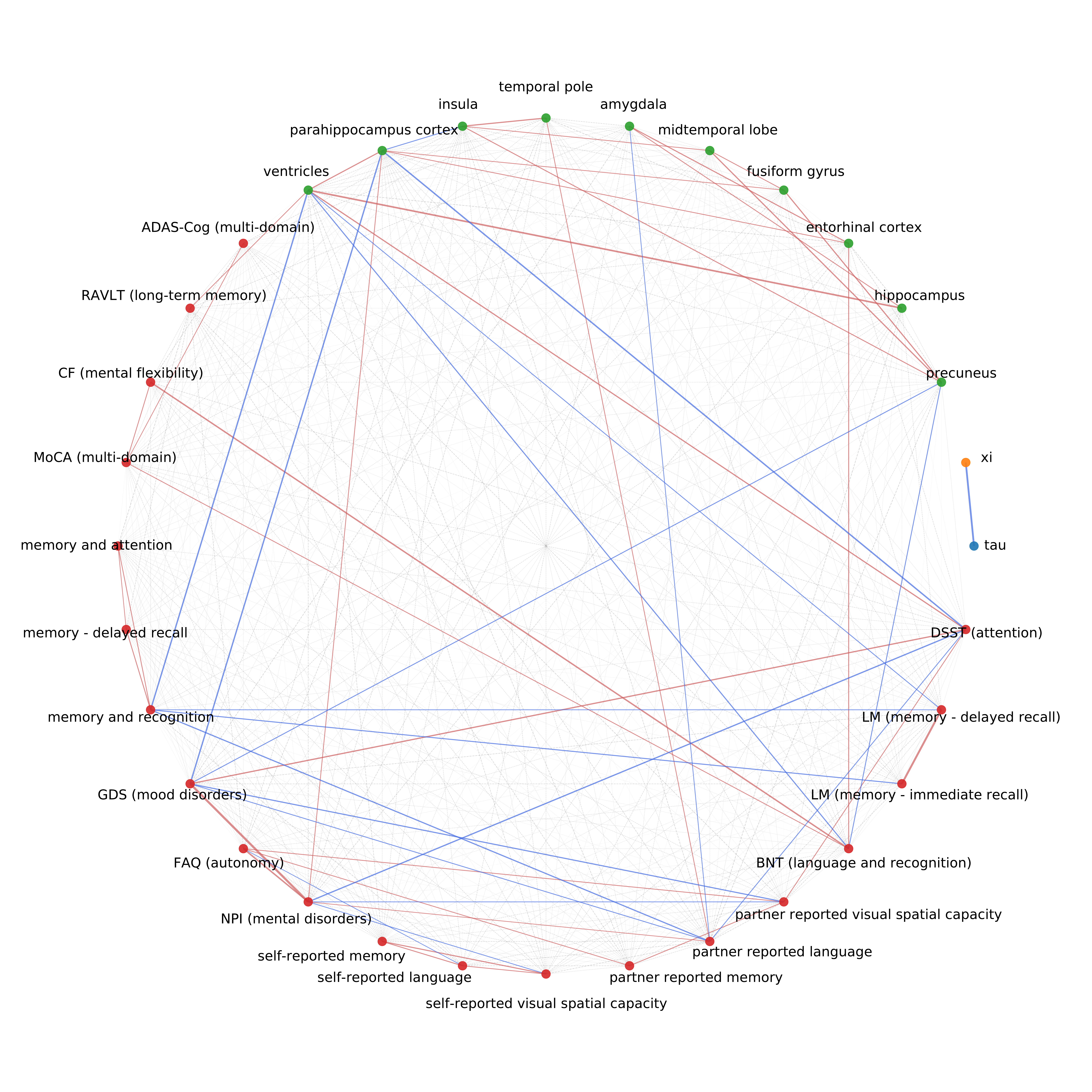}}
    \subfloat{\includegraphics[width=0.49\linewidth]{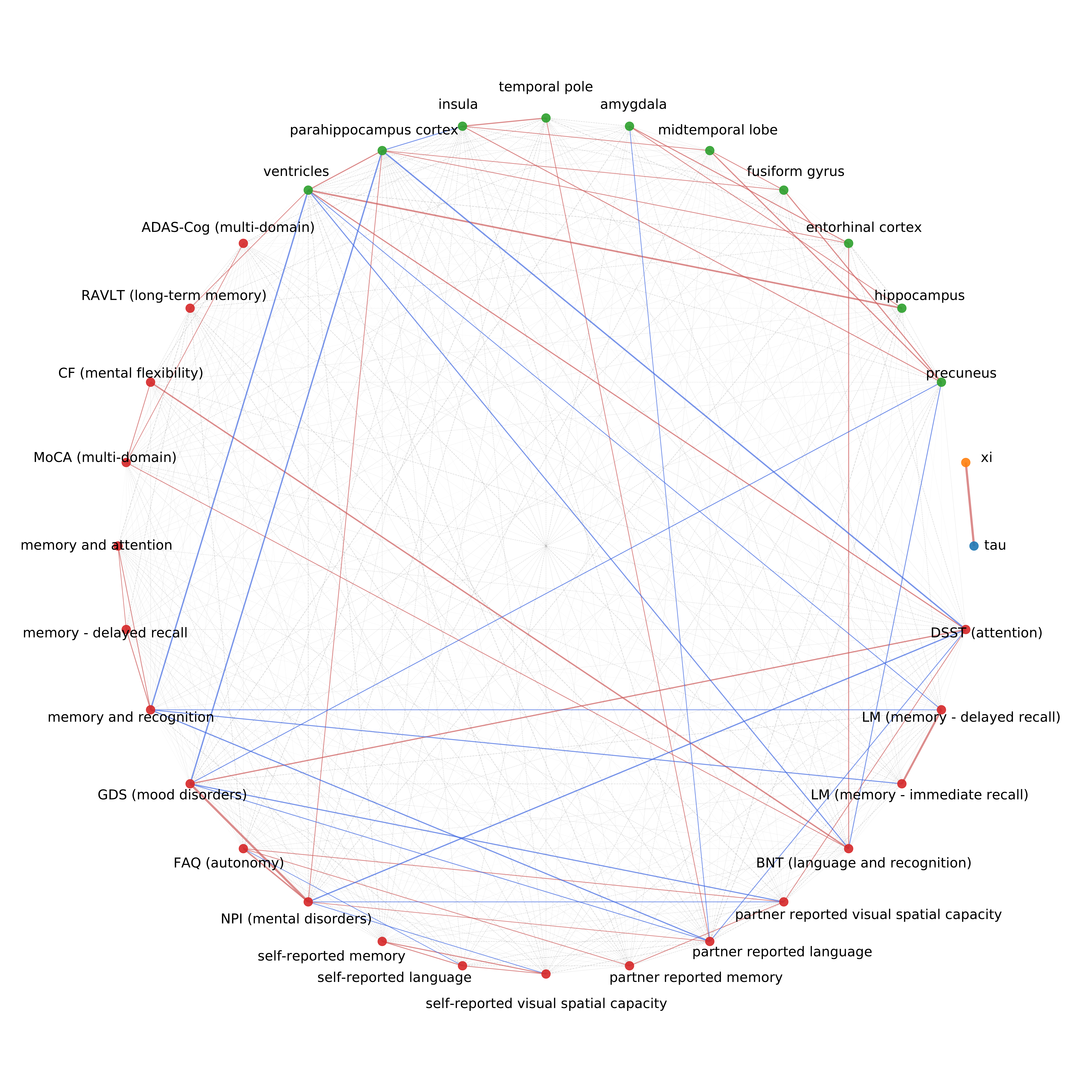}}
\caption{(Left) Conditional correlation graph estimated in a supervised way on the \textbf{"Control" patients}. (Right) Conditional correlation graph estimated in a supervised way on the \textbf{"AD" patients}. }
\label{fig:real_data_graph_oracle}
\end{figure*}

In this section, we confirm our experimental observations with a real, high dimensional, Alzheimer's Disease dataset. We illustrate that the EM with Mixture of CGGM is better suited to identify clusters correlated with the diagnostic than the Mixture of GGM methods. We bring to light the effect of co-features such as the gender and age on the medical features.\\ 
Our dataset is composed of 708 Alzheimer's Disease patients ("AD" patients) and 636 healthy patients ("Control" patients). Each patient makes several visits to the clinic (between 6 and 11 usually) at different moments of their lives. During each visit, the volumes of 10 predetermined cortical regions are measured with MRI, and tests are conducted to evaluate 20 different cognitive faculties (memory, language...). Instead of using this complex raw data, a longitudinal model is then estimated using the methodology of \cite{schiratti2015mixed}. Such a model describes the evolution (or non-evolution) over time of the measured variables with a geodesic trajectory that lives within a Riemannian manifold. Using all patients, a "baseline" trajectory is first estimated. Then, each individual patient's trajectory is described by their deviation from this baseline. In our example, this deviation is fully encoded within $32$ parameters (1 per observed variable + 2 additional parameters). The first model parameter is $\xi$, the time acceleration, which describes if, overall, the disease progresses faster or slower for the patient than for the baseline. The second model parameter is $\tau$, the time shift, which describes whether the disease starts sooner or later for the patient than for the baseline. Finally, there is one parameter $w_i$ for each of the $30$ measured variables ($i=1, ..., 10$ for the MRI regions and $i=11, ..., 30$ for the cognitive scores). Each of these describes whether the variable in question degenerates faster or slower than the baseline, after having taken into account the patient's overall disease acceleration $\xi$.\\ 
For our experiments, these $p=32$ parameters are the features in-between which we will draw graphical models: $Y := (\xi, \tau, (w_i)_{i= 1,..., 30})$. To illustrate the purpose of our approach: within this framework, a positive edge between two features $w_i$ and $w_j$ indicates that the corresponding features tend to co-degenerate together within the population. This is a crucial information for the understanding of the disease.\\
We include in the data three co-features that describe each patient and are often relevant in Alzheimer's Disease studies: the gender, the age at the first visit ("age baseline"), and the number of years of education. With the addition of the constant co-feature $=1$, the vector of co-features is 4-dimensional, $X \in \R^4$. The $n=1344$ (708 AD + 636 Control) patients constitute the unlabelled heterogeneous population. The "AD" patients were diagnosed with the Alzheimer's Disease, either from the start or after a few visits. Before running any method, The data $(Y,X)$ is centred and normalised over the entire population.\\
We run unsupervised methods with $K=2$ on this dataset in order to separate the data into 2 sub-populations and estimate a graphical model for each of them. As before, the three compared algorithms are EM, EM residual and C-EM. There is only one dataset available, hence, in order to check the stability of the results over several different runs, we implement a bootstrap procedure that only uses $70\%$ of the data each time. We generate 10 such bootstrapped dataset. We initialise the algorithms with a KMeans on the $Y^{(i)}$ data points. Since KMeans is not deterministic, we try 3 different KMeans initialisations for each bootstrapped dataset. Like previously, for the sake of fairness, the EM and C-EM are always provided with the same initialisation, and the residual EM is initialised with a KMeans on the residual of $Y$ after subtracting the prediction by the $X$, a more relevant initialisation for this method. We make all these runs with four different feature sets. First with no space shift variable $Y = \brace{\xi, \tau}, p=2$, then we add only the MRI space shifts $Y = \brace{\xi, \tau, (w_i)_{i= 1,..., 10}}, p=12$, then only the Cognitive Scores space shifts $Y = \brace{\xi, \tau, (w_i)_{i= 1,..., 10}}, p=22$, and finally, with all the features $Y = \brace{\xi, \tau, (w_i)_{i= 1,..., 30}}, p=22$. This is a total of $10 \textit{ bootstrapped datasets } \times 3 \textit{ initialisations } \times 4 \textit{ feature sets } = 120$ runs for each method. The classification results are summarised in Table \ref{tab:longiudinal_classification}. With two balanced classes, the classification error of a uniform random classifier is $50\%$. On the smallest dataset, $p=2$, we can see that the discovered cluster are somewhat correlated with the diagnostic, with classification errors below $30\%$. The Mixture of GGM on the uniformly residualised data and the Mixture of CGGM achieve similar levels of error, they are both better than the regular Mixture of GGM. When the MRI features are added, all the discovered cluster become more correlated with the diagnostic. The regular Mixture of GGM achieves in average $16\%$ of hard classification error, the residualised Mixture of GGM is at $11\%$ of error, and the Mixture of CGGM even below, at $7\%$. The results with only the Cognitive Scores are very similar, simply a bit worse for every method. However, when both the MRI and Cognitive Scores feature are included, the performances of both GGM mixtures decrease, with both higher average error and higher variance. On the other hand, the Mixture of CGGM achieves here its best level of performance. This stability of the Mixture of CGGM's performance as the size of the feature set increases indicates that our model is the best suited to properly identify clusters correlated with the diagnostic in high dimension.\\
We analyse the estimated Mixture of CGGM parameters on the full feature set $p=32$. First, since $\mathbb{E}[Y | X, z=k] = - \Lambda_k^{-1} \Theta_k^T X$ in the CGGM, we display on \figurename~\ref{fig:real_data_beta_CEM} the two estimated $\hat{\beta}_k :=  - \widehat{\Lambda}_k^{-1} \widehat{\Theta}_k$ (averaged over the bootstrap), as well as their difference. They play the role of linear regression coefficients in the model. The last column is the constant coefficient, while the first three are the gender, age baseline and years of education coefficients respectively. Since the data is centered, negative and positive values correspond to below average and above average values respectively. The cluster $k=1$ is the one very correlated with the Control patients sub-population. Similarly, the cluster $k=2$ is the one very correlated with the AD patients.\\
The most noticeable difference between the two $\hat{\beta}_k$ are the constant vectors, who have opposite effects on all features. In particular, the "AD cluster" is very correlated with high $\xi$, low $\tau$, an earlier atrophy of the ventricles, as well as high space shift for the two logical memory tests (immediate and delayed). The exact opposite being true for the "Control cluster". These are the expected effects: a high $\xi$ corresponds to a quickly progressing disease, and a low $\tau$ to an early starting disease.\\
The non-constant linear regression coefficients are also different between the clusters, although these differences are often in intensity and not in sign. In order to visualise more clearly the differences in intensity, we represent on the leftmost sub-figure of \figurename~\ref{fig:real_data_beta_CEM}, with the same conventions, the difference $\hat{\beta}_2-\hat{\beta}_1$. In particular, within the AD cluster, we observe stronger positive effect of the Age at the first visit on the space shifts corresponding to the Amygdala, entorhinal cortex, hippocampus and parahippocampus cortex atrophies. On the contrary, there is a stronger positive effect of the education level on all the space shifts of MRI atrophies for the Control patients. The age at the first visit has a stronger negative impact on $w_i$ corresponding to the scores self reported memory, language and visual spatial capacity for the AD patients, and a stronger negative impact on the two logical memory scores for the control patients. These differences, although less intense than the differences between the two constant terms, consequently impact the clustering, since without these three columns, the Mixture of CGGM would be equivalent to the less performing Mixture of GGM.\\
With the same conventions and scale, we display on \figurename~\ref{fig:real_data_beta_oracle}, the linear regression coefficients estimated by maximum likelihood on the Control and AD patients, with the oracle knowledge of the diagnostic this time. \figurename~\ref{fig:real_data_beta_CEM} and \ref{fig:real_data_beta_oracle} are very similar, in particular when it comes to their more potent coefficients. This shows that the C-EM, in addition to identifying clusters very correlated with the hidden diagnosis, also managed to recover the correct linear tendencies between features and co-features.\\
Finally, we display on 
\figurename~\ref{fig:real_data_graph}  
the average conditional correlation graphs estimated for the two clusters estimated by the Mixture of CGGM. Their only noticeable difference is the negative conditional correlation between $\xi$ and $\tau$ in the "Control cluster", which is reversed in the "AD cluster". For the AD patients, this means that a disease that appears later tends to also progress faster, which is in line with medical observations. Apart from this edge, the rest of the connections are almost identical in-between clusters. This suggests that the, cluster dependent, prediction $\mathbb{E}_{\hat{\theta}_k}[Y^{(i)} | X^{(i)}, z=k] =  - \widehat{\Sigma}_k \widehat{\Theta}_k^T X^{(i)}$ takes into account enough of the cluster-specific effects so that the remaining unexplained variance has almost the same form in both clusters. Hence, the conditional correlations pictured in these graphs correspond to very general effects, such as the positive correlations between related cognitive tests or areas of the cortex.\\ 
More strikingly, there are no conditional correlation between $\xi$ or $\tau$ and any of the space shifts $w_i$. This as consequent medical implications, since it suggests that the earliness ($\tau$) and speed ($\xi$) of the disease are conditionally independent from the succession of degradation that the patient's imagery and cognitive scores undergo. In other words, these graphs support the idea that the disease is the same regardless of whether it appears early/late and progresses fast/slowly.\\
As previously, we also estimated these conditional correlations graphs with the oracle knowledge of the diagnosis, they are displayed on 
\figurename~\ref{fig:real_data_graph_oracle}. 
As with the transition matrices, the oracle covariance graphs exhibit the same links as the ones recovered by the C-EM. This once again highlights that the Mixture of CGGM is an appropriate model to handle both cluster identification and structure recovery.\\
For the sake of comparison, the model parameters (both $\hat{\beta}$ and the graph) estimated with a single class CGGM are displayed in Appendix.
\begin{table*}
    \centering
    {\footnotesize
    \caption{\footnotesize Recovery of the diagnostic labels (AD or control) with unsupervised methods on real longitudinal data. The three compared methods are the EM, EM residual (both GGM) and the C-EM (CGGM). Four different feature sets are tried: only $\brace{\tau, \xi}$, adding the MRI space shift coefficients $w_i$, adding the Cognitive Score (CS) space shift coefficients $w_i$, and adding both the MRI and CS space shift coefficients. The table presents the average and standard deviation of the misclassification error over 10 bootstrap iteration, with 5 different KMeans initialisation each. The best results are in \textbf{bold}.}
    \label{tab:longiudinal_classification}
    \begin{tabular}{c|c|ccc}
        \toprule
                      & metric & EM          & EM resid.            & C-EM \\
        \midrule \\
        no CS, no MRI & soft misclassif.  & 0.31 (0.02) & 0.22 (0.03) & \textbf{0.21 (0.01)} \\
        $p=2$         & hard misclassif.  & 0.31 (0.03) & \textbf{0.18} (0.05) & 0.19 \textbf{(0.01)}\\
        \midrule \\
        only MRI      & soft misclassif.  & 0.15 (0.01) & 0.13 (0.01)          & \textbf{0.08 (0.01)}\\
        $p=12$        & hard misclassif.  & 0.12 (0.01) & 0.10 (0.01)          & \textbf{0.07 (0.01)}\\
        \midrule \\
        only CS       & soft misclassif.  & 0.17 (0.02) & 0.15 (0.02)          & \textbf{0.09 (0.01)}\\
        $p=22$        & hard misclassif.  & 0.14 (0.03) & 0.13 (0.03)          & \textbf{0.08 (0.01)}\\
       \midrule \\        
        CS and MRI    & soft misclassif.  & 0.24 (0.09) & 0.17 (0.04)          & \textbf{0.08 (0.01)}\\
        $p=32$        & hard misclassif.  & 0.21 (0.10) & 0.15 (0.05)          & \textbf{0.07 (0.01)}\\
        \bottomrule
    \end{tabular}
    }
\end{table*}

\section{Discussion and Conclusion}
\subsection{Discussion} 
In this work, we assumed that the number $K$ of classes was known and fixed. This corresponds for instance to the case where there is prior medical expertise that provides in advance the number of classes. However, as with all unsupervised methods, the question of the empirical estimation of $K$ is of great importance. Several approaches are possible. Since we have a Maximum Likelihood-type Estimator, a classical technique would be to use theoretical model selection criterion that penalises the number of degrees of freedom in the model. Famous examples include the Akaike information criterion (AIC), \cite{akaike1974new}, the Bayesian information criterion, \cite{schwarz1978estimating}, or the corrected AIC (AICc), \cite{hurvich1989regression}. There are $df(K, p, q) = K(\frac{p(p+1)}{2} + qp)$ degrees of freedom in a Mixture of CGGM with $K$ classes. After running an EM to completion for a value of $K$ and getting the estimates $\hat{\theta}^{K}, \hat{\pi}^{K}$, the loss function associated with class number $K$ is:
\begin{equation} \label{eq:CGGM_aic_bic}
\begin{split}
    - \sum_{i=1}^n ln\parent{\sum_{k=1}^K  \hat{\pi}^{K}_k \,   p_{\hat{\theta}^{K}_k}\parent{Y^{(i)} | X^{(i)}}} &+ n\, pen(\hat{\theta}^{K}, \hat{\pi}^{K})\\
    &+ \text{crit}(df(K, p, q), n)\, .
\end{split}
\end{equation}
Where $\text{crit}(df, n) = df$ for AIC, $\text{crit}(df, n) = df \frac{ln\, n}{2}$ for BIC and $\text{crit}(df, n) = df + \frac{df(df+1)}{n-df-1}$ for AICc. A preferred $\hat{K}$ is then chosen as the one to minimise \eqref{eq:CGGM_aic_bic} among a predefined grid of different $K$. This is notably the approach adopted by \cite{kim2016sparse}. Since this model selection technique requires one full run of the EM algorithm for each value of $K$, it can be computationally demanding. Moreover, the empirical validity of such theoretical criteria is not always guaranteed for all problems, see for instance \cite{lee2009performance} for an such an empirical analysis.\\
As a result, for EM algorithms like ours, one could consider a more data-driven selection procedure such as the \textit{Robust EM}, see \cite{figueiredo2002unsupervised, yang2012robust}. \textit{Robust EM} starts with a high number of classes that decreases and converges along the optimisation thanks to an entropy penalty term. However, this method makes each run of the EM more computationally demanding and can lead to very long execution times.\\
The relative model selection performances and computational burdens of these two different procedure will be studied in future works.\\
\\
One can also question the heterogeneity of the effect of the co-features, and wonder how the model fares in the presence of homogeneous co-feature effects. Indeed, in our empirical study we compare three models : (i) the Mixture of GGM which does not take into account the co-features; (ii) the Mixture of GGM with prior residualisation, which assume an homogeneous effect of the co-features across all hidden sub-populations; and (iii) the Mixture of CGGM, which allows the co-feature to have heterogeneous effects on the hidden sub-population. However, one can wonder how is the Mixture of CGGM able to describe, and detect, the situation where the/several co-features have a homogeneous effect on the/some features? Within the Mixture of CGGM, the heterogeneity assumption is translated as extra degrees of freedom. The homogeneous case is included in this model: it corresponds to the scenario where coefficients have equal values across classes. When it comes to the estimation procedure, one way to enforce truly equal values is to replace the Group Graphical Lasso (GGL) penalty used in our examples by the Fused Graphical Lasso (FGL) penalty. Like GGL, FGL is one of the penalties proposed for the supervised Hierarchical GGM estimation in \cite{danaher2014joint}. As such, it is one of the many penalties with which our EM is tractable. This penalty - thanks to terms of the form: $\sum_{l<k} \det{\Theta_l^{(ij)} - \Theta_k^{(ij)}}$ - is designed to equalise the coefficients whose estimated values across different classes are close. This allows the model to detect homogeneous effects of the co-features.\\
Another, more hands on, possibility that does not require the use of fused-type penalties is to conduct statistical testing during the optimisation procedure. In order to detect homogeneity, such tests would look for significantly close coefficients. Then, if applicable, constrain them to be equal from this point on, hence removing this degree of freedom from the model. Thanks to this reduced complexity of the resulting model, such an step could help make the procedure more stable in the presence of co-features with homogeneous effect.\\
Those approaches will also be tested as part of future works in order to evaluate the performances of Mixtures of CGGM in the case of homogeneous co-feature effects.

\subsection{Conclusion}
We introduced the Mixture of Conditional Gaussian Graphical Models in order to guide the cluster discovery when estimating different Gaussian Graphical Models for an unlabelled heterogeneous population in the presence of co-features. We motivated its usage to deal with the potential in-homogeneous and class-dependent effect of the co-features on the observed data that would otherwise disrupt the clustering effort. To estimate our Mixture model, we proposed a penalised EM algorithm ("Conditional EM" or  "C-EM") compatible with a wide array of penalties. Moreover, we provided detailed algorithmic steps in the specific case of the popular Group Graphical LASSO penalty, and made the corresponding code publicly available. 
Then, we demonstrated the interest of the method with experiments on synthetic and real data. First, we showed on a toy example - with a 2-dimensional feature space and a 1-dimensional co-feature - that the regular Mixture of GGM methods were inadequate to deal with even the most simple in-homogeneous co-feature. We confirmed on a more complex simulation, in higher dimension, that Mixtures of CGGM could identify much better the clusters in the feature space, and recover the actual GGM structure of the data. Finally, we tested all the methods on a real data set, with longitudinal model parameters describing the evolution of several Alzheimer's Disease patients. We demonstrated that our method was the best at identifying the diagnostic with an unlabelled dataset. We unearthed some in-homogeneous effect of co-features on the longitudinal parameter and recovered the conditional correlation graphs by cluster. These graphs hint at a conditional independence between the earliness and speed of the disease and the order in which the many degradation appear. This hypothesis will be tested in future studies. 

%
%

\section*{Appendix: Single class CGGM on the real data}
\begin{figure*}[tbhp]
    \centering
    \subfloat{\includegraphics[width=0.35\linewidth]{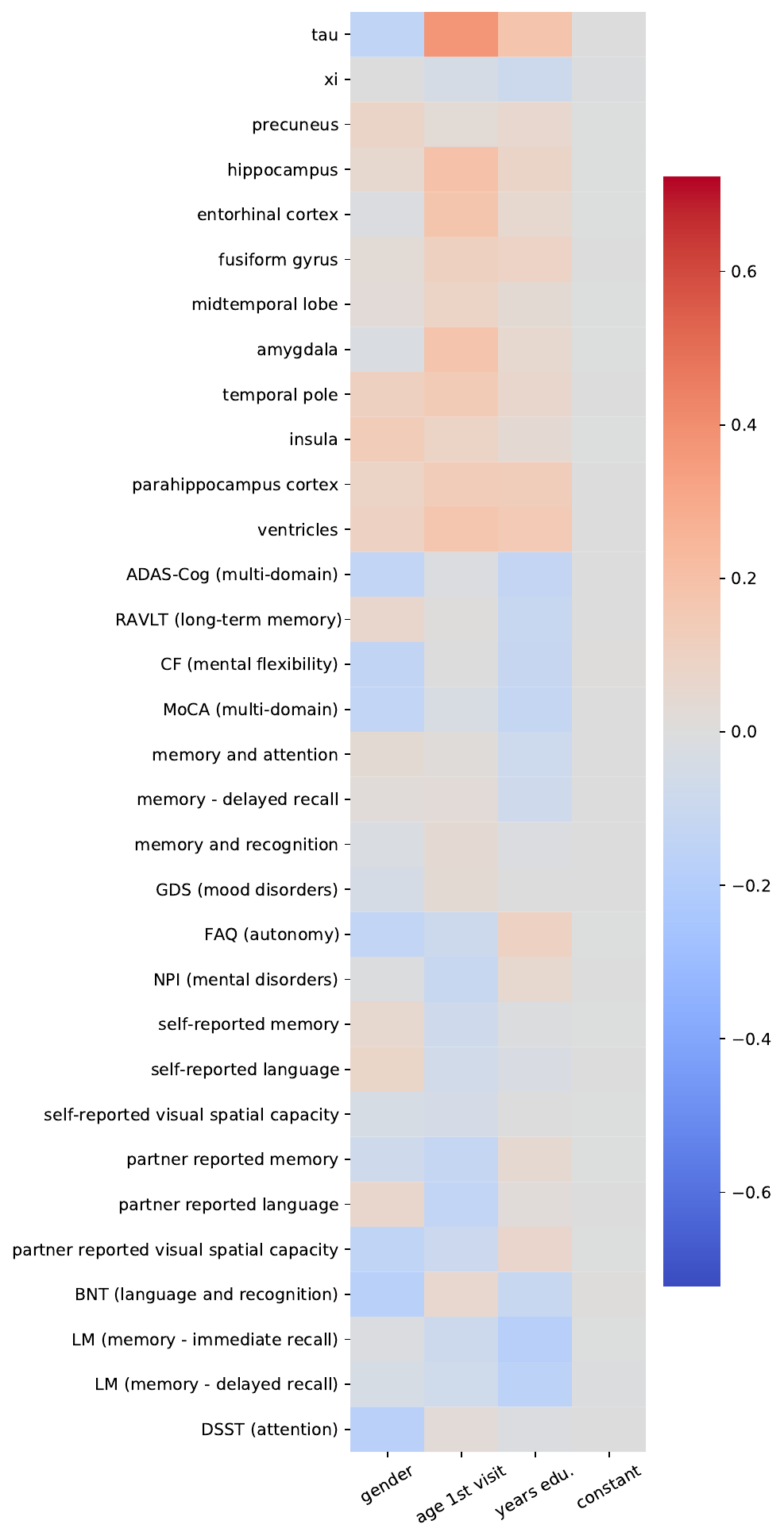}}
    \subfloat{\includegraphics[width=0.65\linewidth]{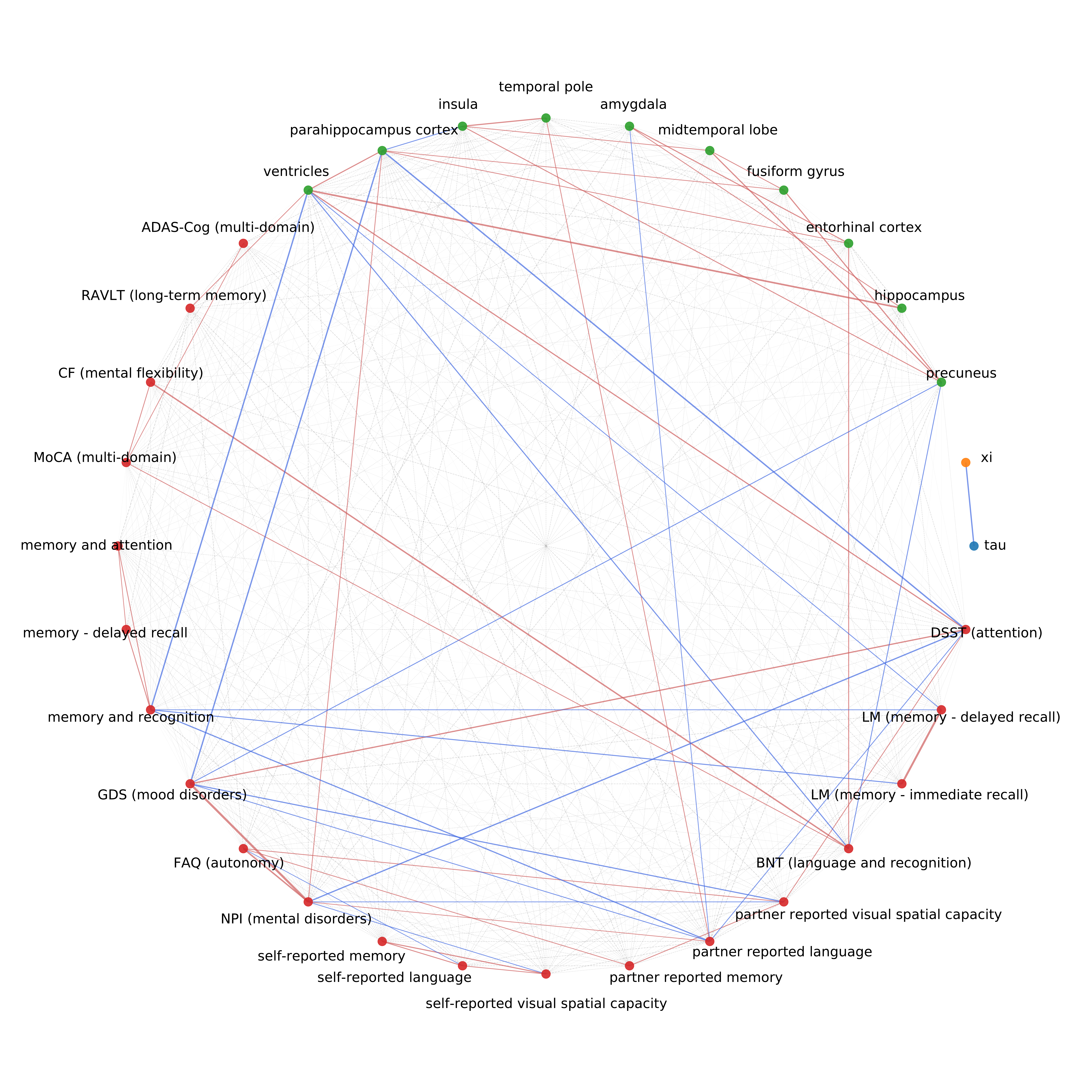}}\\
\caption{Parameters estimated with a simple CGGM on all data. (Left) $\hat{\beta} =  - \widehat{\Sigma} \widehat{\Theta}$. (Right) Conditional correlations graph.}
\label{fig:real_data_uni_class}
\end{figure*}

In this appendix, we take a look at the parameters (averaged over several bootstrap folds) estimated by fitting single CGGM on the real data. On \figurename~\ref{fig:real_data_uni_class}, we display both the estimated $\hat{\beta} =  - \widehat{\Sigma} \widehat{\Theta}$ between $X$ and $Y$ and the estimated conditional correlation graph in-between the components of $Y$. The constant term in $\hat{\beta}$ is 0 since the data is overall centred. Other than that, the coefficient intensities appear to be weaker than in the multi-class parameters. The conditional correlation graph on the other hand displays the negative correlation between disease earliness $\tau$ and speed $\xi$ that was characteristic of the Control patients on \figurename~\ref{fig:real_data_graph} and \ref{fig:real_data_graph_oracle}. This is despite the Controls ($n=636$) being slightly less numerous than the AD ($n=708$) patients on this dataset.

\section*{Declarations}
{
\small
\textbf{Funding.} The research leading to these results has received funding from the European Research Council (ERC) under grant agreement No 678304, European Union’s Horizon 2020 research and innovation program under grant agreement No 666992 (EuroPOND) and No 826421 (TVB-Cloud), and the French government under management of Agence Nationale de la Recherche as part of the "Investissements d'avenir" program, reference ANR-19-P3IA-0001 (PRAIRIE 3IA Institute) and reference ANR-10-IAIHU-06 (IHU-A-ICM).\\
\\
\textbf{Conflicts of interest.} On behalf of all authors, the corresponding author states that there is no conflict of interest.\\
\\
\textbf{Code availability.} Code for our algorithm, as well as a toy example that reproduces some of the results of this paper, publicly available at: \href{https://github.com/tlartigue/Mixture-of-Conditional-Gaussian-Graphical-Models}{https://github.com/tlartigue/Mixture-of-Conditional-Gaussian-Graphical-Models}\\
\\
\textbf{Authors' contributions.} Conceptualisation, T.L. and S.A.; methodology, T.L. and S.A.; software, T.L.; data curation, S.D. and T.L.; validation, T.L. and S.A.; visualisation, T.L.; result analysis, S.D., S.A. and T.L.; writing---original draft preparation; T.L.; writing---review and editing, T.L., S.D. and S.A.; supervision, S.D and S.A.. All authors have read and agreed to the published version of the manuscript. 
}

\bibliographystyle{spbasic}      
\bibliography{references_full}   

\end{document}